\documentclass[11pt]{article}

\usepackage[marginparwidth=35pt,textheight=23cm,footskip=1.3cm,voffset=0.2cm]{geometry}

\usepackage[toc,page]{appendix}

\usepackage{amsmath, amssymb, amsthm}

\usepackage{graphicx}

\newcommand{\longdowndownarrows}{
  \mathrel{\vcenter{\hbox{\scalebox{1}[2]{$\downdownarrows$}}}}
}

\usepackage{hyperref}

\usepackage{xcolor}

\usepackage{mathrsfs}

\usepackage[shortlabels]{enumitem}

\usepackage[T2A]{fontenc}
\usepackage[utf8]{inputenc}

\usepackage[activate={true,nocompatibility}
final,tracking=true,kerning=true,factor=1100,stretch=10,
shrink=10]{microtype}

\usepackage{hyperref}

\newcommand{\Adm}[3][{}]{\operatorname{Adm}^{#1}_{#2,#3}}
\newcommand{\RAdm}[3][{}]{\operatorname{RAdm}^{#1}_{#2,#3}}
\newcommand{\action}[2]{\left\langle #1,#2\right\rangle}

\newcommand{\Z}{{\ensuremath{\mathbb{Z}}}}

\newcommand{\R}{{\ensuremath{\mathbb{R}}}}
\newcommand{\ptd}{\mathcal{P}^1(\td)}

\newcommand{\pj}{\operatorname{p}}

\DeclareMathOperator{\Lip}{Lip}

\newcommand{\eps}{\varepsilon}

\DeclareMathOperator*{\spt}{supp}

\newcommand{\dtorus}{{\ensuremath{\mathbb{T}^d}}}
\newcommand{\onetorus}{{\ensuremath{\mathbb{T}^1}}}
\newcommand{\td}{{\ensuremath{\mathbb{T}^d}}}
\newcommand{\rd}{{\ensuremath{\mathbb{R}^d}}}

\newcommand{\wassonespace}{{\ensuremath{\mathcal{P}^1(\mathbb{T}^d)}}}
\newcommand{\EE}{{\ensuremath{\operatorname{\mathbb{E}}}}}

\newcommand{\tripleprod}{{\ensuremath{[0,T]\times\dtorus\times\mathbb{B}_C}}}

\usepackage[backend=biber,style=alphabetic]{biblatex}

\newtheorem{theorem}{Theorem}[section]

\newtheorem{definition}[theorem]{Definition}
\newtheorem{proposition}[theorem]{Proposition}
\newtheorem{lemma}[theorem]{Lemma}
\newtheorem{corollary}[theorem]{Corollary}

\newtheoremstyle{remarkstyle}
  {3pt}
  {3pt}
  {\normalfont}
  {}
  {\bfseries}
  {.}
  {.5em}
  {}

\theoremstyle{remarkstyle}
\newtheorem{remark}[theorem]{Remark}

\usepackage{etoolbox}
\AtEndEnvironment{remark}{\hspace{0.2em}$\triangle$}

 \title{Mean field game value multifunction: dynamic programming principle}
 \author{Yurii AVERBOUKH\footnote{Krasovskii Institute of Mathematics and Mechanics,
 and HSE University. E-mail: \href{mailto:ayv@imm.uran.ru}{ayv@imm.uran.ru}},
 Sergio MAYORGA\footnote{Innopolis University. E-mails: \href{mailto:s.mayorga@innopolis.ru}{s.mayorga@innopolis.ru},
 \href{mailto:me@mayorga.ru}{me@mayorga.ru}}}
 \date{\today}

\begin{document}

\maketitle

\begin{abstract}
We study mean field games (MFGs) on the torus with non-degenerate diffusion, in
a regime (no monotonicity) where multiple solutions exist for the same data.
In this setting the value of a representative player is not a function of
the initial time and distribution, and the classical master equation loses
its meaning. We describe the dependence of the value on the initial data by
a multivalued mapping, the value multifunction, assigning to each initial
time and distribution the set of value functions of all solutions
of the MFG starting from
that data. Our main result gives a sufficient condition for a closed multifunction to be
a value multifunction, namely, it must satisfy a dynamic programming principle formulated
on proper time subintervals, expressed as invariance of its graph under the forward
operators generated by equilibria of the game. This replaces the global requirement
that each value be produced by a solution of the forward-backward system on the whole
time interval by a local, time-consistent condition together with closedness, which
carries the terminal condition. The result is the second-order analogue of the
viability approach previously developed for first-order mean field games.
We reformulate solutions of the mean field game in terms of relaxed equilibria,
show that the closure of a value multifunction is again a value multifunction,
and prove that the forward operators compose as a semigroup. We also present
an explicit example in which a single terminal cost produces
at least two solutions for every initial time in an interval,
while solutions are unique for initial times sufficiently close to the horizon.
\end{abstract}

\tableofcontents

\paragraph{Acknowledgement} This article was partially prepared within the
framework of the HSE University Basic Research Program.

	\section{Introduction}

In this article we study a mean field game (MFG, from now on) on the torus with a non-degenerate
diffusion, in the case where the possibility for multiple solutions is open.
For such games the value of a representative player is not a function of the
initial time and the initial distribution of the players, because the same
initial data may admit several solutions with different values. We
describe this dependence by a multivalued mapping, called the value
multifunction, which assigns to an initial time and an initial
distribution a set of value functions of solutions. Our main result,
Theorem \ref{thm:aetoDPP}, gives a sufficient condition for a given
multivalued mapping to be a value multifunction, namely, that
it should be closed and satisfy a dynamic programming principle, in turn a
condition for proper time subintervals (i.e.~\([s,r]\) with $r<T$), which
is enough to verify only on intervals of a fixed small length (Corollary \ref{cor:localDPP}).

In this mean field game, the dynamics of a representative player is
given by the stochastic differential equation (SDE, from now on)
\begin{equation}\label{intro:sys:dynamics}
	dX_t=w(t,X_t)\, dt+\sqrt{2}\,dW_t, \qquad t_0\leq t\leq T,
\end{equation}
and his or her payoff is
\begin{equation}\label{intro:sys:payoff}
	\mathbb{E}\Big[G(X_T,m(T))+\int_{t_0}^T\big(F(X_t,m(t))-L(X_t,w(t,X_t))\big)dt
	\Big].
\end{equation}
We work on the torus $\td\triangleq \rd/\mathbb{Z}^d$, which serves as
the phase space. Here $m(t)$ or $m_t$ denotes the distribution of the
players at time $t$, a Borel probability measure on the phase space,
while feedback controls are Borel functions
$w:[t_0,T]\times\dtorus\to\R^d$, written $w_t(x)$ or $w(t,x)$; $L$ is a
Lagrangian, and $F$ and $G$ are the running and terminal interaction
terms.
The mean field game system of equations is the coupled
backward-forward system
\begin{equation}\label{intro:eq:mfg_system}
	\begin{split}
		\partial_t u+\triangle u+H(x,D u)+F(x,m) = & \ 0 \qquad \qquad
		\textrm{ in } \  (t_0,T)\times \dtorus,
		\\
		\partial_t m-\triangle m+\operatorname{div}\big(D_p H(x,D u)m\big)
		=
		& \ 0
		\qquad \qquad \textrm{ in } \ \mathcal{D}'((t_0,T)\times \dtorus),
		\\
		u(T,x) = & \ G(x,m(T)) \quad  \text{ on } \dtorus,
		\\
		m(t_0) = & \  m_*,
	\end{split}
\end{equation}
it consists of a backward Hamilton--Jacobi--Bellman (HJB, from now on,
or just Bellman) equation with terminal value $G(\cdot,m_T)$ and a
forward Fokker--Planck--Kolmogorov (FPK, from now on, or just
Kolmogorov) equation with initial value $m_*$; the latter is understood
in the sense of distributions. The Hamiltonian $H$ is defined for
$x\in\td$ and $p\in\R^d$ by
\begin{equation}\label{intro:def:Hamiltonian}
	H(x,p)\triangleq \max_{\alpha\in\rd} \big[p\cdot \alpha-L(x,\alpha)\big];
\end{equation}
$D_pH$ denotes its gradient with respect to $p$.

Mean field games were introduced independently by Lasry and Lions
\cite{LasryLions2006a,LasryLions2006b,lionsjapanese} and by Huang,
Malham\'e and Caines \cite{huang2006large,huangcainesmalhame},
the latter in the engineering community. They describe Nash equilibria
of differential games with a very large number of small, interchangeable
players, each of whom reacts to the distribution of the
whole ``conglomerate'' of players. The
theory was developed further in the lectures of Lions at the Coll\`ege
de France \cite{courscollegefrance}; for systematic accounts we refer to
\cite{bensoussanbook,CarmonaDelarue2018I,CarmonaDelarue2018II,cardaliaguetnotes,AchdouEtAl2020}.

Under the monotonicity in $m$ condition on $F$ and \(G\) of Lasry and Lions, the system
\eqref{intro:eq:mfg_system} has exactly one solution for every initial
time $t_0$ and initial distribution $m_*$. The value of the
representative player at time $t_0$ is then a function
$U(t_0,\cdot,m_*)$ of the initial data, and this function is time
consistent, by which we mean that if $(u,m)$ is the solution starting from $(t_0,m_*)$ and
$t_0<r<T$, then $u(t_0,\cdot)$ is obtained by solving the Bellman
equation on $[t_0,r]$ with terminal value $U(r,\cdot,m_r)$. This is a
dynamic programming principle for the value, and, formally, it leads to
an equation for $U$ on the space of measures, the master equation,
proposed by Lions \cite{courscollegefrance}. The master equation has been the
subject of intensive work. Classical solutions were constructed by
Chassagneux, Crisan and Delarue \cite{ChassagneuxCrisanDelarue2022} and,
for second order games with or without common noise, by Cardaliaguet,
Delarue, Lasry and Lions \cite{fourhorses}, who in fact also used it to prove
the convergence of Nash equilibria of $N$-player games to the solution
of the mean field game; see also
\cite{BuckdahnLiPengRainer2017,CarmonaDelarue2018II}. Short time results
and results for first order games were obtained by Gangbo and
\'Swi\k{e}ch \cite{mfgmain}, Mayorga \cite{mayorgajde1} and by Cardaliaguet, Cirant and
Porretta \cite{CardaliaguetCirantPorretta2023}. Global results under the
displacement monotonicity condition, an alternative to the condition of
Lasry and Lions, are due to Gangbo and M\'esz\'aros
\cite{GangboMeszaros2022} and to Gangbo, M\'esz\'aros, Mou and Zhang
\cite{GangboMeszarosMouZhang2022}; see also
\cite{AmbroseMeszaros2023,GraberMeszaros2023}. Mou and Zhang
\cite{MouZhang} treated data of low regularity, Bertucci
\cite{Bertucci2021,Bertucci2023} introduced a notion of monotone
solution, and Cecchin and Delarue \cite{CecchinDelarue2022} studied weak
solutions of the master equation of potential games beyond the regime
of uniqueness. The master equation has also been used to describe the
fluctuations and large deviations of $N$-player games
\cite{DelarueLackerRamanan2019,DelarueLackerRamanan2020}.

Without monotonicity the system \eqref{intro:eq:mfg_system} may have
several solutions; see, for instance,
\cite{bardi2019nonuniqueness,cirant2019timedependent,BrianiCardaliaguet2018,CecchinDaiPraFischerPelino2019,DelarueFoguenTchuendom2020}.
Then the value of the representative player is no longer a function of
$(t_0,m_*)$, and the master equation, in its classical form, loses its
meaning. The value multifunction replaces the function $U$ in this case
(Definition \ref{sol_mdf:def:value_multifunction}). It assigns to an
initial time $s$ and an initial distribution $\mu$ a set of functions on
the phase space, each of which is the value function, at time $s$, of
some solution of the MFG on $[s,T]$ starting from $\mu$.
Among value multifunctions there is a largest one, $\mathcal V^{\max}$:
its value at $(s,\mu)$, $s<T$, is the set of the functions $u(s,\cdot)$
over all solutions $(u,m)$ of the MFG on $[s,T]$ starting from $\mu$, and \(\mathcal V^{\max}(T,\mu)=\{G(\cdot,\mu)\}\).
It is a value multifunction, and every value multifunction $\mathcal W$ satisfies
$\mathcal W(s,\mu)\subset\mathcal V^{\max}(s,\mu)$ for all \((s,\mu)\).
When the solution is unique for all initial data, $\mathcal V^{\max}$ is
single-valued and its value is \(U(s,\cdot,\mu)\). The values are
functions of the state variable, not numbers, since an equilibrium
is shared by all players, as it depends on $(s,\mu)$ but not on the
initial position of an individual player.

The definition of a value multifunction is global in time. This means that
each of its values has to be produced by a solution of the MFG on the whole interval
$[s,T]$, that is, by a solution of a forward-backward system in which
the two equations are coupled through the terminal condition at time $T$. Theorem \ref{thm:aetoDPP}
replaces this global requirement by two conditions of different types. The first is the dynamic programming
principle of Definition \ref{defn:DPP}. It involves only intervals $[s,r]$ with $r<T$. On such an interval the terminal
datum is not $G$, it is a value provided by the multifunction itself at time $r$, and the principle asks that the multifunction be consistent with itself under
the forward operators $\Phi^{r,s}$ of Definition \ref{mfg_dyn:def:dynamics}, which ``move'' an initial pair value-measure
along an equilibrium of the game on $[s,r]$. In other words, the graph of the multifunction must be invariant under the mean field game
dynamics. The second condition is closedness. The terminal cost $G$ plays a role only in
the requirement that $\mathcal W(T,\mu)=\{G(\cdot,\mu)\}$. Since the dynamic programming principle involves only times $r<T$,
it is the closedness of $\mathcal W$ that forces the values obtained from it
to satisfy the terminal condition in the limit $r\to T$.
Moreover, because the forward operators compose as a semigroup, it is
enough to verify the principle on intervals $[s,r]$ of length smaller
than a fixed $h>0$ (Corollary \ref{cor:localDPP}). In this way the
principle is a condition on short time steps, which is the form needed
for an infinitesimal version.

To characterize a value function by the invariance of a closed set under
the dynamics, rather than by an equation, is an old idea in deterministic
control and differential games. Its basic case is Nagumo's theorem on
the invariant sets of an ordinary differential equation
\cite{Nagumo1942}; we refer to \cite{Aubin1991,JDE1995viability} for
viability theory and its applications to Hamilton--Jacobi equations, and
to \cite{CardaliaguetQuincampoixSaintPierre1999} for numerical schemes
based on it.

The present paper continues the program started by Averboukh in
\cite{Averboukh2020} for first order (deterministic) mean field games;
viability methods had earlier been applied to mean field type control
in \cite{Averboukh2018}. In \cite{Averboukh2020} the value multifunction
was introduced, viability with respect to the dynamics generated by the
game was shown to be a sufficient condition for a multifunction to be a
value multifunction, and an infinitesimal form of this condition was derived. Here we treat second order games, in which each player is
subject to an independent noise. The presence of the individual noise changes
the machinery, since the Bellman equation is parabolic and solved in the classical sense, the
Kolmogorov equation is understood in the sense of distributions, and
the equilibria on short intervals, which define the forward operators,
are those of a stochastic control problem. To obtain the compactness needed in the passage to the limit we
use relaxed controls. The approach of relaxed controls
that we have taken in our work
closely parallels that which is found
in \cite{ShapShat26}, but with a different purpose.

The
infinitesimal form of the dynamic programming principle in the second
order case is still the subject of an ongoing project.

Closest to the present paper is the work of Iseri and Zhang \cite{IseriZhang}, who were the first to study the dynamics of the set
of equilibrium values of a mean field game without uniqueness. For each
initial time \(t\) and initial distribution $\mu$ they consider the set of values over all equilibria, which they call the set value. They prove
that it satisfies a dynamic programming principle and that the set
values of the corresponding $N$-player games converge to it. They also show, in
discrete and continuous models, how much the set value depends on the class of admissible
controls (closed loop, state or path dependent, relaxed), and they point to a description by set-valued
partial differential equations \cite{IseriZhangSetPDE}. As in the
present paper, their set value consists of functions of the state
variable, for the reason explained above.

The problem solved in \cite{IseriZhang} is, nevertheless, different from ours in two respects. The first is the object. The set value of
\cite{IseriZhang} is defined through $\eps$-equilibria, with the limit \(\varepsilon\to0\)
taken in the definition, and two functions are
compared through integrals against $\mu$, so that only their values on
the support of $\mu$ matter. Our value multifunction is built from exact
equilibria, namely classical solutions of \eqref{intro:eq:mfg_system} or,
equivalently, relaxed equilibria, and its elements are value functions
defined on the whole torus. The equilibrium condition itself is the same: optimality on average with respect to the initial distribution
(Lemma \ref{lemma:characterizationofrelaxedequilibria}) is the case
$\eps=0$ of the condition of \cite{IseriZhang}. In their
terminology, \(\mathcal{V}^{\max}\) corresponds to the raw set
value, the set of values over exact equilibria. Iseri and Zhang workwith the
\(\eps\)-version because the raw set value is, in general,
sensitive to perturbations of the data; in our setting the set of exact
equilibrium values is closed (Proposition \ref{sol_mfg:prop:closure}),
and this is what allows us to work with exact equilibria. The second
respect, which is the more important one, is the direction of the main
result. In \cite{IseriZhang} the dynamic programming principle is a
property of one given object, the set value, and it is proved that this
object possesses this property. Here the dynamic programming principle is a condition
imposed on an arbitrary closed multifunction, and Theorem \ref{thm:aetoDPP} shows that the condition is sufficient for the
multifunction to consist of equilibrium values. The first kind of
result says that the set of equilibrium values is time consistent; the
second says that time consistency, together with closedness and the
terminal condition, identifies equilibrium values.
The two results are
complementary. The settings also differ: \cite{IseriZhang} works in
\(\R^d\) with bounded Lipschitz coefficients and controls taking values in
a fixed set, and addresses the convergence of $N$-player games, which we
do not consider; we work on the torus, with a convex and coercive
Lagrangian and classical solutions of the MFG system.

Besides Theorem \ref{thm:aetoDPP}, our work includes the following results. We reformulate the notion of solution of the MFG in terms of
relaxed equilibria (Proposition \ref{prop:solofMFGusingrelaxedequilibrium}), which is the form used in
the proof of the theorem. We show that the closure of a value multifunction is a
value multifunction (Proposition \ref{sol_mfg:prop:closure}); consequently
$\mathcal V^{\max}$ is closed, so that the closedness required in
Theorem \ref{thm:aetoDPP} is a property of the main object of interest
(Remark \ref{rem:closednessnatural}). We prove that the forward operators
compose as a semigroup, $\Phi^{\theta,s}=\Phi^{\theta,r}\circ\Phi^{r,s}$
(Proposition \ref{mfg_dyn:prop:superposition},
which yields Corollary \ref{cor:localDPP}).

Finally, in Section
\ref{section:nonuniquenessexample} we give an explicit example that
satisfies all our hypotheses \textbf{(HS)} below. In
this example all solutions are found explicitly, for an arbitrary
initial distribution, as fixed points of a scalar equation (Proposition \ref{prop:allsolutions}). A single terminal cost produces at least two
solutions for every initial time in an interval $[0,T-\varepsilon)$ (Corollary \ref{cor:allt0}), while for initial times close enough to $T$
the solution is unique for every initial distribution (Proposition
\ref{prop:uniquenessnearT}). The Lipschitz constant of the terminal cost
needed for the first statement grows like $1/\varepsilon$, and this order cannot be improved within the class
of terminal costs considered
(Corollary \ref{cor:allt0}).

Explicit examples of non-uniqueness for finite horizon mean field games
with a non-degenerate diffusion are not numerous; the best known are
those of Bardi and Fischer \cite{bardi2019nonuniqueness}. Our example
differs from theirs.
Namely, in \cite{bardi2019nonuniqueness} the state space is the real line, the
costs depend on the distribution of the players through its mean, and the
Hamiltonian is convex and either not differentiable at the origin, in
which case two solutions exist for every time horizon, or smooth near the origin and linear away from it, in which case two solutions exist
for horizons longer than a certain threshold. In both cases the
Hamiltonian comes from controls with values in a bounded interval, and
the corresponding Lagrangian does not satisfy our hypotheses \ref{prel:cond:L_convex}
and \ref{prel:cond:L_coercitive}. Our example is posed on the torus, its
Hamiltonian $H(p)=\frac12|p|^2$ is smooth and strictly convex, and the dependence of the terminal cost on the distribution is of a different
kind. Also, the notions of solution are different. In
\cite{bardi2019nonuniqueness} a solution is a classical solution of the
system of partial differential equations, with $m(t,\cdot)$ a
probability density for every $t$ and a density as initial datum (a
probabilistic formulation is treated there as well). Here the Bellman
equation is required to hold in the classical sense, but the Kolmogorov
equation only in the sense of distributions, for a continuous path of probability measures starting from an arbitrary initial distribution
(Section \ref{subsection:notionofclassicalsolutionofmfgsystem}). For the
pairs of distinct solutions exhibited in Section
\ref{section:nonuniquenessexample} this makes no difference, since they
are smooth, but it makes a difference for the two statements that
concern all solutions, Propositions \ref{prop:allsolutions} and
\ref{prop:uniquenessnearT}, because they hold for every initial distribution in
\(\mathcal P^1(\mathbb T^1)\).

Bardi and Fischer also prove that, when the Hamiltonian is smooth and the time horizon is short, the MFG system has at most one classical
solution, without any monotonicity of the costs \cite[Theorem~4.7]{bardi2019nonuniqueness}; the result goes back to a
lecture of Lions, and its proof applies to periodic data as well
\cite[Remark~4.8]{bardi2019nonuniqueness}. After the change of sign
$u\mapsto-u$, which turns our maximization into their minimization, our
example satisfies the assumptions of that theorem for initial
distributions with a bounded density.
Uniqueness for initial times close to \(T\) was,
therefore, something to be expected. Proposition \ref{prop:uniquenessnearT}
establishes it in the example with an explicit threshold on $T-t_0$;
unlike the one in \cite{bardi2019nonuniqueness} this threshold does not depend on the
supremum norm of the initial density and applies to every initial distribution. Bardi and Fischer also observe
that uniqueness holds for every horizon when the Lipschitz constant of the
gradient of the terminal cost with respect to the distribution is small
\cite[Remark~4.11]{bardi2019nonuniqueness}; this is consistent with the
growth of the Lipschitz constant in Corollary \ref{cor:allt0}.

The article is organized like this: Section 2 collects notation,
hypotheses and the notion of classical solution. Section 3 contains the
example. Section 4 introduces control processes, relaxed controls and
stochastic controls, and proves the properties of the value functions
used later. Section 5 introduces relaxed equilibria, the value
multifunction, and its closure. Section 6 introduces the forward
operators, the dynamic programming principle, and proves Theorem
\ref{thm:aetoDPP}. Standard properties of relaxed controls are proved in
the Appendix.

\paragraph{Artificial Intelligence (AI) use disclosure} We used DeepSeek to generate
an example of non-uniqueness, which we manually verified, and then
further worked out and generalized
(see Section \ref{section:nonuniquenessexample} and Remark \ref{rem:originalexample}
for details). We also used it to spot mistakes in the draft, missing
steps in proofs, and insufficiency in some of the justifications. DeepSeek's
feedback led in some degree to the reorganization of the draft,
but this task was done entirely and manually by us.

	\section{Preliminaries}
	\subsection{General notation}
	If $(\Omega,\mathcal{F})$, $(\Omega',\mathcal{F}')$ are measurable
	spaces, $\mathbb{P}$
is a probability on $\mathcal{F}$, $h:\Omega\rightarrow\Omega'$ is an
$\mathcal{F}/\mathcal{F}'$-measurable mapping, then $h\sharp\mathbb{P}$ is a
 probability on $\mathcal{F}'$ defined by the rule:
 for $E \in\mathcal{F}'$,
	\[
    (h\sharp\mathbb{P})(E) \triangleq \mathbb{P}(h^{-1}(E));
    \]
	here and throughout the text we use the symbol
    $\triangleq$ to denote definition.
	Further, if $(\Omega,\mathcal{F},\mathbb{P})$ is a probability space,
	$\xi:\Omega\rightarrow\mathbb{R}^d$ is a $d$-dimensional random variable,
	then $\operatorname{Law}(\xi)$ stands for
    the distribution of the random vector $\xi$
    with respect to $\mathbb{P}$, i.e.~\(\operatorname{Law}(\xi)\triangleq \xi\sharp\mathbb{P}\).
	If $n$ is a natural number, $X_1,\ldots,X_n$ are sets,
	$i_1,\ldots,i_k$ are indices from $\{1,\ldots,n\}$, then
	the projection operator $\operatorname{p}^{i_1,\ldots,i_k}:
	X_1\times\ldots\times X_n\rightarrow X_{i_1}\times\ldots\times
	X_{i_k}$ is defined by the rule:
	\[
    \operatorname{p}^{i_1,\ldots,i_k}(x_1,\ldots,x_n)\triangleq
	(x_{i_1},\ldots,x_{i_k}).
    \]
	If $(X,\rho_X)$ is a Polish space, then
	\begin{itemize}[-~~]
		\item $\mathbb{B}_r(x)$ denotes the closed ball of
		the radius $r$ centered at the point $x$; if, additionally,
		$X$ is a Banach space, $x$ is the origin, we omit the argument;
		\item $\mathcal{C}_b(X)$ stands for the set of all bounded continuous
		functions with values in $\mathbb{R}$;
		\item For \(K>0\), $\operatorname{Lip}_K(X)$ denotes
		the set of $K$-Lipschitz
		continuous functions from $X$ to $\mathbb{R}$;
		\item the set of all Borel measures on $X$ is denoted by
		$\mathcal{M}(X)$;
		\item the set of Borel probabilities on $X$ is denoted by
		$\mathcal{P}(X)$;
        \item if \(\phi\in \mathcal{C}_b(X)\) or some other class
        of functions on $X$, and \(m\in \mathcal{P}(X)\), then we will often denote
		\begin{equation}\label{eq:actionnotation}
		\action{m}{\phi}\triangleq \int_X \phi(x)m(dx).
		\end{equation}
		In this notation the domain of integration for the value
		\(\action{m}{\phi}\) will always be implicit in the space to which
		\(m\) belongs. We may also some times omit the domain of integration
        when this domain is made clear by the measure against which the integral
        is carried out.
	\end{itemize}

\subsection{Optimal transport notions}

A sequence of probability
measures $m_n$ converges narrowly to $m$ by definition
if, for each $\phi\in \mathcal{C}_b(X)$, one has
	\[
	\action{m_n}{\phi}\rightarrow \action{m}{\phi}.
	\]
	We will work with the space of probabilities
with finite first moment $\mathcal{P}^1(X)$; this consists of all
Borel probabilities $m$ such that, for some (equivalently, every)
$x_*\in X$
	\[\int_X\rho_X(x,x_*)m(dx)<\infty.\]  This space is metricized
	using the 1-Wasserstein (Kantorovich-Rubinstein) distance
	defined by the rule: for $m_1,m_2\in\mathcal{P}^1(X)$:
	\[
    W_1(m_1,m_2)\triangleq
	\inf_{\pi\in\Pi(m_1,m_2)}\int_{X\times X}\rho_X(x_1,x_2)\pi(dx_1dx_2).
    \]
	Here $\Pi(m_1,m_2)$ is the set of probabilities $\pi\in
	\mathcal{P}(X\times X)$ such that
	\[\operatorname{p}^1\sharp \pi=m_1,\ \
	\operatorname{p}^2\sharp \pi=m_2.\]
	Elements of $\Pi(m_1,m_2)$ are called plans between $m_1$ and
	$m_2$, while the set of
    plans between $m_1$ and $m_2$ that are optimal
    for the $W_1$ distance is denoted by $\Pi_0(m_1,m_2)$.
    Due to the Kantorovich-Rubinstein duality
	(see \cite[p.~88]{santambrogiobook}),
	\[
	W_1(m_1,m_2)=\sup\Bigg\{\int_X\psi(x)\big[m_1(dx)-m_2(dx)]:
	\psi\in\operatorname{Lip}_1(X)\Bigg\}.
	\]
\subsection{Data hypotheses \textbf{(HS)}}\label{subsection:datahyp}
In this subsection we fix the operating assumptions
that will be used in our work.

	We fix an arbitrary
	number \(T\in (0,\infty)\), the ``time horizon'', for the rest of the paper.

\subsubsection{Hypotheses on the Lagrangian}

The function $L:\dtorus\times \R^d\to \R$ is subjected
 to the following conditions:
\begin{enumerate}[label=\textbf{(L\arabic*)}]
	\item\label{prel:cond:L_lipschitz} The function
    $x\mapsto L(x,v)$ is Lipschitz
	continuous with a constant $K_L$ for any $v\in\R^d$.
	\item\label{prel:cond:L_convex}
$L$ is twice differentiable in $v$ at any \((x,v)\in \dtorus\times \R^d\),
with
\[
 D_{vv}L(x,v) >0.
 \]
In particular, $L$ is strictly convex in $v$.
	\item\label{prel:cond:L_coercitive}
    (Coercivity) There exist
	constants $q$, $C_L^0$,
	$C^1_L$ with $q>1$, $C_L^0>0$ such that
	\[
    L(x,v)\geq C_L^0|v|^{q}+C_L^1.
    \]
\end{enumerate}
We note for later use that \ref{prel:cond:L_lipschitz} and
\ref{prel:cond:L_convex} imply that $L$ is jointly continuous, and that
$M_L^0\triangleq\max_{x\in\td}|L(x,0)|$ is finite.

\subsubsection{Hypotheses on the interaction terms $F$ and $G$}

The functions $F$ (the running interaction)
and $G$ (the terminal interaction)
are subjected to fulfill the properties:
\begin{enumerate}[label=\textbf{(\Alph{enumi}-Lip)}]
	\setcounter{enumi}{5}
	\item\label{sol_mfg:cond:F} $|F(x_1,m_1)-F(x_2,m_2)|
	\leq K_F(|x_1-x_2|+W_1(m_1,m_2))$
	for some constant $K_F$,
	\item\label{sol_mfg:cond:G} $|G(x_1,m_1)-G(x_2,m_2)|
	\leq K_G(|x_1-x_2|+W_1(m_1,m_2))$
	for some constant $K_G$.
\end{enumerate}
These conditions and the compactness of
$\td\times\mathcal{P}^1(\td)$
imply that
\begin{equation}\label{sol_mfg:ineq:F_bound}
	\big|F(x,m)\big|\leq M_F,
\end{equation}
\begin{equation}\label{sol_mfg:ineq:G_bound}
	\big|G(x,m)\big|\leq M_G
\end{equation}
for all \(x\) and \(m\) and
some constants $M_F$ and $M_G$.
Additionally, we impose:
\begin{enumerate}[label=\textbf{(\Alph{enumi}-(2+$\alpha$))}]
	\setcounter{enumi}{6}
	\item\label{sol_mfg:cond:G2plusalpha}
    There is $\alpha\in (0,1)$ such that
    \(G(\cdot,m)\) is of class \(\mathcal C^{2+\alpha}\)
    uniformly in $m\in \mathcal{P}^1(\dtorus)$; that is,
    for any $m$
    the second derivative $D_{xx}G(x,m)$ exists,
    with the function $x\mapsto D_{xx}G(x,m)$ being
    H\"older with exponent $\alpha$ and with
    a H\"older norm that is independent of $m$.
\end{enumerate}
Since $\dtorus$ is compact,
the latter condition implies that the functions \(G(\cdot,m)\)
are of class $C^{1,1}$, in the sense that they are
Lipschitz, with Lipschitz derivative.

We will collectively refer to the hypotheses just listed
as \emph{the Hypotheses Set} \textbf{(HS)}.

\subsection{The Hamiltonian}\label{subsection:Hamiltonian}

The following elementary properties of the Hamiltonian
\eqref{intro:def:Hamiltonian} follow from
\ref{prel:cond:L_lipschitz}--\ref{prel:cond:L_coercitive} alone.
They will be used repeatedly; in particular, part (iii) makes
the term \(D_pH(x,Du)\) in the system \eqref{intro:eq:mfg_system}
well defined.

\begin{lemma}\label{prel:lm:Hamiltonian}
For every $(x,p)\in\dtorus\times\rd$:
\begin{enumerate}[label=(\roman*)]
\item the maximum in \eqref{intro:def:Hamiltonian} is attained at a
unique point $v^*(x,p)\in\rd$; consequently
\begin{equation}\label{eq:Young}
H(x,p)\geq p\cdot v-L(x,v)\quad\text{for all }v\in\rd,
\quad\text{with equality if and only if }v=v^*(x,p);
\end{equation}
\item for every $R>0$ there is a constant $C_*[R]$, depending only on
$R$, $C_L^0$, $C_L^1$, $q$ and $M_L^0$, and nondecreasing in $R$, such
that $|v^*(x,p)|\leq C_*[R]$ whenever $|p|\leq R$;
\item the map $(x,p)\mapsto v^*(x,p)$ is continuous, and $H(x,\cdot)$
is differentiable at every $p$, with
\[
D_pH(x,p)=v^*(x,p);
\]
\item $|H(x,p)-H(x',p)|\leq K_L|x-x'|$ for all $x,x'\in\td$, $p\in\rd$,
and $|H(x,p)-H(x,p')|\leq C_*[R]\,|p-p'|$ whenever $|p|,|p'|\leq R$.
In particular, $H$ is jointly continuous.
\end{enumerate}
\end{lemma}
\begin{proof}
(i) By \ref{prel:cond:L_convex} and \ref{prel:cond:L_coercitive}, the
function $v\mapsto p\cdot v-L(x,v)$ is continuous, strictly concave and
tends to $-\infty$ as $|v|\to\infty$; hence its supremum is attained at a
unique point $v^*(x,p)$. Statement \eqref{eq:Young}
follows.

(ii) Comparing the values at $v^*(x,p)$ and at $v=0$,
\[
|p|\,|v^*(x,p)|-C_L^0|v^*(x,p)|^{q}-C_L^1
\ \geq\ p\cdot v^*(x,p)-L\big(x,v^*(x,p)\big)
\ \geq\ -L(x,0)\ \geq\ -M_L^0 .
\]
The function $\varrho\mapsto
R\varrho-C_L^0\varrho^q-C_L^1$ tends to $-\infty$ as
$\varrho\to\infty$; taking $C_*[R]$ to be the largest $\varrho$ at
which its value is at least as large as $-M_L^0$ yields the claim for $|p|\leq R$.

(iii) Let $(x_n,p_n)\to (x,p)$. Using (ii) with $R=\sup_n|p_n|$, the sequence $v_n\triangleq v^*(x_n,p_n)$ is bounded.
If a subsequence $v_{n_k}\to \bar{v}$ as $k\to \infty$ then, passing to
the limit in the inequalities $p_{n_k}\cdot
v_{n_k}-L(x_{n_k},v_{n_k})\geq p_{n_k}\cdot v-L(x_{n_k},v)$,
$v\in\R^d$ we
obtain $p\cdot\bar{v}-L(x,\bar{v})\geq p\cdot v-L(x,v)$ for all $v$,
so, by uniqueness, $\bar{v}=v^*(x,p)$. Hence the whole sequence converges to $v^*(x,p)$. For the
differentiability, let $h\in\rd$; from the definition of
$H(x,p+h)$ and $H(x,p)$ applied with the maximizers $v^*(x,p)$ and
$v^*(x,p+h)$ respectively we get
\[
h\cdot v^*(x,p) \leq H(x,p+h)-H(x,p) \leq h\cdot v^*(x,p+h),
\]
whence $0\leq H(x,p+h)-H(x,p)-h\cdot
v^*(x,p)\leq|h|\,|v^*(x,p+h)-v^*(x,p)|=o(|h|)$ as $h\rightarrow0$, by
the continuity just proved.

(iv) For the first estimate, using \eqref{eq:Young} at $(x',p)$ with
$v=v^*(x,p)$ and \ref{prel:cond:L_lipschitz},
\[
H(x,p)-H(x',p)\leq
\big[p\cdot v^*(x,p)-L(x,v^*(x,p))\big]
-\big[p\cdot v^*(x,p)-L(x',v^*(x,p))\big]\leq K_L|x-x'|.
\]
For the second, if $|p|,|p'|\leq R$, then
\[
H(x,p)-H(x,p')\leq (p-p')\cdot v^*(x,p)\leq C_*[R]\,|p-p'|.
\]
The inequalities for the opposite differences are obtained
in the same way.
\end{proof}

\subsection{Notion of classical solution of the MFG system of
equations}\label{subsection:notionofclassicalsolutionofmfgsystem}

By a \textit{classical solution of the system
\eqref{intro:eq:mfg_system} with initial time $t_0$
on the time interval
$[t_0,T]$ and
initial distribution
$m_*$},
we will mean
a pair $(u,m_{\cdot})$, where
$u:[t_0,T]\times\dtorus\to\R$
and $m_{\cdot}:[t_0,T]\rightarrow\mathcal{P}^1(\dtorus)$,
such that:
\begin{itemize}[-~~]
\item $u\in\mathcal{C}^0([t_0,T]\times\dtorus)$; $u$ is continuously
differentiable in $t$ on $(t_0,T)\times\dtorus$; for every
$t\in[t_0,T]$ the function $u(t,\cdot)$ is twice differentiable, and
the derivatives $D u$, $D_{xx}u$ are continuous on the closed
cylinder $[t_0,T]\times\dtorus$;
\item $u$ satisfies the Bellman equation of \eqref{intro:eq:mfg_system}
pointwise on $(t_0,T)\times\dtorus$, and $u(T,\cdot)=G(\cdot,m_T)$;
\item $m_{\cdot}$
belongs to $\mathcal{C}^0([t_0,T];\mathcal{P}^1(\dtorus))$,
i.e.~it is a continuous path in $\mathcal{P}^1(\dtorus)$,
satisfies the FPK equation of \eqref{intro:eq:mfg_system} in the weak
(distributional) sense on $(t_0,T)\times\dtorus$
(the drift $D_pH(x,Du)$ being well defined by
Lemma \ref{prel:lm:Hamiltonian}), and $m_{t_0}=m_*$.
\end{itemize}

\begin{remark}\label{rem:dtextension}
If $(u,m_\cdot)$ is as above, then, by the Bellman equation,
\[
\partial_t u(t,x)=-\triangle u(t,x)-H(x,Du(t,x))-F(x,m_t)
\qquad\text{on }(t_0,T)\times\dtorus,
\]
and the right-hand side is continuous on the closed cylinder
$[t_0,T]\times\dtorus$ (by the regularity of $Du$, $D_{xx}u$, the joint
continuity of $H$ from Lemma \ref{prel:lm:Hamiltonian}(iv), the
continuity of $m_\cdot$ and \ref{sol_mfg:cond:F}). Hence $\partial_t u$
extends continuously to $[t_0,T]\times\dtorus$, so that
$u\in\mathcal{C}^{1,2}([t_0,T]\times\dtorus)$, and It\^o's formula may
be applied to $u(t,X_t)$ on any subinterval of $[t_0,T]$. The same
observation applies word for word to classical solutions of the
terminal-value problem \eqref{eq:HJBfsigma} considered in Lemma
\ref{prel:lm:verification} below.
\end{remark}

It is known (see, e.g.~\cite{lionsjapanese}) that
under the Hypotheses Set \textbf{(HS)}, the MFG
system \eqref{intro:eq:mfg_system} admits a classical
solution, with precisely the regularity required above.
However, as mentioned in the introduction,
we do not impose any of the conditions that are typically further
used in statements of uniqueness of this solution.
We will provide a specific example exhibiting non-uniqueness
in the next section.

\textit{Note}. The font style difference between
$C^{1,1}$, on one hand, and $\mathcal{C}^0$, $\mathcal{C}^{1,2}$ etc.~on the other
is not a typo. The class $C^{1,1}$ is the class of Lipschitz functions
with a Lipschitz derivative.

\section{A non-uniqueness example}\label{section:nonuniquenessexample}

The essence of the example presented here was discovered
by DeepSeek, when we prompted it to find a concrete example
of non-uniqueness that would satisfy our hypotheses.
We have somewhat generalized and systematized the exampled
and established some additional features, contained
in Lemma \ref{prop:allsolutions} and Propositions
\ref{prop:nonuniqueness} and \ref{prop:uniquenessnearT} below.

Remark \ref{rem:originalexample} below reports, within in our generalized context,
what the case was of the example found by DeepSeek.

Set
\begin{align*}
	F \equiv 0, \qquad G(x,\mu) = 2 \log\big(1+ \mathcal{G}(\mu)\cos(2\pi x)\big),
	\qquad L(x,v) \equiv \frac{1}{2}|v|^2,
\end{align*}
where $\mathcal{G}$ is a $W_1$-Lipschitz
function taking values in $[\theta_1,\theta_2]\subset (-1,1)$.
Such data satisfy our hypotheses \textbf{(HS)}. Consider the HJB
terminal condition problem
\begin{equation}\label{eq:HJBwithterminalG}
	\left.
	\begin{aligned}
		\partial_t u(t,x) + \frac{1}{2}|Du(t,x)|^2 + \Delta u(t,x) = & \ 0
		\qquad \text{ in } (t_0,T)\times \onetorus,
		\\
		u(T,x) = & \ G(x,\mu) \qquad \text { on } \onetorus
	\end{aligned}
	\right\}
\end{equation}
and let
\begin{align*}
p_s(x) \triangleq\frac{1}{\sqrt{4\pi s}}\sum_{j\in \Z}e^{-\frac{(x-j)^2}{4s}},
\quad s\in (0,T], \quad x\in\mathbb{R}^1, \quad p_0 \triangleq\delta_0.
\end{align*}
The solution of \eqref{eq:HJBwithterminalG} is unique and given by
\begin{align*}
u(t,x) = & \ 2\log\big(e^{\frac{G(\cdot,\mu)}{2}} \star p_{T-t}(x)\big)
=  2 \log\big(  1+ \mathcal{G}(\mu)e^{-4\pi^2(T-t)}\cos(2\pi x)\big)
\end{align*}
where $\star$ denotes the convolution on the torus in the $x$ variable.
Let us denote
\begin{align*}
	q^\mu(t,x) \triangleq 1+ \mathcal{G}(\mu)e^{-4\pi^2(T-t)}\cos(2\pi x),
\end{align*}
so that $u(t,x) = 2\log q^\mu(t,x)$ and $q^\mu$ satisfies
the linear backward heat equation
$$
\partial_t q^\mu(t,x) + \Delta q^\mu(t,x) = 0
$$
with $q^\mu(T,x) = e^{G(x,\mu)/2}$. (The solution to the non-linear one is
actually found by first making the change of variables $w=e^{u/2}$, which
gives the linear equation solved by $q^\mu$.)

We now cast the MFG system in this example into a fixed-point
form. Fix an
initial probability measure
$m_*\in\mathcal{P}^1(\onetorus)$.
Given $\mu\in\mathcal{P}^1(\mathbb{T}^1)$, the HJB equation \eqref{eq:HJBwithterminalG} with
terminal condition $G(\cdot,\mu)$ has the unique classical solution
$u^\mu(t,x) \triangleq 2\log q^\mu(t,x)$, where $q^\mu$ is as above. We then solve
the FPK equation forward in time,
\begin{equation}\label{eq:FPKforward}
	\left.
	\begin{aligned}
		\partial_t m(t,x) - \Delta m(t,x) + \operatorname{div}\big(Du^\mu(t,x)\,m(t,x)\big)
		= & \ 0  \qquad \text{ in } (t_0,T)\times\onetorus
		\\
 m(t_0,\cdot) = & \ m_*.
	\end{aligned}
	\right\}
\end{equation}
Set
\[
	n_*\triangleq\frac{m_*}{q^\mu(t_0,\cdot)},
	\qquad n(t,x) \triangleq (n_*\star p_{t-t_0})(x) \quad
	\textrm{ for } t > t_0,
\]
(where by $n_{*}$ we mean the measure
$n_*(A) := \int_A (q^\mu(t_0,x))^{-1}n_*(dx)$ for any
Borelian \(A\subset\onetorus\))
and
\[
	m(t,x) \triangleq q^\mu(t,x)\,n(t,x), \qquad t>t_0, \quad x\in\onetorus.
\]
Since $n(t,x) = (n_*\star p_{t-t_0})(x)$ solves the heat equation
$\partial_t n = \partial_{xx} n$ on $(t_0,T)\times\onetorus$ with
$n(t,\cdot)\to n_*$ weakly as $t\to t_0$, and since $q^\mu$ solves the
backward heat equation, a direct computation shows that $m$ satisfies the
FPK equation \eqref{eq:FPKforward} in the classical sense, with
$m(t,\cdot)\to q^\mu(t_0,\cdot)m_*$ weakly as $t\to t_0$.
By uniqueness of solutions to
\eqref{eq:FPKforward} in the class of probability measures
(see, e.g.~\cite[Thm.~9.8.5]{bogachev2015fokker}, this is the
unique solution; we denote it by $m_\cdot[m_*,\mu]$.
Thus it is natural to define the terminal-map operator
\begin{align*}
	\mathcal{M}\colon \mathcal{P}^1(\mathbb{T}^1)\times\mathcal{P}^1(\mathbb{T}^1)
	& \ \longrightarrow \mathcal{P}(\onetorus),
	\\
	(m_*,\mu)& \ \longmapsto \mathcal{M}(m_*,\mu) \triangleq m_T[m_*,\mu].
\end{align*}
With this notation,
the MFG system on $[t_0,T]$ with starting probability measure $m_*$ is, by
definition, satisfied if and only if there exists $\mu\in\mathcal{P}^1(\mathbb{T}^1)$ such
that
\begin{equation}\label{eq:fixedpointabstract}
	\mathcal{M}(m_*,\mu) = \mu.
\end{equation}
Indeed, if \eqref{eq:fixedpointabstract} holds, then $u = u^\mu$ solves the
HJB equation with terminal condition
$G(\cdot,\mu)=G(\cdot,m_T[m_*,\mu])$, and $m =m_\cdot[m_*,\mu]$ solves the
FPK equation with drift $Du$ and initial condition $m_*$; hence $(u,m)$ is
a solution of the MFG system. Conversely, if $(u,m)$ is a solution of the
MFG system with $m(t_0,\cdot)=m_*$, then, setting $\mu:=m(T,\cdot)$,
uniqueness of the HJB solution with terminal condition $G(\cdot,\mu)$ gives
$u=u^\mu$, and uniqueness of the FPK solution then gives
$m=m_\cdot[m_*,\mu]$; in particular
$\mathcal{M}(m_*,\mu)=m_T[m_*,\mu]=m(T,\cdot)=\mu$, so
\eqref{eq:fixedpointabstract} holds. This establishes the equivalence.
Equivalently, the MFG equation (without a prescribed initial measure) is
satisfied with an absolutely continuous starting measure
if and only if there exists a pair
$(m_*,\mu)\in\mathcal{P}^1(\onetorus)\times\mathcal{P}^1(\mathbb{T}^1)$ satisfying
\eqref{eq:fixedpointabstract}.
Let us write down the operator $\mathcal{M}$ explicitly.
From the formula for $m_\cdot[m_*,\mu]$ we obtain
\begin{align*}
	& \mathcal{M}(m_*,\mu)(x) = m(T,x) = q^\mu(T,x)\,n(T,x)
	\\
	= & \ \big(1+\mathcal{G}(\mu)\cos(2\pi x)\big)
	\bigg(\frac{m_*}{1+\mathcal{G}(\mu)e^{-4\pi^2(T-t_0)}\cos(2\pi\,\cdot)}
	\star p_{T-t_0}\bigg)(x).
\end{align*}
Therefore, the fixed-point condition \eqref{eq:fixedpointabstract}, with
$\mu=m(T,\cdot)$ and hence $\mathcal{G}(\mu)=\mathcal{G}(m(T,\cdot))$, is
precisely
\begin{equation}\label{eq:fixedpoint1}
	m(T,x) =  \big(1+\mathcal{G}(m(T,\cdot))\cos(2\pi x)\big)
	\bigg(\frac{m_{\star}}{1+\mathcal{G}(m(T,\cdot))e^{-4\pi^2(T-t_0)}\cos(2\pi\,\cdot)}
	\star p_{T-t_0}\bigg)(x).
\end{equation}
Before narrowing the class of solutions down to the one
where the family of examples (for different $t_0$) is
exhibited, let us collect the
findings above into a lemma.
For this purpose, for $\theta\in[\theta_1,\theta_2]$ define
\begin{align*}
q^\theta(t,x)&\triangleq1+\theta e^{-4\pi^2(T-t)}\cos(2\pi x),\qquad
u^\theta(t,x)\triangleq2\log q^\theta(t,x),\\
\bar{m}^\theta(t,x)&\triangleq q^\theta(t,x)\int_{\onetorus}
p_{t-t_0}(x-y)\,
\frac{m_*(dy)}{q^\theta(t_0,y)},\qquad t\in(t_0,T],\qquad
\bar{m}^\theta(t_0,\cdot)\triangleq m_* .
\end{align*}
With these notations, the observations given above
constitute the proof of the following lemma.
\begin{lemma}\label{prop:allsolutions}
Let, as above, $F \equiv 0$, $G(x,\mu) \triangleq
2 \log(1+ \mathcal{G}(\mu)\cos(2\pi x))$,
$L(x,v) \equiv \frac{1}{2}|v|^2$, where
$\mathcal{G}$ is a $W_1$-Lipschitz
function taking values in $[\theta_1,\theta_2]\subset (-1,1)$.
Let $t_0\in[0,T)$ and $m_*\in\mathcal P^1(\onetorus)$.
A pair $(u,m)$ is a
solution of the MFG on $[t_0,T]$ with initial distribution $m_*$ if and
only if $(u,m)=(u^\theta,\bar m^\theta)$ for some
$\theta\in[\theta_1,\theta_2]$ satisfying
\begin{equation}\label{eq:generalfixedpoint}
\mathcal{G}\big(\bar m^\theta(T,\cdot)\big)=\theta.
\end{equation}
\end{lemma}

We are going to restrict now to initial measures
(and therefore, to paths of measures) that have a specific form.
Let $0\leq \beta < 1$. If we set
$$
m_* = \frac{\big(1+\mathcal{G}(m(T,\cdot))
e^{-4\pi^2 (T-t_0)}\cos(2\pi x)\big)
(1+\beta\cos(2\pi x))}{1+\frac{\beta\mathcal{G}(m(T,\,\cdot))}{2}e^{-4\pi^2(T-t_0)}},
$$
which is the probability density that makes
the fraction in front of $\star$ simplify exactly to
$(1+\beta\cos(2\pi x))(1+\frac{\mathcal{G}(m(T,\,\cdot))\beta}{2}e^{-4\pi^2(T-t_0)})^{-1}$,
the path of measures becomes
\begin{align*}
	m(t,x) = \frac{(1+ \mathcal{G}(m(T,\,\cdot)) e^{-4\pi^2(T-t)}\cos(2\pi x))
	(1+\beta e^{-4\pi^2(t-t_0)}\cos(2\pi x))}
	{1+\frac{\mathcal{G}(m(T,\,\cdot))\beta}{2}e^{-4\pi^2(T-t_0)}}.
\end{align*}
For the function $\mathcal{G}$ we let
\begin{align*}
	\mathcal{G}(\mu) \triangleq f\bigg( 2\int_{\onetorus} \cos(2\pi x)\,\mu(dx)\bigg),
\end{align*}
where $f:\R\to\R$ is to be determined later,
and we will consider the set $S^{t_0,\theta,\beta}$ (\(0<\theta<1,0<\beta<1\))
consisting
of probability measure paths (in fact, densities)
\begin{align*}
	m^{t_0,\theta,\beta}(t,x)=
	 \frac{(1+ \theta e^{-4\pi^2(T-t)}\cos(2\pi x))
	(1+\beta e^{-4\pi^2(t-t_0)}\cos(2\pi x))}
	{1+\frac{\theta\beta}{2}e^{-4\pi^2(T-t_0)}}
\end{align*}
and restrict the study to MFG solutions whose probability measure path
belongs to this set. Looking back at condition \eqref{eq:fixedpoint1}
if $t_0$ and $\beta$ are fixed, the problem
of finding solutions reduce to the fixed-point problem of finding $\theta$
such that
\begin{align*}
	\mathcal{G}(m^{t_0,\theta,\beta}(T,\cdot)) = \theta.
\end{align*}
Explicitly,
\begin{equation}\label{eq:explicitscalarfixedpoint1}
	\theta = f\bigg( \frac{\theta+\beta e^{-4\pi^2(T-t_0)}}{
		1+\frac{\theta\beta}{2}e^{-4\pi^2(T-t_0)}
	} \bigg).
\end{equation}
Below we will look for a convenient $f$ that allows for distinct pairs $(\theta,\beta)$
and $(\theta',\beta')$ to yield the same initial measure,
i.e.
\begin{align*}
	m^{t_0,\theta,\beta}(t_0,\cdot) = m^{t_0,\theta',\beta'}(t_0,\cdot)
\end{align*}
while simultaneously satisfying the scalar fixed-point equation
\eqref{eq:explicitscalarfixedpoint1}.
If such pairs exist, then the same cost data, horizon, and
initial distribution admit at least two distinct MFG solutions,
which realizes non-uniqueness for a fixed initial datum.
At \(t=t_0\), with
\[
A_{t_0}:=e^{-4\pi^2(T-t_0)},
\]
we have
\begin{align*}
 m^{t_0,\theta,\beta}(t_0,x)  & \ =
\frac{(1+\theta A_{t_0}\cos(2\pi x))(1+\beta\cos(2\pi x))}{1+\frac{\theta\beta A_{t_0}}{2}}
\\
& \ =
\frac{1+(\theta A_{t_0}+\beta)\cos(2\pi x)+\frac{\theta\beta A_{t_0}}{2}(1+\cos(4\pi x))}{1+\frac{\theta\beta A_{t_0}}{2}}
\\
& \ =
1
+
\frac{\theta A_{t_0}+\beta}{1+\frac{\theta\beta A_{t_0}}{2}}\cos(2\pi x)
+
\frac{\frac{\theta\beta A_{t_0}}{2}}{1+\frac{\theta\beta A_{t_0}}{2}}\cos(4\pi x).
\end{align*}
At this point one can see the symmetry that will be used to construct non-uniqueness:
for a given
pair $(\theta,\beta)$,
if we let
$$
\theta' \triangleq\beta/A_{t_0}, \qquad \beta' \triangleq A_{t_0}\theta,
$$
then $m^{t_0,\theta,\beta}(t_0,\cdot)
=m^{t_0,\theta',\beta'}(t_0,\cdot)$
while for $t>t_0$ the two paths differ, provided $\beta\neq A_{t_0}\theta$.
In this case we need, of course, $\beta<A_{t_0}$.
Whether both pairs satisfy the scalar fixed-point equation is a
separate condition to be addressed below.
If the resulting expression for $m^{t_0,\theta,\beta}(t_0,x)$ is denoted by
\begin{align*}
	1+ a_{t_0}\cos(2\pi x) + b_{t_0}\cos(4\pi x),
\end{align*}
then
\begin{equation}\label{eq:prequadratictheta}
	a_{t_0} = \frac{\theta A_{t_0}+\beta}{1+\frac{\theta\beta A_{t_0}}{2}},
	\qquad
	b_{t_0} = \frac{\frac{\theta\beta A_{t_0}}{2}}{1+\frac{\theta\beta A_{t_0}}{2}}.
\end{equation}
If the coefficients of \(\cos(2\pi x)\) and \(\cos(4\pi x)\), that
is $a_{t_0}$ and $b_{t_0}$ are considered for the moment as fixed
parameters, one can solve \eqref{eq:prequadratictheta} for
$\theta$ and $\beta$:
\begin{equation*}
	\theta A_{t_0} + \beta = \frac{a_{t_0}}{1-b_{t_0}}, \qquad
	\theta\beta = \frac{2b_{t_0}}{A_{t_0}(1-b_{t_0})}.
\end{equation*}
Eliminating $\beta$ we obtain a quadratic equation for $\theta$:
\begin{equation}\label{eq:quadratictheta}
	A_{t_0}^2(1-b_{t_0})\theta^2 - a_{t_0} A_{t_0}\theta + 2b_{t_0} =0,
\end{equation}
so we are interested in this equation having two distinct real roots
in the interval \((0,1)\), with the additional constraint that
equation \eqref{eq:explicitscalarfixedpoint1}
be satisfied (for the yet unspecified $f$).

If $a_{t_0}$ and $b_{t_0}$ are considered fixed, and $\theta$
is a solution of \eqref{eq:quadratictheta}, and
$\beta \triangleq  2b_{t_0}/(A_{t_0}\theta(1-b_{t_0}))$, then
the pair \((\theta,\beta)\) satisfies the relations
\eqref{eq:prequadratictheta}. But the pair
\((\theta',\beta')\), with
\begin{align*}
	\theta' = \frac{\beta}{A_{t_0}}, \qquad \beta' = \frac{2b_{t_0}}{\beta(1-\beta)},
\end{align*}
also satisfies them. Thus, for a fixed $t_0$, the two pairs $(\theta,\beta)$ and $(\theta',\beta')$ are
the only pairs that can produce the initial measure
$m^{t_0,\theta,\beta}(t_0,\cdot)$
(in the problem class $S^{t_0,\theta,\beta}$); the second also
belongs to $(0,1)\times (0,1)$ if and only if $\beta<A_{t_0}$.

The initial measure is parametrized by the variables
$(t_0,\beta,\theta)$, but notice that the quantity
\begin{equation}\label{eq:defkappa}
	\kappa\triangleq\theta\beta A_{t_0}=\theta\beta e^{-4\pi^2(T-t_0)}
\end{equation}
is the same for the two pairs, since
$\theta'\beta'A_{t_0}=(\beta/A_{t_0})(A_{t_0}\theta)A_{t_0}=\theta\beta A_{t_0}$.
If $\kappa$ is fixed, then $\beta A_{t_0}=\kappa/\theta$,
and in the new variables
$(\kappa,\beta,\theta)$ the argument of
$f$ in \eqref{eq:explicitscalarfixedpoint1} becomes a function of
$\theta$ alone:
\begin{equation}\label{eq:defhkappa}
	\frac{\theta+\beta A_{t_0}}{1+\frac{\theta\beta A_{t_0}}{2}}
	=\frac{\theta+\kappa/\theta}{1+\kappa/2}
	\triangleq h_\kappa(\theta).
\end{equation}
Consequently, a pair $(\theta,\beta)$ with $\theta\beta A_{t_0}=\kappa$ satisfies
\eqref{eq:explicitscalarfixedpoint1} if and only if
\begin{equation}\label{eq:finvertshkappa}
	f\big(h_\kappa(\theta)\big)=\theta .
\end{equation}
Condition \eqref{eq:finvertshkappa} does not involve $t_0$
explicitly,
$t_0$ is implicit in the relation between the two values of $\theta$:
\begin{equation}\label{eq:thetathetaprime}
	\theta\theta'=\frac{\theta\beta}{A_{t_0}}=\frac{\kappa}{A_{t_0}^2}
	=\kappa e^{8\pi^2(T-t_0)}.
\end{equation}
If $\kappa$ is kept fixed, and $\theta$ is allowed to vary,
this means $t_0$ is allowed to vary, due to
\eqref{eq:defkappa}.
Hence, if $f\circ h_\kappa$ is the identity on a whole interval of values
of $\theta$, then non-uniqueness occurs for every $t_0$ for which that
interval contains two distinct numbers $\theta,\theta'$ satisfying
\eqref{eq:thetathetaprime}. Since
$h_\kappa(\theta)=h_\kappa(\kappa/\theta)$, the function $h_\kappa$ is not
injective on $(0,\infty)$; however,
\[
h_\kappa'(\theta)=\frac{1-\kappa/\theta^2}{1+\kappa/2}>0
\qquad\text{for } \theta>\sqrt{\kappa},
\]
so $h_\kappa$ is strictly increasing on $[\sqrt\kappa,\infty)$,
therefore invertible on any subinterval of $(\sqrt\kappa,1)$.

Suppose now that $\kappa, \theta_1, \theta_2$ are such that
\begin{equation}\label{eq:parameterchoice}
0 < \sqrt{\kappa} < \theta_1 < \theta_2 < 1,
\end{equation}
and
\begin{equation}\label{eq:deff}
\begin{aligned}
	f(s) \triangleq
	\begin{cases}
	\theta_1, & \qquad 	s\leq h_{\kappa}(\theta_1),
	\\
	\frac{1}{2}\big( (1+\kappa/2)s +
	\sqrt{(1+\kappa/2)^2s^2 -4\kappa}\,\big), & \qquad h_{\kappa}(\theta_1)\leq s
	\leq h_{\kappa}(\theta_2),
	\\
	\theta_2, & \qquad s\geq h_{\kappa}(\theta_2).
	\end{cases}
\end{aligned}
\end{equation}
\begin{lemma}\label{lemma:propertiesoff}
For such $\kappa$, $\theta_1$, $\theta_2$ and $f$
as just declared,
the following hold.
\begin{enumerate}[(i)]
\item\label{item:f:inverse} $f$ is well defined, continuous and
nondecreasing, takes values in $[\theta_1,\theta_2]$, and
$f(h_\kappa(\theta))=\theta$ for every $\theta\in[\theta_1,\theta_2]$.
\item\label{item:f:lip} $f$ is Lipschitz continuous with constant
\[
\Lip(f)=\Big(1+\frac{\kappa}{2}\Big)\frac{\theta_1^2}{\theta_1^2-\kappa}.
\]
\item\label{item:f:HS} The function
$\mathcal{G}(\mu)=f\big(2\int_{\onetorus}\cos(2\pi x)\,\mu(dx)\big)$ is
Lipschitz continuous on $\mathcal P^1(\onetorus)$ with respect to $W_1$,
with constant $4\pi\Lip(f)$, and takes values in
$[\theta_1,\theta_2]\subset(0,1)$. In particular, the data $F\equiv0$,
$L(x,v)=\frac12|v|^2$ and
$G(x,\mu)=2\log\big(1+\mathcal{G}(\mu)\cos(2\pi x)\big)$ satisfy the
standing assumptions stated at the beginning of this subsection.
\end{enumerate}
\end{lemma}
\begin{proof}
	As noted above, the function $\theta\mapsto h_{\kappa}(\theta)$
	has a minimum at $\theta = \sqrt{\kappa}$ and is increasing
	from this point. If $\theta_1 \leq \theta\leq \theta_2$,
	and $s=f_{\kappa}(\theta)$, then $\theta^2 - s(1+\kappa/2)\theta + \kappa =0$,
	and $\kappa/\theta$ is another root of this equation,
	which is less than $\theta$. So the inverse of $f_{\kappa}$ on the
	chosen interval corresponds to the positive root of the
	discriminant of this quadratic equation, which is the second case
	in the formula for $f(s)$. All the other statements follow
	with direct computations; part \ref{item:f:HS}
	uses the Kantorovich-Rubenstein duality and part \ref{item:f:lip}.
\end{proof}
Since
$f\circ h_\kappa$ is now the identity on a whole interval of values
of $\theta$, then, as noted above,
non-uniqueness occurs for every $t_0$ for which that
interval contains two distinct numbers $\theta,\theta'$ satisfying
\eqref{eq:thetathetaprime}, i.e.~satisfying $\theta\theta' A_{t_0}^2 = \kappa$.
This will now impose conditions on
the $\kappa$, $\theta_1$, $\theta_2$ that we are to fix.
Let us now justify that these conditions exist such that one can find
$\sigma\in (0,\sqrt{\kappa})$ and
$t_0\in [0,T)$ with
\begin{equation}\label{eq:conditiont0sigma}
\theta_1\leq\sigma e^{4\pi^2(T-t_0)}, \qquad
\frac{\kappa}{\sigma}e^{4\pi^2(T-t_0)}\leq\theta_2;
\end{equation}
these two numbers will then be \(\theta\) and \(\theta'\)
(and note that $\sigma<\sqrt{\kappa}$ gives $\sigma<\kappa/\sigma$,
so in that case, $\theta_1\leq\sigma e^{4\pi^2(T-t_0)}$
\(<\tfrac{\kappa}{\sigma}e^{4\pi^2(T-t_0)}\leq\theta_2\)).
We claim that this occurs if and only if
\[
\theta_1<\sqrt\kappa\,e^{4\pi^2(T-t_0)}<\theta_2.
\]
Indeed, conditions \eqref{eq:conditiont0sigma}
read
$\sigma\geq \theta_1 A_{t_0}$ and
\(\sigma\geq \tfrac{\kappa}{A_{t_0}\theta_2}\),
and such a $\sigma$ in the interval \((0,\sqrt{\kappa})\)
exists if and only \(\theta_1A_{t_0}<\sqrt{\kappa}\)
and $\tfrac{\kappa}{A_{t_0}\theta_2}<\sqrt{\kappa}$.

For the next proposition,
if $\sigma\in(0,\sqrt\kappa)$ set
\begin{equation}\label{eq:defmusigma}
\mu_\sigma(x)\triangleq
\frac{\big(1+\sigma\cos(2\pi x)\big)
\big(1+\frac{\kappa}{\sigma}\cos(2\pi x)\big)}
{1+\frac{\kappa}{2}},\qquad x\in\onetorus.
\end{equation}

\begin{proposition}\label{prop:nonuniqueness}
Assume that $0<\sqrt{\kappa}<\theta_1 <\theta_2<1$
and $t_0$ is such that
\[
\theta_1<\sqrt\kappa\,e^{4\pi^2(T-t_0)}<\theta_2,
\]
so that, as explained above, there is $\sigma\in (0,\sqrt{\kappa})$
satisfying \eqref{eq:conditiont0sigma}, i.~e.
\begin{equation*}
\theta_1\leq\sigma e^{4\pi^2(T-t_0)},
\qquad
\frac{\kappa}{\sigma}\,e^{4\pi^2(T-t_0)} \leq \theta_2 ;
\end{equation*}
let $f$ be as defined above in \eqref{eq:deff}.
Then $\mu_\sigma$, as
defined in \eqref{eq:defmusigma} is a probability density,
and the MFG on $[t_0,T]$
with initial distribution $\mu_\sigma$ has the two solutions
$(u_1,m_1)$ and $(u_2,m_2)$ given by
\begin{align*}
u_1(t,x)&=2\log\big(1+\sigma e^{4\pi^2(t-t_0)}\cos(2\pi x)\big),
&m_1&=m^{t_0,\theta,\beta},
&&\theta=\sigma e^{4\pi^2(T-t_0)},\ \ \beta=\frac{\kappa}{\sigma},\\
u_2(t,x)&=2\log\Big(1+\frac{\kappa}{\sigma} e^{4\pi^2(t-t_0)}\cos(2\pi x)\Big),
&m_2&=m^{t_0,\theta',\beta'},
&&\theta'=\frac{\kappa}{\sigma}e^{4\pi^2(T-t_0)},\ \ \beta'=\sigma .
\end{align*}
These solutions are distinct, and their values at the initial time,
\[
u_1(t_0,\cdot)=2\log\big(1+\sigma\cos(2\pi\,\cdot)\big),\qquad
u_2(t_0,\cdot)=2\log\bigg(1+\frac{\kappa}{\sigma}\cos(2\pi\,\cdot)\bigg),
\]
are distinct and do not depend on $t_0$.
\end{proposition}
\begin{proof}
We let $\theta=\sigma/A_{t_0}$ and
$\theta'=\kappa/(\sigma A_{t_0})$.
Since $\sigma<\sqrt\kappa<\kappa/\sigma$, we
have $\theta<\theta'$, and \eqref{eq:conditiont0sigma} gives
\[
\theta_1\leq\theta<\theta'\leq\theta_2 .
\]
In particular $\theta,\theta'\in(0,1)$. Moreover $\beta'=\sigma<1$ and
$\beta=\kappa/\sigma=A_{t_0}\theta'\leq\theta'<1$. Hence both pairs lie in
$(0,1)^2$, all four factors in the definitions of $m^{t_0,\theta,\beta}$
and $m^{t_0,\theta',\beta'}$ are positive, and both are paths of
probability densities. The pair $(\theta',\beta')$ is obtained from
$(\theta,\beta)$ by the substitution $\theta'=\beta/A_{t_0}$, $\beta'=A_{t_0}\theta$
considered above, and
\[
\theta\beta A_{t_0}=\theta'\beta'A_{t_0}=\kappa .
\]
Regarding the initial measure, since $\theta A_{t_0}=\sigma$ and
$\theta'A_{t_0}=\kappa/\sigma$,
\[
m^{t_0,\theta,\beta}(t_0,x)
=\frac{\big(1+\sigma\cos(2\pi x)\big)\big(1+\frac{\kappa}{\sigma}\cos(2\pi x)\big)}{1+\frac\kappa2}
=m^{t_0,\theta',\beta'}(t_0,x),
\]
and this is $\mu_\sigma(x)$. In particular $\mu_\sigma$ is a probability
density.

As for the fixed point equation, by \eqref{eq:defhkappa} the argument of
$f$ in \eqref{eq:explicitscalarfixedpoint1} for the pair $(\theta,\beta)$
equals $h_\kappa(\theta)$, and for the pair $(\theta',\beta')$ it equals
$h_\kappa(\theta')$. Since $\theta,\theta'\in[\theta_1,\theta_2]$,
Lemma \ref{lemma:propertiesoff}\ref{item:f:inverse} yields
$f(h_\kappa(\theta))=\theta$ and $f(h_\kappa(\theta'))=\theta'$, that is,
$\mathcal{G}(m_1(T,\cdot))=\theta$ and $\mathcal{G}(m_2(T,\cdot))=\theta'$.

For the claim on solutions, recall that by the construction at the beginning of this subsection
(with $q(t,x)=1+\theta e^{-4\pi^2(T-t)}\cos(2\pi x)$ and
$n(t,x)=(1+\beta e^{-4\pi^2(t-t_0)}\cos(2\pi x))/(1+\frac\kappa2)$),
the function
\[
u_1(t,x)=2\log\big(1+\theta e^{-4\pi^2(T-t)}\cos(2\pi x)\big)
=2\log\big(1+\sigma e^{4\pi^2(t-t_0)}\cos(2\pi x)\big)
\]
is the solution of the Bellman equation with terminal value
$2\log(1+\theta\cos(2\pi x))=G(x,m_1(T,\cdot))$, and $m_1=qn$ solves the
FPK equation with drift $Du_1$ in the classical sense on
$(t_0,T)\times\onetorus$ and hence in the sense of distributions;
moreover $m_1$ is smooth on $[t_0,T]\times\onetorus$, so that
$t\mapsto m_1(t,\cdot)$ is continuous into $\mathcal P^1(\onetorus)$, and
$m_1(t_0,\cdot)=\mu_\sigma$. Thus $(u_1,m_1)$ is a solution of the MFG
on $[t_0,T]$ with initial distribution $\mu_\sigma$. The same argument
applies to $(u_2,m_2)$.

Finally, $u_1(t_0,\cdot)$ and $u_2(t_0,\cdot)$ are the
functions displayed in the statement; they differ because
$\sigma\neq\kappa/\sigma$, and they do not depend on $t_0$.
\end{proof}

\begin{corollary}[A more specific example]\label{cor:allt0}
Let $\varepsilon\in(0,T)$ and choose
\[
\sqrt\kappa=\tfrac12e^{-4\pi^2T},\qquad
\theta_1=\sqrt\kappa\,e^{4\pi^2\varepsilon}=\tfrac12e^{-4\pi^2(T-\varepsilon)},\qquad
\theta_2=\tfrac34 .
\]
Then \eqref{eq:parameterchoice} holds, and for the resulting terminal cost
$G$ the following is true: for every $t_0\in[0,T-\varepsilon)$ there is an
initial distribution from which the MFG on $[t_0,T]$ has at least two
distinct solutions. The Lipschitz constant of $f$ is
\[
\Lip(f)=\frac{1+\kappa/2}{1-e^{-8\pi^2\varepsilon}} .
\]
\end{corollary}
\begin{proof}
Clearly $0<\sqrt\kappa<\theta_1<\frac12<\theta_2<1$. For
$t_0\in[0,T-\varepsilon)$,
\[
\sqrt\kappa\,e^{4\pi^2(T-t_0)}=\tfrac12e^{-4\pi^2t_0}\in
\big(\tfrac12e^{-4\pi^2(T-\varepsilon)},\tfrac12\big]\subset(\theta_1,\theta_2),
\]
so the condition \(\theta_1<\sqrt\kappa\,e^{4\pi^2(T-t_0)}<\theta_2\)
of Proposition \ref{prop:nonuniqueness} holds, implying
that \eqref{eq:conditiont0sigma} is fulfilled and therefore
the conclusion the proposition applies.
The formula for $\Lip(f)$ follows from Lemma
\ref{lemma:propertiesoff}\ref{item:f:lip} and
$\theta_1^2=\kappa e^{8\pi^2\varepsilon}$.
\end{proof}

\begin{remark}\label{rem:originalexample}
The example found by DeepSeek corresponds to $T=1$, $t_0=0$,
$\theta=0.4$, $\theta'=0.6$, hence
$\kappa=\theta\theta'e^{-8\pi^2}=0.24\,e^{-8\pi^2}$ by
\eqref{eq:thetathetaprime}, $\sigma=0.4\,e^{-4\pi^2}$ and
$\kappa/\sigma=0.6\,e^{-4\pi^2}$.
There $f$ was a
continuous piecewise function,
constant on either side of an interval, chosen so that
\eqref{eq:finvertshkappa} holds at the two points $\theta=0.4$ and
$\theta=0.6$ only. It was clear in the example that the idea was
to obtain solutions in which the initial measure is the product of two functions
in which the coefficients could be switched, without altering
the initial measure. We took it upon ourself to take
the example a bit further, in a way that would make it work
for different values of $t_0$.

Namely, by \eqref{eq:thetathetaprime}, a different initial
time requires a different pair $\theta,\theta'$, and therefore new points
at which \eqref{eq:finvertshkappa} must hold. So the terminal
cost would have to be modified if one wanted to change
the initial time. The
choice \eqref{eq:deff} makes \eqref{eq:finvertshkappa} hold on the whole
interval $[\theta_1,\theta_2]$.

Thus, the general framework of this section, together
with all the statements (lemmas, propositions, remarks) are the work
of the authors.
\end{remark}

We conclude this section with an observation on uniqueness
for times $t_0$ close to $T$. Let us return to the more
general context of Lemma \ref{prop:allsolutions}.

When $\mathcal G(\mu)=f\big(2\int\cos(2\pi x)\,\mu(dx)\big)$, the
fixed-point equation \eqref{eq:generalfixedpoint} can also be written
explicitly. With \(A_{t_0}=e^{-4\pi^2(T-t_0)}\) as before, Fubini's theorem, the
symmetry of the function \(p_s$, the identities
$2\cos(2\pi x)q^\theta(T,x)=2\cos(2\pi x)+\theta\big(1+\cos(4\pi x)\big)\)
and $\cos(2\pi k\,\cdot)\star p_s=e^{-4\pi^2k^2s}\cos(2\pi k\,\cdot)$
give us
\begin{equation}\label{eq:terminalcoefficient}
	\begin{aligned}
2\int_{\onetorus}\cos(2\pi x)\,\bar m^\theta(T,x)\,dx
= & \
\int_{\onetorus}\frac{2A\cos(2\pi y)+\theta\big(1+A^4\cos(4\pi y)\big)}
{1+\theta A\cos(2\pi y)}\,m_*(dy)
\\
\triangleq  & \ B[\theta;m_{*}],
	\end{aligned}
\end{equation}
so that \eqref{eq:generalfixedpoint} reads
\begin{displaymath}
\theta=f(B[\theta;m_*]).
\end{displaymath}
For
$m_*=m^{t_0,\theta,\beta}(t_0,\cdot)$ one recovers the argument of $f$ in
\eqref{eq:explicitscalarfixedpoint1}.

\begin{proposition}\label{prop:uniquenessnearT}
Let $\mathcal{G}(\mu)=f\big(2\int_{\onetorus}\cos(2\pi x)\,\mu(dx)\big)$, where
$f:\R\to[\theta_1,\theta_2]\subset(-1,1)$ is Lipschitz continuous, and set
$\bar\theta\triangleq\max\{|\theta_1|,|\theta_2|\}$. If $t_0\in[0,T)$ satisfies
\begin{equation}\label{eq:uniquenesscondition}
\Lip(f)\,\big(1-e^{-16\pi^2(T-t_0)}\big)<(1-\bar\theta)^2,
\end{equation}
then for every $m_*\in\mathcal P^1(\onetorus)$ the MFG on $[t_0,T]$ with
initial distribution $m_*$ has exactly one solution.
\end{proposition}
\begin{proof}
By Lemma \ref{prop:allsolutions}, the solutions are in one-to-one
correspondence with the solutions $\theta\in[\theta_1,\theta_2]$ of
$\theta=f(B[\theta;m_*])$, with distinct \(\theta\) giving
distinct $u^\theta$. The
map \(\theta\mapsto f(B[\theta,m_*])\) takes the interval $[\theta_1,\theta_2]$ into
itself; the integrand in \eqref{eq:terminalcoefficient} is smooth in
$\theta$, and a direct computation using $2\cos^2(2\pi y)=1+\cos(4\pi y)$ results in
\[
\frac{\partial}{\partial\theta}\,
\frac{2A_{t_0}\cos(2\pi y)+\theta\big(1+A_{t_0}^4\cos(4\pi y)\big)}{1+\theta A_{t_0}
\cos(2\pi y)}
=
\frac{(1-A_{t_0}^2)\big(1-A_{t_0}^2\cos(4\pi y)\big)}{\big(1+\theta A_{t_0}\cos(2\pi y)\big)^2}.
\]
The absolute value of the right-hand side is not larger than
$(1-A_{t_0}^2)(1+A_{t_0}^2)/(1-\bar\theta A_{t_0})^2\leq(1-A_{t_0}^4)/(1-\bar\theta)^2$.
Since $m_*$ is a probability measure,
\[
|B[\theta;m_*]-B[\tilde\theta;m_{*}]|
\leq \frac{1-A_{t_0}^4}{(1-\bar\theta)^2}\,|\theta-\tilde\theta|,
\qquad \text{ any } \theta, \tilde\theta \in [\theta_1,\theta_2]
\]
obtains, so $\theta\mapsto B[\theta;m_*]$ is a contraction map of
$[\theta_1,\theta_2]$ under the condition \eqref{eq:uniquenesscondition}, since
$A_{t_0}^4=e^{-16\pi^2(T-t_0)}$. By the Banach fixed-point theorem it has
exactly one fixed point.
\end{proof}

Thus, Proposition \ref{prop:uniquenessnearT} shows that, for a fixed terminal
cost of the form considered here cannot produce
non-uniqueness for all initial times in $[0,T)$ because, since
$1-e^{-16\pi^2(T-t_0)}\to0$ as $t_0\to T$, condition
\eqref{eq:uniquenesscondition} holds for all $t_0$ sufficiently close to
$T$, whatever the value of $\Lip(f)$ is.

	\section{Control theory: relaxed controls,
	value functions and estimates}\label{section:relaxedcontrols}

As we justify below in Section \ref{section:MFGsolsandvaluemf},
the MFG system of equations \eqref{intro:eq:mfg_system} can be
seen as the optimality conditions
of a control problem. In this
section we present three formulations of this problem that
are equivalent.

\textit{Note}. The functions $f=f(t,x)$ and $\sigma=\sigma(x)$ that will appear below,
for the time being, are functions that
satisfy the minimal conditions
for the expressions in which they appear to be meaningful;
concrete conditions on them will be placed on an as-needed basis.
After this section the functions
$(t,x)\mapsto F(x,m_t)$ will take the place of $f(t,x)$.
For the moment $f$ and $\sigma$ may be assumed, say, Borel and bounded.

	\subsection{Feedback formulation, control
	processes}\label{subsection:feedbackformulation}
	We begin by introducing the class of stochastic control problems:
	\begin{equation}\label{prel:aux_prob:criterion}\text{maximize }
		\mathbb{E}\Bigg[\int_s^r\big[f(t,X_t)-L(X_t,w(t,X_t))\big]dt
		+\sigma(X_r)\Bigg]\end{equation}
	subject to the SDE
	initial value problems:
	\begin{equation}\label{prel:aux_prob:dynamics}
		\left.
			\begin{aligned}
		dX_t = & \ w(t,X_t)dt+\sqrt{2} dW_t,
		\\
		 X_s =  & \ \xi.
			\end{aligned}
		\right\}
	\end{equation}
	The time instants $s,r$ are taken from the interval $[0,T]$,
	with \(s\leq r\).

   The function $w:[s,r]\times\td\rightarrow \rd$,
    a ``velocity vector field'', that appears
    in the SDE \eqref{prel:aux_prob:dynamics}
    is always assumed to be Borel measurable and bounded
    and we will call such a function
    a \textit{feedback control}.

    When we refer to a \emph{solution} of the
    SDE \eqref{prel:aux_prob:dynamics} we
    assume that such a solution is set in \emph{some}
	filtered probability space
	$(\Omega,\mathcal{F},(\mathcal{F}_t)_{s\leq t\leq r},\mathbb{P})$
	with a \(d\)-dimensional \((\mathcal{F}_t)_{s\leq t\leq r}\)-Wiener process $(W_t)_{t\geq s}$
	on it (in particular, \(W\) is \((\mathcal{F}_t)\)-adapted and \(W_t-W_u\) is
	independent of \(\mathcal{F}_u\) for \(s\leq u\leq t\)),
    with the equality in \eqref{prel:aux_prob:dynamics}
    understood
    in the $\mathbb{P}$-a.e.~sense; we don't assume
    that the probability space,
    filtration and Wiener process are
    fixed (which is why we some times favor
    using the notation \(\operatorname{Law}(\xi)\) over
    writing this distribution as the pushforward of a concrete
    probability measure by the mapping $\xi$). The solution \((X_t)_{s\leq t\leq r}\)
    is by definition
    $(\mathcal{F}_t)_{s\leq t\leq r}$-adapted.
	The (deterministic) functions
	$f:[s,r]\times\dtorus\to \R$
	and
	\(L:\dtorus\times \R^d\to \R\) (Lagrangian)
	represent some running payoff and
	cost respectively,
	while $\sigma$ represents a terminal payoff.

It is known \cite{Veretennikov1981} that
for each feedback control $w$ and initial
distribution $\xi$, the SDE \eqref{prel:aux_prob:dynamics} has a solution
$(X_t)_{t\geq s}$.
Its law (distribution) $m_t \triangleq \operatorname{Law}(X_t)$
satisfies the FPK
equation
\begin{equation}\label{prel:eq:FPK}
	\partial_tm_t+\operatorname{div}(w(t,x)m_t)=\triangle m_t
\end{equation} with the initial condition
\begin{equation*}
	m_s=m_*,
\end{equation*}
where
\begin{equation}
m_* = \operatorname{Law}(\xi).
\end{equation}
Equation \eqref{prel:eq:FPK} is understood in the sense that
\begin{equation}\label{eq:FPKsolcond}
\int_{s}^r \int_{\dtorus} \big( \partial_t \phi(t,x) + w(t,x)\cdot D \phi(t,x)
 + \Delta \phi(t,x) \big) m_t(dx) dt = 0
\end{equation}
for every \(\phi\in \mathcal{C}_c^{\infty}((s,r)\times \dtorus)\).
We say in this case that the FPK equation is satisfied
in the sense of distributions or in the weak sense.
A key and well-known observation is that for the functional
in \eqref{prel:aux_prob:criterion} we have
\[
\begin{split}
	\mathbb{E}\Bigg[\int_s^r\big[f(&t,X_t)-L(X_t,w(t,X_t))\big] dt
	+\sigma(X_r)\Bigg]\\&=
	\int_s^r\int_{\td}[f(t,x)-L(x,w(t,x))\big]m_t(dx)dt
	+\int_{\td}\sigma(x)m_r(dx),
\end{split}
\]
which leads us to the following definition.
\begin{definition}\label{prel:def:admissible}
	We say that a pair $(m_\cdot,w)$, where $m:[s,r]\rightarrow
	\mathcal{P}^1(\td)$ is weakly measurable
and
$w:[s,r]\times \td\rightarrow\rd$ is a measurable and bounded
(i.e.~a feedback control),
is a \emph{control process on $[s,r]$}, if it satisfies
\eqref{prel:eq:FPK} in the weak sense. The set of control
processes is denoted by $\Adm{s}{r}$.
The set of control processes with velocities that are essentially bounded by
a specific constant $C$ is denoted by $\Adm[C]{s}{r}$, i.e.,
\[
\Adm[C]{s}{r}=\bigg\{(m_\cdot,w)\in\Adm{s}{r}: \ \
|w(t,x)|\leq C\text{ for }m_t\text{-a.e. }x\in\td\text{ and a.e. }
s\leq t\leq r\bigg\}.
\]
\end{definition}
We note that the meaning of our term ``control process''
differs from that which is usually assigned in
stochastic control theory, as our emphasis here
is on the FPK equation.

If $(m_\cdot,w)$ is a control process, then we denote
\begin{equation}\label{eq:defnofIa}
	I^a[s,r,m_\cdot,w,f,\sigma]\triangleq \int_s^r\int_{\td}
[f(t,x)-L(x,w(t,x))\big]m_t(dx)dt+\int_{\td}\sigma(x)m_r(dx).
\end{equation}
Since the law of any solution of the SDE \eqref{prel:aux_prob:dynamics}
solves \eqref{prel:eq:FPK},
maximizing the functional
in \eqref{prel:aux_prob:criterion}
over all feedback controls $w$
bounded by a constant $C$
yields a value
at least as large as that of maximizing
$I^a[s,r,m_{\cdot},w,f,\sigma]$ over all control processes $(m_{\cdot},w)$
in \(\Adm[C]{s}{r}\).
On the other hand, due to the
superposition principle (see \cite{Figalli2008}),
given a narrowly continuous
solution $m_{\cdot}$ to the FPK equation
\eqref{prel:eq:FPK} with initial measure $m_*$, one can construct (we omit the details)
a solution $(X_t)_{s\leq t\leq r}$ to the SDE
\eqref{prel:aux_prob:dynamics} with an initial
random variable $\xi$ whose distribution (law)
is $m_*$, i.e.~satisfying
$\xi\sharp\mathbb{P}=m_*$ and
such that the distribution of $X_t$ is $m_t$,
i.e.~$X_t\sharp\mathbb{P}=m_t$, $s\leq t\leq r$.
Therefore the two maximization problems yield
the same value:
the stochastic control problem
\eqref{prel:aux_prob:criterion}-\eqref{prel:aux_prob:dynamics}
with feedback controls bounded by $C$
is equivalent
to the maximization of
$I^a[s,r,m_\cdot,w,f,\sigma]$ over the set $\Adm[C]{s}{r}$.

Below, in Lemma \ref{prel:lm:verification} and
Proposition \ref{prop:prel:boundness}, we show that,
if $f$ and $\sigma$ are regular enough
(Lipschitz continuous w.r.t.~$x$, with an additional
mild regularity requirement stated there),
then the supremum of $I^a$ is attained by a feedback control
bounded by a constant $C$
that depends only on the Lipschitz constants of $f$, $\sigma$
and on the data of $L$; hence one can a priori consider the
functional $I^a$
(since we are concerned in its maximization)
only on $\Adm[C]{s}{r}$.
We introduce the value function
\begin{equation}\label{prel:intro:value}
	\operatorname{Val}[s,r,m_*,f,\sigma]\triangleq
	\sup\Big\{I^a[s,r,m_\cdot,w,f,\sigma]:\, (m_\cdot,w)
	\in\Adm{s}{r},\ \ m_s=m_*\Big\}.
\end{equation}
Thus, by the results just referenced,
given the conditions on $f$ and $\sigma$ named there,
eventually one can write this definition
with \((m_{\cdot},w)\in\Adm[C]{s}{r}\) rather than
\((m_{\cdot},w)\in\Adm{s}{r}\).

Abusing notation somewhat we will denote, for $y\in\td$,
\begin{equation*}
	\operatorname{Val}[s,r,y,f,\sigma]\triangleq
	\operatorname{Val}[s,r,\delta_y,f,\sigma].
\end{equation*}
The case $s=r$ is $\operatorname{Val}[s,s,m_{*},f,\sigma]=
\action{m_*}{\sigma}$, in particular,
$\operatorname{Val}[s,s,x,f,\sigma]=\sigma(x)$.

\subsection{Relaxed control processes}\label{subsection:relaxedcontrolprocesses}
It will be convenient to \emph{relax} the control processes such that
the set of relaxed control processes,
endowed with a suitable topology, enjoys compactness.
We will now introduce relaxed controls and show
that the value function
defined by \eqref{prel:intro:value} can be equivalently
obtained through a maximization on relaxed controls
(Lemma \ref{lemma:equality:value_RAdm} below).
\begin{definition}\label{prel:def:relaxed_control_procs}
	Let \(0\leq s<r\leq T\), \(C>0\). The
set of \emph{relaxed controls on \([s,r]\)},
\(\RAdm[C]{s}{r}\)
is the set of Borel measures \(\gamma\) on
\([0,T]\times \mathbb{T}^d \times \mathbb{B}_C\)
for which
\begin{enumerate}
  \item \(\spt\gamma \subset [s,r]\times\dtorus\times\mathbb{B}_C\),
\item  \(\gamma([s,r]\times\dtorus\times\mathbb{B}_C) = r-s\),
\item \(\operatorname{p}^1\sharp \gamma = \mathscr{L}^1_{[s,r]}\),
where \(\mathscr{L}^1_{[s,r]}\) is defined by the equalities
\(\mathscr{L}^1_{[s,r]}(B) = \mathscr{L}^1(B\cap[s,r])\),
with \(B\subset [0,T]\) a Borel subset and
\(\mathscr{L}^1\) the Lebesgue measure on the real line,
\item for all
\(\varphi \in \mathcal{C}_c^{1,2}((s,r)\times \mathbb{T}^d)\):
\begin{equation}\label{eq:weakcontrolFPK}
  \int\limits_{[0,T]\times \dtorus\times \mathbb{B}_C}
  \big(  \partial_t \varphi(t,x) + D \varphi (t,x) \cdot v +
  \Delta_x \varphi(t,x) \big) \gamma(dt,dx,dv)  = 0;
\end{equation}
here \(\mathcal{C}_c^{1,2}((s,r)\times\dtorus)\) is
the set of functions on
\([0,T]\times\dtorus\)
whose support is included in a compact subset of
\((s,r)\times \dtorus\),
with continuous derivative in \(t\) and continuous second-order
derivatives in the variable \(x \in \dtorus\).
\end{enumerate}
By \(\RAdm{s}{r}\) we will mean
\begin{align*}
    \RAdm{s}{r} \triangleq \bigcup_{C>0} \RAdm[C]{s}{r}.
\end{align*}
\end{definition}
If $\gamma\in\RAdm{s}{r}$, then we denote by $\mathfrak{v}[\gamma]$
the vector field defined by
\[
\mathfrak{v}[\gamma](t,x)\triangleq \int_{\mathbb{B}_C}v\gamma(dv|t,x).
\]
Here $\gamma(\cdot|t,x)$ denotes the disintegration of $\gamma$
w.r.t.~$\vartheta\triangleq \operatorname{p}^{1,2}\sharp\gamma$.
Furthermore, for $s\leq t\leq r$, we define the operator $\mathfrak{f}_t$
from $\RAdm{s}{r}$ to $\mathcal{P}^1(\td)$ by the following rule:
if $\gamma\in \RAdm{s}{r}$, first set $\vartheta\triangleq
\operatorname{p}^{1,2}\sharp\gamma$, then let $\mathfrak{f}_t[\gamma]$
be equal to $m_t$ where $m_\cdot$ is the disintegration of $\vartheta$
with respect to the Lebesgue measure. Notice that at this point
$\mathfrak{f}$ is defined up
to a set of zero Lebesgue measure on \([s,r]\),
but this is remedied by Proposition
\ref{prel:prop:Holder_m} below.
The proofs, which are rather standard,
of the properties of relaxed controls that are given in
Propositions \ref{prel:prop:link_adm} and
\ref{prel:prop:Holder_m} are deferred to Appendix A.
Control processes, and relaxed control processes
(or, simply, \textit{relaxed controls}) are in correspondence
with one another:
\begin{proposition}\label{prel:prop:link_adm}
	If $(m_\cdot,w)\in\Adm{s}{r}$, then the measure $\gamma$ defined
	by the rule: for each $\phi\in \mathcal{C}_b([s,r]\times\td\times\rd)$,
	\begin{equation}\label{prel:intr:adm_embedding}
		\begin{split}
    \action{\gamma}{\phi} = & \  \int_{[s,r]\times
		\td\times\rd}\phi(t,x,v)\gamma(dt,dx,dv)
        \\
		\triangleq & \
		\int_{s}^r\int_{\td}\int_{\rd}\phi(t,x,v) \delta_{w(t,x)}(dv)
		m_t(dx)dt
		\end{split}
        \end{equation}
        lies in $\RAdm{s}{r}$.
	Conversely, if $\gamma\in\RAdm{s}{r}$ then
	$(\mathfrak{f}_\cdot[\gamma],\mathfrak{v}[\gamma])\in\Adm{s}{r}$.
\end{proposition}
The same holds for those that are bound by the constant \(C\):
\begin{corollary}\label{prel:corollary:link_adm}
If $(m_\cdot,w)\in\Adm[C]{s}{r}$, then $\gamma$,
defined by \eqref{prel:intr:adm_embedding}
lies in $\RAdm[C]{s}{r}$.
Conversely,
$(\mathfrak{f}_\cdot[\gamma],\mathfrak{v}[\gamma])\in\Adm[C]{s}{r}$
provided that $\gamma\in\RAdm[C]{s}{r}$.
\end{corollary}
It will be useful to know that the path of probabilities
$\mathfrak{f}_\cdot[\gamma]$, i.e.~\(t\mapsto
\mathfrak{f}_t[\gamma]\), is continuous:
\begin{proposition}\label{prel:prop:Holder_m}
	For every $\gamma\in\RAdm[C]{s}{r}$,
	there exists a \(\frac{1}{2}\)-H\"older
    continuous version of $\{\mathfrak{f}_t[\gamma]\}_{s\leq t\leq r}$.
    More precisely: there exists a family \(\{\mu_t\}_{s\leq t\leq r}
    \subset \mathcal{P}^1(\dtorus)\) that coincides with
    $\{\mathfrak{f}_t[\gamma]\}_{s\leq t\leq r}$ a.e.~with respect
    to the Lebesgue measure on $[s,r]$ and
	$t\mapsto\mu_t$ is $\frac{1}{2}$-H\"older
	continuous in the metric $W_1$ with H\"older constant
	$C\sqrt{T}+\sqrt{2d}$.
\end{proposition}
Since $\{\mathfrak{f}_t[\gamma]\}_{s\leq t\leq r}$
coincides with the continuous family
\(\{\mu_t\}_{s\leq t\leq r}
    \subset \mathcal{P}^1(\dtorus)\) a.e.,
    this family is also a valid
    disintegration of $\pj^{1,2}\sharp\gamma$ with respect
    to the Lebesgue measure on $[s,r]$.
\emph{From now on $\mathfrak{f}_\cdot[\gamma]$ always denotes
this continuous representative}; in particular $\mathfrak{f}_s[\gamma]$ and
$\mathfrak{f}_r[\gamma]$ are well defined.

Also, Proposition \ref{prel:prop:Holder_m} enables us to consider $\mathfrak{f}_\cdot$
as an operator from the space $\RAdm[C]{s}{r}$ to the space
$\mathcal{C}^{0}([s,r];\mathcal{P}^1(\td))$
or even $\mathcal{C}^{0,\frac{1}{2}}([s,r];\mathcal{P}^1(\td))$ (H\"older continuous in $t$ paths
in \(\mathcal{P}^1(\dtorus)\)), but it will be enough
to stick with the former.
On $\mathcal{C}^{0}([s,r];\mathcal{P}^1(\td))$
we assume the topology generated by the
usual supremum distance.
\begin{proposition}\label{prel:prop:continuity}
	The mapping
	$\mathfrak{f}_\cdot:\RAdm[C]{s}{r}\rightarrow
	\mathcal{C}^0([s,r];\mathcal{P}^1(\td))$ is continuous.
\end{proposition}
The proofs of Propositions \ref{prel:prop:Holder_m} and
\ref{prel:prop:continuity} are provided in the Appendix.
We will now introduce two operations with relaxed control
processes: concatenation and restriction.
\begin{definition}\label{prel:def:concatenation}
	For any \(C>0\), if
	$s<r<\theta$, $\gamma_1\in \RAdm[C]{s}{r}$,
	$\gamma_2\in\RAdm[C]{r}{\theta}$
	are such that $\mathfrak{f}_r[\gamma_1]
	=\mathfrak{f}_r[\gamma_2]$,
	then we denote by $\gamma_1\diamond\gamma_2$
	the concatenation of the
	relaxed control
	processes $\gamma_1$ and $\gamma_2$ defined by the rule: for
	$\phi\in \mathcal{C}_b([s,\theta]\times\td\times\rd)$,
    $\action{\gamma_1\diamond\gamma_2}{\phi}
    \triangleq
    \action{\gamma_1}{\phi}+\action{\gamma_2}{\phi}$,
    that is:
	\[
	\begin{split}
		\int_{[s,\theta]\times
		\td\times\rd}\phi(t,x,v)(&\gamma_1\diamond\gamma_2)(dt,dx,dv)
		\\
        \triangleq
		\int_{[s,r]\times \td\times\rd}&\phi(t,x,v)\gamma_1(dt,dx,dv)
		+\int_{[r,\theta]\times \td\times\rd}\phi(t,x,v)\gamma_2(dt,dx,dv).
	\end{split}
	\]
\end{definition}
By construction, we have that
\[
\gamma_1\diamond\gamma_2\in\RAdm[C]{s}{\theta};
\]
the only point requiring care is the verification of
\eqref{eq:weakcontrolFPK} for test functions whose support crosses the
hyperplane \(t=r\): approximating such a test function by test
functions vanishing near \(t=r\), the two boundary terms that arise are
\(\pm\action{\mathfrak{f}_r[\gamma_i]}{\varphi(r,\cdot)}\), \(i=1,2\),
which cancel by the matching condition
\(\mathfrak{f}_r[\gamma_1]=\mathfrak{f}_r[\gamma_2]\)
and the continuity of the flows (Proposition
\ref{prel:prop:Holder_m}). We omit the details.
\begin{definition}\label{prel:def:restriction} If $s\leq s'<r'\leq r$,
	$\gamma\in\RAdm{s}{r}$, then denote by $\gamma|_{[s',r']}$
	the restriction of the relaxed control process $\gamma$ on $[s',r']$
	that is defined by the rule: for $\phi\in \mathcal{C}_b([s',r']\times\td\times\rd)$,
	\[
	\int_{[s',r']\times \td\times\rd}\phi(t,x,v)(\gamma|_{[s',r']})(dt,dx,dv)
	\triangleq
	\int_{[s,r]\times \td\times\rd}\tilde{\phi}(t,x,v)\gamma(dt,dx,dv),
	\]
	where
	\[
	\tilde{\phi}(t,x,v)\triangleq \left\{\begin{array}{ll}
		\phi(t,x,v), & t\in [s',r'],
		\\
		0,&\text{otherwise}.
	\end{array}\right.
	\]
\end{definition}
Obviously, if $\gamma\in\RAdm{s}{r}$ (respectively,
$\gamma\in \RAdm[C]{s}{r})$, then
$\gamma|_{[s',r']}\in\RAdm{s'}{r'}$ (respectively,
$\gamma|_{[s',r']}\in\RAdm[C]{s'}{r'}$).
Now, for $\gamma\in\RAdm{s}{r}$, set
\begin{equation}\label{eq:defnofIr}
I^r[s,r,\gamma,f,\sigma]
\triangleq
\int\limits_{[s,r]\times\td\times\rd}\Big[f(t,x)-L(x,v)\Big]
\gamma(dt,dx,dv)
+\int\limits_{\td}\sigma(x)\mathfrak{f}_r[\gamma](dx);
\end{equation}
compare with \eqref{eq:defnofIa}.
Using the ``duality'' notation from the Preliminaries,
we may very tersely write this as
$$
I^r[s,r,\gamma,f,\sigma] = \action{\gamma}{f-L} + \action{\mathfrak{f}_r[\gamma]}{\sigma}.
$$
\subsubsection{Equivalence of
maximizing $I^r$ and $I^a$}\label{subsubsection:equivalenceIaIr}
We will now see that the value functions
\eqref{prel:intro:value} can be attained by
maximizing over the functional \(I^r\) on relaxed
controls that we just defined in \eqref{eq:defnofIr}:
see formula \eqref{prel:equality:value_RAdm} below.
Indeed, first observe that if \((m(\cdot),w)\in\Adm{s}{r}\)
with \(m(s) = m_*\) and
\(\gamma\) is constructed as in
Proposition \ref{prel:prop:link_adm}, then
\(\mathfrak{f}_t[\gamma] = m_t\) for a.~e.~\(t\) and
\(\mathfrak{v}[\gamma] = w\). Thus,
for every \((m(\cdot),w)\in\Adm{s}{r}\)
and \(\gamma\) constructed in this way, we have
\[
I^r[s,r,\gamma,f,\sigma]
= I^a[s,r,m(\cdot),w,f,\sigma].
\]
Therefore, if $C>0$ and
\(\operatorname{Val}[s,r,m_*,f,\sigma]\)
is defined as a maximization over
\(\Adm[C]{s}{r}\), then
\[
\operatorname{Val}[s,r,m_*,f,\sigma]
\leq \sup\{I^r[s,r,\gamma,f,\sigma]:\, \gamma\in \RAdm[C]{s}{r},
\mathfrak{f}_s[\gamma] = m_*\}.
\]
In the opposite direction, let \(C>0\) and
\(\gamma\in \RAdm[C]{s}{r}\) such that
\(\mathfrak{f}_{s}[\gamma]=m_*\), and let
\(m(t)\triangleq \mathfrak{f}_t[\gamma], s\leq t\leq r\),
\(w\triangleq \mathfrak{v}[\gamma]\). Then
\((m_{\cdot},w)\in \Adm[C]{s}{r}\) and
\begin{align*}
	I^a&[s,r,m(\cdot),w,f,\sigma]
	=
	\int_s^r\int_{\dtorus}
	\int_{\R^d} [f(t,x)-L(x,y)]\delta_{w(t,x)}(dy)m_t(dx)dt
	+ \int_{\dtorus}\sigma(x)\mathfrak{f}_r[\gamma](dx)
	\\
	= & \ \int_s^r
	\int_{\dtorus}
	\big[f(t,x)-L\big(x,\int_{\mathbb{B}_C}v\gamma(dv|t,x)\big)\big]m_t(dx) dt
	+ \int_{\dtorus} \sigma(x)m_r(dx).
\end{align*}
Since \(\gamma(\mathbb{B}_C|t,x)=1\), Jensen's inequality gives us
\begin{align*}
I^a[s,&r,m(\cdot),w,f,\sigma]
\\\geq & \
\int_s^r \int_{\dtorus}\int_{\mathbb{B}_C}[f(t,x)-L(x,v)]\gamma(dv|t,x)
\mathfrak{f}_t[\gamma](dx)dt  +\int_{\dtorus}\sigma(x)\mathfrak{f}_r
[\gamma](dx)
\\
= & \ \int_s^r\int_{\dtorus}
\int_{\mathbb{B}_C}[f(t,x)-L(x,v)]\gamma(dv|t,x)(\operatorname{p}^{1,2}\sharp\gamma)(dx|t) dt
+\int_{\dtorus}\sigma(x)\mathfrak{f}_r
[\gamma](dx)
\\
= & \int_{\tripleprod}[f(t,x)-L(x,v)]\gamma(dt,dx,dv)
+\int_{\dtorus}\sigma(x)\mathfrak{f}_r
[\gamma](dx)
\\
= & \ I^r[s,r,\gamma,f,\sigma].
\end{align*}
In conclusion, we have:
\begin{lemma}\label{lemma:equality:value_RAdm}
	Let \(C>0\), \(0\leq s<r\leq T\), \(m_{\ast}\in \wassonespace\)
    and consider \(\operatorname{Val}[s,r,m_{\ast},f,\sigma]\) of \eqref{prel:intro:value}
    to be a maximization over $\Adm[C]{s}{r}$. Then the equality
\begin{equation}\label{prel:equality:value_RAdm}
	\operatorname{Val}[s,r,m_*,f,\sigma]
	=\sup\Big\{I^r[s,r,\gamma,f,\sigma]:\,
	\gamma\in \RAdm[C]{s}{r},\, \mathfrak{f}_s[\gamma]=m_*\Big\}
\end{equation}
holds.
\end{lemma}
\subsubsection{The topology on \(\RAdm[C]{s}{r}\)}\label{subsubsection:topologyRAdm}
For brevity let us denote \(E=C(\tripleprod)\)
(note that the space \(E\), like the sets \(\RAdm[C]{s}{r}\)
themselves, depends on the constant \(C\), which is considered
fixed throughout this discussion).
For each \(s,r\) with \(0\leq s<r\leq T\),
every element \(\gamma\in \RAdm[C]{s}{r}\) is put
in correspondence with
an element \(T_{\gamma}\) of the topological dual
\(E^{\ast}\) through the formula
\[
\langle T_{\gamma},f\rangle=\int_{\tripleprod}f(t,x,v)\gamma(dt,dx,dv),
\]
\(f\in E\) (the function $f$ in this short paragraph is not the same
as above, it is just an arbitrary function in the space \(E\)). The normed space
\(E^{\ast}\) (with the dual norm \(\|T\|_{E^{\ast}}=\sup_{\|f\|\leq 1}\langle T,f\rangle\))
is endowed with the
weak-\(\star\) topology \(\mathcal{T}\triangleq\sigma(E^{\ast},E)\), the topology
generated by the functionals
\(E^{\ast}\ni T \mapsto \langle T, f\rangle\), \(f\in E\).
Thus each \(\RAdm[C]{s}{r}\) is identified with a convex subset of the
topological space \(E^{\ast}\).
Since \(E\) is separable, any \(\|\cdot\|_{E^{\ast}}\)-bounded subset of $E^{\ast}$
in the weak-\(\star\) topology
is metrizable \cite[Theorem 5.1]{conway1990course}, so
closedness and compactness in any \(\|\cdot\|_{E^{\ast}}\)-bounded subset of
\((E^{\ast},\mathcal{T})\) are equivalent to sequential
closedness and sequential compactness.
All discussions regarding
\(\RAdm[C]{s}{r}\) are done in the spirit of this identification of
\(\RAdm[C]{s}{r}\) with a bounded convex subset of
\((E^{\ast},\mathcal{T})\). We will call convergence in this topology
also \textit{narrow}.

The set \(\RAdm[C]{s}{r}\) is closed: let \(\{\gamma_n\}_1^{\infty}\)
be a sequence of \(\RAdm[C]{s}{r}\) which converges to
\(T\in E^{\ast}\). By Riesz' representation theorem
\cite{dunford1958linear1}, there is a Borel signed measure
\(\gamma\) on \(\tripleprod\) representing \(T\), that is,
\(\int f(t,x,v)\gamma(dt,dx,dv) = T(f)\) for all
\(f\in C(\tripleprod)\).
Thus
\begin{align*}
  \int f(t,x,v)\gamma_n(dt,dx,dv) \to \int f(t,x,v) \gamma(dt,dx,dv)
\end{align*}
for all \(f\in C(\tripleprod)\). From this it follows that
\(\gamma\) is a measure (non-negative). Properties 1-4 of
the definition of
\(\RAdm[C]{s}{r}\) are easily checked. This shows that \(\RAdm[C]{s}{r}\)
is closed.
Furthermore, since for \(\gamma\in \RAdm[C]{s}{r}\),
\(|\langle T_{\gamma}, f\rangle| \leq (r-s)\|f\|_E\) for all
\(f\in E\) we have that
\begin{align*}
  \RAdm[C]{s}{r} \subset B_{s,r}^{\ast}
  \triangleq \{ T\in E^{\ast} \ : \
  \|T\|_{E^{\ast}} \leq r-s\}.
\end{align*}
By Alaoglu's theorem,
\(B_{s,r}^{\ast}\) is a compact set with respect to
\(\mathcal{T}\). Thus, \(\RAdm[C]{s}{r}\) is a compact, convex
set with the topology \(\mathcal{T}\).

We remark that, despite
the inclusion \(\RAdm[C]{s}{r}\subset B_{s,r}^{\ast}\),
if \(s'\leq s < r \leq r'\), and \(\gamma\in
\RAdm[C]{s}{r}\), then \(\gamma\notin \RAdm[C]{s'}{r'}\) unless
\(s=s'\) and \(r=r'\). However, both spaces
\(\RAdm[C]{s}{r}\) and \(\RAdm[C]{s'}{r'}\) are included
in the larger of the two dual spaces \(B_{s',r'}^{\ast}\),
and certainly in \(B_{0,T}^{\ast}\) too.
This allows us to consider sequences of relaxed controls
that are taken from \(\RAdm{s}{r}\) with different
values of \(s\) and \(r\).
\subsection{Stochastic control
formulation}\label{subsection:stochasticcontrolformulation}
We briefly describe a third way of formulating the same problem,
by means of stochastic controls.
A \emph{stochastic control process} on \([s,r]\)
is a tuple
\((\Omega,\mathcal{F},(\mathcal{F}_t)_{s\leq t\leq r},\mathbb{P},W,\xi,\alpha)\)
consisting of
\begin{itemize}[-~~]
  \item a filtered probability space
    \((\Omega,\mathcal{F},(\mathcal{F}_t)_{s\leq t\leq r},\mathbb{P})\)
    carrying a \(d\)-dimensional \((\mathcal{F}_t)_{s\leq t\leq r}\)-Wiener process
	\((W_t)_{t\geq s}\) (in the sense of Section \ref{subsection:feedbackformulation});
	\item an \(\mathcal{F}_s\)-measurable random variable \(\xi\);
	\item a bounded, \((\mathcal{F}_t)_{s\leq t\leq r}\)-progressively measurable process \(\alpha\),
\end{itemize}
and we let \(X\) denote the solution of the stochastic differential equation
\begin{equation}\label{prel:eq:X__sde_stochastic}
  X_t = \xi + \int_{s}^{t} \alpha_u\,du + \sqrt{2}\,(W_t-W_s), \qquad s\leq t\leq r.
\end{equation}
We will often refer simply to the pair \((X,\alpha)\) as the stochastic
control process, the probabilistic setup being implicit.
The corresponding payoff functional is
\begin{equation*}
  I^s[s,r,X,\alpha,f,\sigma]\triangleq
  \mathbb{E}\Bigg[\int_{s}^{r}\bigl(f(t,X_t)-L(X_t,\alpha_t)\bigr)dt
  +\sigma(X_r)\Bigg].
\end{equation*}
If $w$ is a given feedback control and $X$ satisfies \eqref{prel:aux_prob:dynamics}
with this \(w\),
then $X$ is a solution of \eqref{prel:eq:X__sde_stochastic} with
\[
\alpha_t = w(t,X_t).
\]
That is, every value of the target functional \(I^s\) by a
feedback control (the context of the formulation
\eqref{prel:aux_prob:criterion}) can be attained
by a stochastic control process.
Thus,
\begin{equation}\label{prel:ineq:value_alpha}
\begin{split}
	& \operatorname{Val}[s,r,m_*,f,\sigma]
	\\
    \leq & \
	\sup\Big\{I^s[s,r,X,\alpha,f,\sigma]:
	(X,\alpha)\text{ is a stoch.~control process},
	\, \operatorname{Law}(X(s))=m_*\Big\},
\end{split}
\end{equation}
with the same understanding as before,
that if the left-hand side is obtained
by maximization over feedback controls
bounded by a constant $C$, then the maximization
on the right is done over pairs $(X,\alpha)$
with $|\alpha|\leq C$.
Conversely, let \((X,\alpha)\) be a stochastic control process
with \(\operatorname{Law}(X_s)=m_*\)
and \(|\alpha_t|\leq C\). Define the measure \(\gamma\)
on \([0,T]\times\dtorus\times\mathbb{B}_C\)
by
the rule: for \(\psi\in\mathcal{C}_b([s,r]\times\td\times\rd)\),
\[
\action{\gamma}{\psi}\triangleq
\EE\int_s^r\psi(t,X_t,\alpha_t)\,dt .
\]
Then \(\gamma\in\RAdm[C]{s}{r}\). Indeed, properties 1--3 of
Definition \ref{prel:def:relaxed_control_procs} are immediate, while
for \(\varphi\in\mathcal{C}^{1,2}_c((s,r)\times\dtorus)\), It\^o's
formula applied to \(\varphi(t,X_t)\) gives, after taking expectations,
\[
\EE\int_s^r\big(\partial_t\varphi
+D\varphi\cdot\alpha_t+\Delta\varphi\big)(t,X_t)\,dt=0,
\]
which is precisely \eqref{eq:weakcontrolFPK} for \(\gamma\).
Moreover, disintegrating \(\gamma\) in time shows that
\(\mathfrak{f}_t[\gamma]=\operatorname{Law}(X_t)\) for a.e.~\(t\), and
since both flows are continuous (for the left-hand side, this is the
convention following Proposition \ref{prel:prop:Holder_m}; for the
right-hand side, \(\EE|X_t-X_{t'}|\leq C|t-t'|+\sqrt{2}\,\EE|W_t-W_{t'}|\)),
the equality holds for \emph{every} \(t\); in particular
\(\mathfrak{f}_s[\gamma]=m_*\) and
\(\mathfrak{f}_r[\gamma]=\operatorname{Law}(X_r)\).
Applying the definitions,
\(I^s[s,r,X,\alpha,f,\sigma]=I^r[s,r,\gamma,f,\sigma]\),
and therefore, by Lemma \ref{lemma:equality:value_RAdm},
\[
I^s[s,r,X,\alpha,f,\sigma]\ \leq\
\operatorname{Val}[s,r,m_*,f,\sigma]
\]
(the value being understood as a maximization over
\(\Adm[C]{s}{r}\)). So, every realization of the target functional
using stochastic controls is dominated by the value over control
processes, and the converse of the inequality
\eqref{prel:ineq:value_alpha} holds as well. In conclusion,
	\begin{equation}\label{prel:equality:value_alpha}
    \begin{split}
		& \operatorname{Val}[s,r,m_*,f,\sigma]
		\\
        = & \
		\sup\Big\{I^s[s,r,X,\alpha,f,\sigma]:
		(X,\alpha)\text{ is a stoch.~control process},
		\, \operatorname{Law}(X(s))=m_*\Big\}
        \end{split}
	\end{equation}
	obtains, where we keep the caveat that the same boundedness constant $C$
is used for control processes, and stochastic controls,
on left and right of the equality.
We also remark, although we will not use this, that
\(\mathfrak{v}[\gamma](t,x)=\EE[\alpha_t\mid X_t=x]\) for
\(\mathfrak{f}_t[\gamma]\)-a.e.~\(x\) and a.e.~\(t\)
(this is the classical way to pass from
a stochastic control to a control process,
by the conditional-expectation feedback).

\subsection{Regularity of the target functionals and value
functions}\label{subsection:regularityofvaluefunctions}

The next proposition concerns the
convergence of the functional $I^r$ along
subsequences in all its arguments. Recall that,
by the comments at the end of Section \ref{subsubsection:topologyRAdm},
if the \(\gamma_n\) relate to subintervals
\([s_n,r_n]\subset [s',r']\), they are to be regarded
as elements of the common space
\(B_{s',r'}^{*}\), itself a subspace
of $E^*$ (the topological dual of
\(C([0,T]\times\dtorus\times\mathbb{B}_C)\),
with the topology $\mathcal{T}=\sigma(E,E^*)$).
\begin{proposition}\label{prel:prop:cont_I_r} Let
\begin{enumerate}[label=(\roman*)]
\item the sequences
$\{s_n\}_{n=1}^\infty,\, \{r_n\}_{n=1}^\infty\subset [s',r']$
(\(s_n <r_n\)) converge to $s$ and $r$ respectively (\(s<r\));
\item the functions $f_n:[s_n,r_n]\times \td\rightarrow \mathbb{R}$,
$\sigma_n:\dtorus\to\R$ be continuous
and converge uniformly to $f$ and $\sigma$ respectively;
\item
the relaxed controls $\gamma_n\in\RAdm[C]{s_n}{r_n}$,
regarded as elements of \(B^{\ast}_{s',r'}\),
converge in the topology \(\mathcal{T}\)
to a relaxed control process $\gamma\in\RAdm[C]{s}{r}$.
\end{enumerate} Then,
\[I^r[s_n,r_n,\gamma_n,f_n,\sigma_n]\rightarrow
I^r[s,r,\gamma,f,\sigma]\text{ as }n\rightarrow\infty.\]
\end{proposition}
The convergence of the functions \(f_n\) to \(f\) should be
understood in the sense that extensions of them to the full
interval \([s',r']\) (defining them
to be, say, \(f_n(s_n)\) on \([s',s_n]\) and
\(f_n(r_n)\) on \([r_n,r']\)) converge
to \(f\).
The proof follows from the definition of the functional
and Proposition \ref{prel:prop:continuity}; the convergence of the
terminal terms
\(\action{\mathfrak{f}_{r_n}[\gamma_n]}{\sigma_n}\), where the
evaluation times vary, additionally uses the uniform
\(\tfrac12\)-H\"older bound of Proposition \ref{prel:prop:Holder_m}.
\begin{proposition}\label{prel:prop:val_lip}
	Assume that
	\begin{enumerate}
		\item $f$ is Lipschitz continuous w.r.t.~$x$
		with a Lipschitz constant $K_f$ that does not depend on $t$;
		\item $\sigma$ is Lipschitz continuous w.r.t.~$x$
		with a Lipschitz constant $K_\sigma$.
	\end{enumerate} Then, for each $y',y''\in\td$,
\begin{equation}\label{prel:ineq:Lipschitz_values}
	\big|\operatorname{Val}[s,r,y',f,\sigma]
	-\operatorname{Val}[s,r,y'',f,\sigma]\big|
	\leq
	(K_\sigma+(r-s)(K_f+K_L))|y'-y''|.
\end{equation}
\end{proposition}
This follows directly from the
representation of
the function $\operatorname{Val}$
through the stochastic control processes
\eqref{prel:equality:value_alpha}: if \((X',\alpha)\) is a stochastic
control process with \(X'_s=y'\), let \(X''\) be driven by the
same \(\alpha\) and \(W\) from \(y''\); then
\(X''_t-X'_t\equiv y''-y'\), and comparing
\(I^s\) termwise gives
\(|I^s[s,r,X'',\alpha,f,\sigma]-I^s[s,r,X',\alpha,f,\sigma]|
\leq(K_\sigma+(r-s)(K_f+K_L))|y'-y''|\)
(the Lagrangian terms via \ref{prel:cond:L_lipschitz});
taking suprema yields \eqref{prel:ineq:Lipschitz_values}.
The Lipschitz constant for the value function
\(\operatorname{Val}[s,r,\cdot,f,\sigma]\) that appears in
\eqref{prel:ineq:Lipschitz_values} will be denoted by
	\[
	\widehat{K}_{\sigma,f}
	\triangleq K_\sigma+(r-s)(K_f+K_L).
	\]
\begin{lemma}\label{prel:lm:verification}
Let \(0\leq s<r\leq T\) and assume that
\begin{enumerate}[label=(V\arabic*)]
\item\label{cond:V:f} \(f:[s,r]\times\dtorus\to\R\) is Lipschitz
continuous w.r.t.~\(x\) with a constant \(K_f\) independent of \(t\),
and \(\tfrac12\)-H\"older continuous w.r.t.~\(t\), uniformly in \(x\);
\item\label{cond:V:sigma} \(\sigma\in\operatorname{Lip}_{K_\sigma}(\td)\)
and \(\sigma\in\mathcal{C}^{2+\alpha'}(\td)\) for some
\(\alpha'\in(0,1)\).
\end{enumerate}
Then the following statements hold true.
\begin{enumerate}[label=(\alph*)]
\item\label{lm:verification:a} The terminal-value problem
\begin{equation}\label{eq:HJBfsigma}
\partial_t\varphi+\triangle\varphi+H(x,D\varphi)+f(t,x)=0
\ \text{ on }(s,r)\times\dtorus,
\qquad \varphi(r,\cdot)=\sigma,
\end{equation}
admits a solution \(\varphi\) with the regularity listed in Section
\ref{subsection:notionofclassicalsolutionofmfgsystem}
(\(\varphi\in\mathcal{C}^0([s,r]\times\dtorus)\), continuously
differentiable in \(t\) on the open cylinder, with \(D\varphi\),
\(D_{xx}\varphi\) continuous on the closed cylinder; hence, by Remark
\ref{rem:dtextension},
\(\varphi\in\mathcal{C}^{1,2}([s,r]\times\dtorus)\)).
This solution is unique, and
\[
\varphi(t,y)=\operatorname{Val}[t,r,y,f,\sigma]
\qquad\text{for all }(t,y)\in[s,r]\times\dtorus;
\]
in particular \(|D\varphi|\leq\widehat{K}_{\sigma,f}\)
on \([s,r]\times\dtorus\), by Proposition \ref{prel:prop:val_lip}.
\item\label{lm:verification:b} The feedback
\[
w^*(t,x)\triangleq D_pH\big(x,D\varphi(t,x)\big)
\]
is continuous on \([s,r]\times\dtorus\) and
\(|w^*|\leq C_*[\widehat{K}_{\sigma,f}]\), with \(C_*[\cdot]\) the
function of Lemma \ref{prel:lm:Hamiltonian}(ii). For every
\(t'\in[s,r)\) and every \(m_0\in\ptd\),
\begin{equation}\label{eq:Valismphi}
\operatorname{Val}[t',r,m_0,f,\sigma]=\action{m_0}{\varphi(t',\cdot)},
\end{equation}.
and if \((m_\cdot,w)\in\Adm{t'}{r}\), \(m_{t'}=m_0\), is such that
\(w(t,x)=w^*(t,x)\) for a.e.~\(t\in[t',r]\) and
\(m_t\)-a.e.~\(x\in\dtorus\), then \((m_\cdot,w)\) is a maximizer:
\(I^a[t',r,m_\cdot,w,f,\sigma]=\operatorname{Val}[t',r,m_0,f,\sigma]\).
At least one such pair exists, namely the law flow of the SDE
\eqref{prel:aux_prob:dynamics} driven by \(w^*\) from \(m_0\)
starting at time $t'$.
\item\label{lm:verification:c} Conversely,
for every
\(t'\in[s,r)\) and every \(m_0\in\ptd\),
if
\((\hat{m}_\cdot,\hat{w})\in\Adm{t'}{r}\) with
\(\hat{m}_{t'}=m_0\) is a
maximizer of \(I^a[t',r, \cdot , \cdot  ,f,\sigma]\)
(with initial values constrained at $m_0$), then
\begin{equation}\label{eq:equalityatmaximizervwiththat}
\max_{v\in\R^d}
\big\{
D\varphi(t,x)\cdot v-L(x,v)
\big\}
=
D\varphi(t,x)\cdot\hat{w}(t,x)-L\big(x,\hat{w}(t,x)\big),
\end{equation}
i.e.~\(\hat{w}(t,x)=w^*(t,x)\) (with $w^*$
as defined in part \ref{lm:verification:b})
for a.e.~\(t\in[t',r]\) and
\(\hat{m}_t\)-a.e.~\(x\in\dtorus\). Moreover:
\begin{equation}\label{eq:Lambdazero}
\action{\hat{m}_a}{\varphi(a,\cdot)}
=I^a[a,r,\hat{m}_\cdot,\hat{w},f,\sigma]
\qquad\text{for every }a\in[t',r].
\end{equation}
\end{enumerate}
\end{lemma}
\begin{proof}
Let us dispatch the first claim of part (a).
Under \ref{cond:V:f}--\ref{cond:V:sigma} and
\ref{prel:cond:L_lipschitz}--\ref{prel:cond:L_coercitive},
the terminal-value problem \eqref{eq:HJBfsigma} admits a solution
with the stated regularity: see
\cite[Th.~V.6.1]{ladyzenskaja1967linear}; this solvability, for data
of precisely our class, is the first step of the fixed-point
existence argument for the MFG system in \cite{fourhorses}.
By Remark \ref{rem:dtextension} (applied to \eqref{eq:HJBfsigma},
where the continuity in \(t\) of the zeroth-order term is part of
\ref{cond:V:f}), \(\varphi\in\mathcal{C}^{1,2}([s,r]\times\dtorus)\),
and \(D\varphi\) is bounded on \([s,r]\times\dtorus\), being
continuous on a compact set. We will prove the other claims in
the three items across several steps.

Let \(t'\in[s,r)\),
\(m_0\in\ptd\), and let \((m_\cdot,w)\in\Adm{t'}{r}\) with
\(m_{t'}=m_0\).

\emph{Step 1}. The superposition principle (see the comment following
\eqref{eq:defnofIa}) provides a filtered probability space carrying a
solution \((X_t)_{t'\leq t\leq r}\) of
\(dX_t=w(t,X_t)\,dt+\sqrt{2}\,dW_t\) with
\(\operatorname{Law}(X_t)=m_t\), \(t'\leq t\leq r\). Applying It\^o's
formula to \(\varphi(t,X_t)\) on \([t',r]\) (legitimate by the
regularity established above), taking expectations
and using
Fubini's theorem,
\begin{equation}\label{eq:itoidentity}
\action{m_r}{\sigma}-\action{m_0}{\varphi(t',\cdot)}
=\int_{t'}^r\int_{\dtorus}
\big[\partial_t\varphi+D\varphi\cdot w+\triangle\varphi\big](t,x)\,
m_t(dx)\,dt .
\end{equation}
By \eqref{eq:HJBfsigma}, the integrand equals
\(D\varphi(t,x)\cdot w(t,x)-H(x,D\varphi(t,x))-f(t,x)\), which, by
\eqref{eq:Young}, is bounded above by \(L(x,w(t,x))-f(t,x)\), with
equality at \((t,x)\) if and only if \(w(t,x)=w^*(t,x)\).
Adding \(\int_{t'}^r\action{m_t}{f(t,\cdot)-L(\cdot,w(t,\cdot))}dt\)
to both sides of \eqref{eq:itoidentity} we conclude:
\begin{equation}\label{eq:verificationineq}
I^a[t',r,m_\cdot,w,f,\sigma]\ \leq\ \action{m_0}{\varphi(t',\cdot)},
\qquad\text{with equality if }w=w^*\ m_t\otimes dt\text{-a.e.}
\end{equation}

\emph{Step 2}.
The feedback \(w^*\) is Borel measurable and bounded: indeed,
\(w^*(t,x)=v^*(x,D\varphi(t,x))\) is continuous on
\([s,r]\times\dtorus\) by Lemma \ref{prel:lm:Hamiltonian}(iii) and the
continuity of \(D\varphi\), and bounded by Lemma
\ref{prel:lm:Hamiltonian}(ii) with \(R=\max|D\varphi|\).
By \cite{Veretennikov1981}, the SDE \eqref{prel:aux_prob:dynamics}
with \(w=w^*\) and \(\operatorname{Law}(\xi)=m_0\) has a solution on
\([t',r]\); its law flow \(m^*_\cdot\) solves \eqref{prel:eq:FPK}, so
\((m^*_\cdot,w^*)\in\Adm{t'}{r}\), and the equality case of
\eqref{eq:verificationineq} gives
\(I^a[t',r,m^*_\cdot,w^*,f,\sigma]=\action{m_0}{\varphi(t',\cdot)}\).
Together with the inequality of \eqref{eq:verificationineq} for
arbitrary \((m_\cdot,w)\), this proves \eqref{eq:Valismphi} and the
maximizer property claimed in \ref{lm:verification:b}. In particular,
taking \(m_0=\delta_y\),
\(\varphi(t',y)=\operatorname{Val}[t',r,y,f,\sigma]\) for every
\((t',y)\); Proposition \ref{prel:prop:val_lip} then gives
\(\varphi(t,\cdot)\in\operatorname{Lip}_{\widehat{K}_{\sigma,f}}(\td)\)
for every \(t\), whence \(|D\varphi|\leq\widehat{K}_{\sigma,f}\) and,
by Lemma \ref{prel:lm:Hamiltonian}(ii),
\(|w^*|\leq C_*[\widehat{K}_{\sigma,f}]\). This completes the proof of
\ref{lm:verification:a} and \ref{lm:verification:b}, except that only
the uniqueness claim in \ref{lm:verification:a} remains. But
if \(\varphi'\) is another
solution with the stated regularity, then Steps 1--2 apply to it
word by word,
so \(\varphi'(t,y)=\operatorname{Val}[t,r,y,f,\sigma]
=\varphi(t,y)\).

\emph{Step 3}. Let
\((\hat{m}_\cdot,\hat{w})\in\Adm{t'}{r}\), \(\hat{m}_{t'}=m_0\), be a
maximizer as declared in the statement of \ref{lm:verification:c},
and set
\[
\Lambda(\theta)\ \triangleq\
\action{\hat{m}_{\theta}}{\varphi({\theta},\cdot)}
-I^a\big[{\theta},r,\hat{m}_\cdot,\hat{w},f,\sigma\big],
\qquad {\theta}\in[t',r].
\]
Essentially the same
computation of Step 1 (before applying the Legendre inequality),
but carried out with $(\hat{m},\hat{w})$
on \([{\theta},r]\) instead of
\([t',r]\), gives
\begin{equation}\label{eq:Lambdaidentity}
\Lambda({\theta})=\int_{\theta}^r\int_{\dtorus}
\Big[H\big(x,D\varphi(t,x)\big)
-\big(D\varphi(t,x)\cdot\hat{w}(t,x)
-L(x,\hat{w}(t,x))\big)\Big]\,\hat{m}_t(dx)\,dt ,
\end{equation}
whose integrand is non-negative by \eqref{eq:Young}. On the other
hand,
\[
\Lambda(t')=\action{m_0}{\varphi(t',\cdot)}
-\operatorname{Val}[t',r,m_0,f,\sigma]=0
\]
by the maximizer hypothesis and \eqref{eq:Valismphi}. Since
\(0\leq\Lambda({\theta})\leq\Lambda(t')=0\) for every \({\theta}\in[t',r]\) (the
integrand in \eqref{eq:Lambdaidentity} being non-negative and the
domain of integration decreasing in \({\theta}\)), the integrand vanishes for
a.e.~\(t\) and \(\hat{m}_t\)-a.e.~\(x\), which is
\eqref{eq:equalityatmaximizervwiththat}; the identification
\(\hat{w}=w^*\) follows from the uniqueness of the maximizer in Lemma
\ref{prel:lm:Hamiltonian}(i).
The very last formula, \eqref{eq:Lambdazero}, is precisely \(\Lambda\equiv0\),
which we just obtained.
\end{proof}
The existence of maximizers in \eqref{prel:intro:value}, for data as
in Lemma \ref{prel:lm:verification}, is part of statement
\ref{lm:verification:b} of that lemma. The next proposition
provides a bound for all maximizers.
\begin{proposition}\label{prop:prel:boundness}
	Let $f$ and $\sigma$ satisfy the assumptions
	\ref{cond:V:f}--\ref{cond:V:sigma} of Lemma
	\ref{prel:lm:verification}.
	Set
	\[
	C\triangleq C_*\big[\widehat{K}_{\sigma,f}\big] =C_*\big[K_\sigma+(r-s)(K_f+K_L)\big],
	\]
	with \(C_*[\cdot]\) the function of Lemma \ref{prel:lm:Hamiltonian}(ii); thus \(C\) depends only on
	$K_{\sigma}$, $K_f$, $K_L$, $C_L^0$, $C^1_L$, $q$, $M^0_L$ and
	$r-s$. Then for each $m_*\in\mathcal{P}^1(\td)$ and
	control process $(\hat{m}(\cdot),\hat{w})\in\Adm{s}{r}$
	fulfilling the conditions
	\begin{itemize}[-~~]
		\item $\hat m_s=m_*$,
		\item $\operatorname{Val}[s,r,m_*,f,\sigma]
		=
		I^a[s,r,\hat{m}(\cdot),\hat{w},f,\sigma]$,
	\end{itemize}
	one has
	\[
	|\hat{w}(t,x)|\leq C\text{ for a.e. }
	s\leq t\leq r\text{ and }\hat{m}_t \text{-a.e. }x\in\td.
	\]
\end{proposition}
\begin{proof}
By Lemma part \ref{lm:verification:c} of Lemma
\ref{prel:lm:verification},
\(\hat{w}(t,x)=w^*(t,x)\) for a.e.~\(t\) and
\(\hat{m}_t\)-a.e.~\(x\), and by part
\ref{lm:verification:b} of the same lemma,
\(|w^*|\leq C_*[\widehat{K}_{\sigma,f}]=C\).
\end{proof}

\begin{proposition}\label{prel:prop:cont_Val} Assume that
	\begin{itemize}[-~~]
		\item the sequences
		$\{s_n\}_{n=1}^\infty,\, \{r_n\}_{n=1}^\infty
		\subset [s',r']$
		converge to $s$ and $r$
		respectively ($s_n<r_n$, $s<r$);
		\item the functions $f_n:[s',r']\times\td\rightarrow \R$
		are uniformly Lipschitz continuous w.r.t.~$x$ with a constant
		$K_f$, and each \(f_n\) is \(\tfrac12\)-H\"older continuous
		w.r.t.~\(t\);
		\item the functions
		$\sigma_n:\td\rightarrow \R$
		are uniformly Lipschitz
		continuous with a constant $K_\sigma$, and each \(\sigma_n\)
		belongs to \(\mathcal{C}^{2+\alpha_n}(\td)\) for some
		\(\alpha_n\in(0,1)\);
		\item the sequences of
		functions
		$\{f_n\}_{n=1}^\infty$
		and $\{\sigma_n\}_{n=1}^\infty$
		converge uniformly to functions
		$f:[s,r]\times\dtorus\to\R$
		and $\sigma:\dtorus\to\R$ respectively.
	\end{itemize} Then, the sequence of functions
$\{\operatorname{Val}[s_n,r_n,\cdot,f_n,\sigma_n]\}_{n=1}^\infty$
converges to $\operatorname{Val}[s,r,\cdot,f,\sigma]$, uniformly on
\(\dtorus\).
\end{proposition}
\begin{proof}
	It suffices to use the representation of the value function
	through the stochastic control
    processes \eqref{prel:equality:value_alpha},
    take into account the fact that,
    by Lemma \ref{prel:lm:verification}\ref{lm:verification:b} and
    Proposition \ref{prop:prel:boundness}, for each \(n\) the optimal
    control $\alpha$ can be
	chosen to satisfy the inequality
	\[
    |\alpha_t|\leq C,
    \]
    where $C$
	is the constant of
	Proposition \ref{prop:prel:boundness} for the pair
	\((K_f,K_\sigma)\), which is uniform in \(n\), and use the
	continuity
    of the functional $I^s$ in the functions $f$, $\sigma$ and
    the time endpoints. The uniformity of the convergence on
    \(\dtorus\) follows from the pointwise convergence and the
    equi-Lipschitz continuity of the family
    \(\{\operatorname{Val}[s_n,r_n,\cdot,f_n,\sigma_n]\}\)
    guaranteed by Proposition \ref{prel:prop:val_lip}.
\end{proof}

\section{Mean field game solutions and value multifunction}\label{section:MFGsolsandvaluemf}
Below we consider the mean field game on $[t_0,T]$
with the dynamics of each player given by
\eqref{intro:sys:dynamics}
and his/her payoff equal to
\eqref{intro:sys:payoff}.
The initial distribution of players is given by $m_*$.

From now on,
we assume that $C$ is the constant given by Proposition
\ref{prop:prel:boundness} with $K_f=K_F$, $K_\sigma=K_G$
and $r-s=T$.
We note that the data that will occupy the places of \(f\) and
\(\sigma\) from now on, namely \(f(t,x)=F(x,\hat{m}_t)\) and
\(\sigma=G(\cdot,\hat{m}_T)\) for \(\tfrac12\)-H\"older continuous
paths \(\hat{m}_\cdot\), satisfies the assumptions
\ref{cond:V:f}--\ref{cond:V:sigma} of Lemma
\ref{prel:lm:verification}, by \ref{sol_mfg:cond:F},
\ref{sol_mfg:cond:G} and \ref{sol_mfg:cond:G2plusalpha}.
We additionally set
\begin{equation}\label{sol_mfg:intro:M_L}
	L_C\triangleq \sup\big\{|L(x,v)|:\, x\in\td,\,
	v\in\mathbb{B}_C \big\}.
\end{equation}
For an arbitrary $\sigma\in \mathcal{C}(\dtorus)$, we let
 $A_{\hat{m}_\cdot}^{s,r}\sigma$ be the
 continuous on $\dtorus$ function
\begin{equation}\label{sol_mfg:intro:Lax_Oleink}
(A_{\hat{m}_\cdot}^{s,r}\sigma)(y )\triangleq \operatorname{Val}[s,r,y,F(\cdot,\hat{m}_{\cdot}),\sigma],
\end{equation}
which is to say,
\((A_{\hat{m}}^{s,r}\sigma)(y)\)
is \(\varphi(s,y)\)
where
\(\varphi(s,y)
=\max\{I^a[s,r,m_{\cdot},w,F(\cdot,\hat{m}),\sigma]\}\)
over all $(m_{\cdot},w)\in \Adm[C]{s}{r}$ with
$m_s=\delta_y$;
equivalently,
due to Lemma \ref{lemma:equality:value_RAdm},
\begin{align*}
(A_{\hat{m}}^{s,r}\sigma)(y)   = & \
	\max\Bigg\{ \action{\mathfrak{f}_r[\gamma]}{\sigma} +
	\int_{\tripleprod}
	\big( F(x,\hat{m}(t)) - L(x,v)\big) \gamma(dt,dx,dv)
	\\
	& \qquad \qquad \qquad \qquad   :
	\gamma \in \RAdm[C]{s}{r}, \ \mathfrak{f}_s[\gamma]=\delta_y \Bigg\}
\end{align*}
(when \(\sigma\) is in the class \ref{cond:V:sigma}, the maxima are
attained, by Lemma \ref{prel:lm:verification}\ref{lm:verification:b};
for general continuous \(\sigma\) the right-hand sides are understood
as suprema).
Thus, $(A^{s,r}_{\hat{m}}\sigma)(y)$ is the value function
of the control problem (in any of the three above versions)
on $[s,r]$ with running cost $L$, running interaction term $F(\cdot,\hat{m}_t)$
and terminal cost $\sigma$ for the initial condition $y\in\dtorus$.

Recall now the definition of a solution of the MFG system
of equations \eqref{intro:eq:mfg_system};
this was introduced in Section \ref{subsection:notionofclassicalsolutionofmfgsystem}.
The hypotheses \textbf{(HS)} are standing.
\begin{proposition}\label{sol_mfg:def:solution} Let $t_0\in [0,T)$,
	$m_*\in\mathcal{P}^1(\td)$.
	A pair $(\hat{\varphi},\hat{m}_\cdot)$,
	where $\hat{\varphi}:[t_0,T]\times\td\rightarrow \mathbb{R}$ and
	$\hat{m}_\cdot\in\mathcal{C}^0([t_0,T];\mathcal{P}^1(\td))$,
is a solution
(see Section \ref{subsection:notionofclassicalsolutionofmfgsystem})
of the MFG with initial time $t_0$
on the time interval $[t_0,T]$ and
initial distribution
$m_*$, if and only if there exists a bounded,
Borel function $\hat{w}:[t_0,T]\times\td\rightarrow \rd$
such that
	\begin{enumerate}[label=(\roman*)]
		\item\label{sol_mfg:def:solution:cond:t_0} $\hat{m}_{t_0}=m_*$;
		\item\label{sol_mfg:def:solution:cond:T}
		$\hat{\varphi}(T,\cdot)=G(\cdot,\hat{m}_T)$;
		\item\label{sol_mfg:def:solution:cond:value} \(
        \hat\varphi(s,y)=\big(A^{s,T}_{\hat{m}_{\cdot}}
        G(\cdot,\hat m_T)\big)(y)\)
        for all \(s\in[t_0,T)\) and all \(y\in\td\);
		\item\label{sol_mfg:def:solution:cond:admissible}
		$(\hat m_\cdot,\hat w)\in\Adm{t_0}{T}$;
		\item\label{sol_mfg:def:solution:cond:link}
		\(
        \action{m_{*}}{\hat\varphi(t_0,\cdot)}=
		I^a[t_0,T,\hat{m}_\cdot,\hat{w},F(\cdot,\hat{m}_{\cdot}),
		G(\cdot,\hat{m}_T)]\).
	\end{enumerate}
\end{proposition}
	In this context, we
	call $\hat{w}$ an \emph{equilibrium feedback strategy}.

\begin{remark}\label{rem:on_the_statement}
The following comments are in order. First, the continuity of the path
\(\hat{m}_\cdot\) is a hypothesis of the ``if'' direction and cannot
be dismissed, because conditions
\ref{sol_mfg:def:solution:cond:T}--\ref{sol_mfg:def:solution:cond:link}
evaluate \(\hat{m}_{\cdot}\) at the individual time instants \(T\) and \(t_0\)
while condition
\ref{sol_mfg:def:solution:cond:admissible} by itself constrains
\(\hat{m}_{\cdot}\) only up to a set of times of measure zero.
Second, part of the outcome of the proof will be that,
for a solution $(\hat{\varphi},\hat{m}_{\cdot})$
of the MFG on \([t_0,T]\), one actually has
\begin{equation}\label{sol_mfg:def:solution:cond:linkpropagated}
\action{\hat{m}_{t}}{\hat\varphi(t,\cdot)}=
		I^a[t,T,\hat{m}_\cdot,\hat{w},F(\cdot,\hat{m}_{\cdot}),
		G(\cdot,\hat{m}_T)]
\end{equation}
for all \(t\in [t_0,T]\). Third, condition
\ref{sol_mfg:def:solution:cond:T} follows from
condition
\ref{sol_mfg:def:solution:cond:value} in the limit
\(s\to T\) by the elementary estimate
\(|\operatorname{Val}[s,T,y,F(\cdot,\hat{m}_\cdot),\sigma]
-\sigma(y)|\) \(\leq (T-s)(M_F+L_C)+K_\sigma(C(T-s)
\) \(+\sqrt{2d(T-s)})\), valid for the maximization over
\(\Adm[C]{s}{T}\), together with the continuity of \(\hat{m}\) and
\ref{sol_mfg:cond:G}. However,
we keep the two conditions separate for readability and referencing convenience.
\end{remark}

\begin{proof}[Proof of Proposition \ref{sol_mfg:def:solution}]
Throughout the proof we write
\[
f(t,x)\triangleq F(x,\hat{m}_t),\qquad
\sigma\triangleq G(\cdot,\hat{m}_T),
\]
which is legitimate in both directions, as the pair
\((\hat{\varphi},\hat{m}_\cdot)\), and in particular the path
\(\hat{m}_\cdot\), is given in both. We first record the following
observation, used in both directions.

\emph{Observation}. If \(\hat{m}_\cdot\) is continuous and
\((\hat{m}_\cdot,\hat{w})\in\Adm{t_0}{T}\) for some feedback control
\(\hat{w}\) with \(\|\hat{w}\|_\infty\triangleq c<\infty\), then
\(\hat{m}_\cdot\) is \(\tfrac12\)-H\"older continuous with constant
\(c\sqrt{T}+\sqrt{2d}\): by Proposition
\ref{prel:corollary:link_adm}, the pair is an element of
\(\RAdm[c]{t_0}{T}\), whose flow admits a \(\tfrac12\)-H\"older continuous version by Proposition \ref{prel:prop:Holder_m},
so the path is this continuous version. Consequently,
by \ref{sol_mfg:cond:F}, the function \(f\) is Lipschitz continuous
w.r.t.~\(x\) with constant \(K_F\), uniformly in \(t\), and
\(\tfrac{1}{2}\)-H\"older continuous w.r.t.~\(t\), uniformly in \(x\);
and, by \ref{sol_mfg:cond:G} and \ref{sol_mfg:cond:G2plusalpha},
\(\sigma\in\operatorname{Lip}_{K_G}(\td)\cap
\mathcal{C}^{2+\alpha}(\td)\). That is, \((f,\sigma)\) satisfies the
assumptions \ref{cond:V:f}--\ref{cond:V:sigma} of Lemma
\ref{prel:lm:verification} on \([t_0,T]\).

We now prove the \underline{\smash{sufficiency}} statement,
so let $(\hat{\varphi},\hat{m}_\cdot)$ be a pair
that satisfies the five items of the statement, with
\(\hat{w}\) the corresponding feedback control.
By condition \ref{sol_mfg:def:solution:cond:admissible} and the
Observation, Lemma \ref{prel:lm:verification} applies to
\((f,\sigma)\); let \(\varphi\) be the classical solution of the
terminal-value problem \eqref{eq:HJBfsigma} provided by part
\ref{lm:verification:a}, that is, of
	\begin{equation}\label{eq:bellmanwithmhat}
	\partial_t \hat\varphi+\triangle \hat\varphi+H(x,D
	\hat\varphi)+F(x,\hat{m}_t)=0
	\end{equation}
    with terminal condition
    \(\hat{\varphi}(T,\cdot)=G(\cdot,\hat{m}_T)\)
    (we display the problem in the notation of the statement, since
	in a moment we will
	show that \(\hat{\varphi}\) solves it). By Lemma
\ref{prel:lm:verification}\ref{lm:verification:a} and the definition
\eqref{sol_mfg:intro:Lax_Oleink} of the operator \(A\),
\[
\varphi(s,y)=\operatorname{Val}[s,T,y,f,\sigma]
=\big(A^{s,T}_{\hat{m}_\cdot}G(\cdot,\hat{m}_T)\big)(y),
\qquad (s,y)\in[t_0,T]\times\dtorus ,
\]
so conditions \ref{sol_mfg:def:solution:cond:value} and
\ref{sol_mfg:def:solution:cond:T} say precisely that
\(\hat{\varphi}=\varphi\) on \([t_0,T]\times\dtorus\). In particular,
\(\hat{\varphi}\) possesses the regularity listed in Section
\ref{subsection:notionofclassicalsolutionofmfgsystem} and satisfies
the Bellman equation \eqref{eq:bellmanwithmhat} pointwise on
\((t_0,T)\times\dtorus\) with the terminal condition
\ref{sol_mfg:def:solution:cond:T}.
Moving on to the FPK equation part,
by conditions \ref{sol_mfg:def:solution:cond:t_0},
\ref{sol_mfg:def:solution:cond:link} and \eqref{eq:Valismphi},
\[
I^a[t_0,T,\hat{m}_\cdot,\hat{w},f,\sigma]
=\action{m_*}{\hat{\varphi}(t_0,\cdot)}
=\action{m_*}{\varphi(t_0,\cdot)}
=\operatorname{Val}[t_0,T,m_*,f,\sigma],
\]
i.e.~\((\hat{m}_\cdot,\hat{w})\) is a maximizer
with the path starting from \(m_*\). By Lemma
\ref{prel:lm:verification}\ref{lm:verification:c},
\begin{equation}\label{eq:whatisDpH}
\hat{w}(t,x)=w^*(t,x)=D_pH\big(x,D\hat{\varphi}(t,x)\big)
\qquad\text{for a.e.~}t\in[t_0,T]
\text{ and }\hat{m}_t\text{-a.e.~}x\in\td .
\end{equation}
By condition \ref{sol_mfg:def:solution:cond:admissible},
\((\hat{m}_\cdot,\hat{w})\) satisfies \eqref{prel:eq:FPK} in the weak
sense; since the drift enters the distributional formulation
\eqref{eq:FPKsolcond} only through the measure
\(\hat{m}_t(dx)\,dt\), the identity \eqref{eq:whatisDpH} lets us
replace \(\hat{w}\) by \(D_pH(x,D\hat{\varphi})\) there, so that
\(\hat{m}_\cdot\) satisfies the FPK equation of the system
\eqref{intro:eq:mfg_system} in
\(\mathcal{D}'((t_0,T)\times\dtorus)\). As, by hypothesis, \(\hat{m}_\cdot\) is
continuous with \(\hat{m}_{t_0}=m_*\), the pair
\((\hat{\varphi},\hat{m}_\cdot)\) is a classical solution of the MFG
in the sense of Section
\ref{subsection:notionofclassicalsolutionofmfgsystem}. Finally, the
identity \eqref{sol_mfg:def:solution:cond:linkpropagated} for all
\(t\) is \eqref{eq:Lambdazero} of Lemma
\ref{prel:lm:verification}\ref{lm:verification:c}.

Now let us prove the
\underline{\smash{necessity}} part of the
proposition, so let
$(\hat{\varphi},\hat{m}_\cdot)$
be a solution
of the MFG with initial time $t_0$
on the time interval $[t_0,T]$ and
initial distribution
$m_*$, in the sense of Section
\ref{subsection:notionofclassicalsolutionofmfgsystem}. Conditions
\ref{sol_mfg:def:solution:cond:t_0} and
\ref{sol_mfg:def:solution:cond:T}, and the continuity of
\(\hat{m}_\cdot\), are part of that definition. Define
\[
\hat{w}(t,x)\ \triangleq\
D_pH\big(x,D\hat{\varphi}(t,x)\big),
\qquad (t,x)\in[t_0,T]\times\dtorus .
\]
By the regularity required of a classical solution, \(D\hat{\varphi}\)
is continuous on the closed cylinder, so
\(\bar{K}\triangleq\max_{[t_0,T]\times\td}|D\hat{\varphi}|\) is
finite, and, by Lemma \ref{prel:lm:Hamiltonian}(ii)--(iii),
\(\hat{w}\) is continuous on \([t_0,T]\times\dtorus\) --- hence Borel
--- and bounded by \(C_*[\bar{K}]\).

\emph{Condition \ref{sol_mfg:def:solution:cond:admissible}}. By the
definition of a classical solution, \(\hat{m}_\cdot\) satisfies the
FPK equation of \eqref{intro:eq:mfg_system} in
\(\mathcal{D}'((t_0,T)\times\dtorus)\), whose drift is precisely
\(\hat{w}\); this is equation \eqref{prel:eq:FPK} with \(w=\hat{w}\)
in the weak sense \eqref{eq:FPKsolcond}. Since \(\hat{m}_\cdot\) is
continuous (hence weakly measurable) and \(\hat{w}\) is bounded and
Borel, \((\hat{m}_\cdot,\hat{w})\in\Adm{t_0}{T}\).

\emph{Condition \ref{sol_mfg:def:solution:cond:value}}. By the
Observation, Lemma \ref{prel:lm:verification} applies to
\((f,\sigma)\). The function \(\hat{\varphi}\) is a solution of the
terminal-value problem
\eqref{eq:HJBfsigma} with this data,
possessing the regularity required in
part \ref{lm:verification:a},
thus fulfilling definition of a classical solution of the MFG
in regards to the Bellman equation; by the uniqueness statement of
\ref{lm:verification:a},
\[
\hat{\varphi}(s,y)=\operatorname{Val}[s,T,y,f,\sigma]
=\big(A^{s,T}_{\hat{m}_\cdot}G(\cdot,\hat{m}_T)\big)(y)
\qquad\text{for all }(s,y)\in[t_0,T]\times\dtorus,
\]
which is condition \ref{sol_mfg:def:solution:cond:value}. (In
particular, a posteriori,
\(\bar{K}\leq\widehat{K}_{\sigma,f}\) and
\(|\hat{w}|\leq C_*[\widehat{K}_{\sigma,f}]=C\), by Proposition
\ref{prel:prop:val_lip} and Lemma \ref{prel:lm:Hamiltonian}(ii).)

\emph{Condition \ref{sol_mfg:def:solution:cond:link} and the identity
\eqref{sol_mfg:def:solution:cond:linkpropagated}}. The pair
\((\hat{m}_\cdot,\hat{w})\) belongs to \(\Adm{t_0}{T}\) and its
feedback coincides --- everywhere on \([t_0,T]\times\dtorus\), by
construction --- with the optimal feedback
\(w^*=D_pH(\cdot,D\varphi(\cdot,\cdot))\) of Lemma
\ref{prel:lm:verification}\ref{lm:verification:b} (here
\(\varphi=\hat{\varphi}\) by the uniqueness just invoked). Hence, by
\ref{lm:verification:b} applied on the interval \([t,T]\) with
\(m_0=\hat{m}_t\), for every \(t\in[t_0,T)\),
\[
I^a[t,T,\hat{m}_\cdot,\hat{w},f,\sigma]
=\operatorname{Val}[t,T,\hat{m}_t,f,\sigma]
=\action{\hat{m}_t}{\hat{\varphi}(t,\cdot)},
\]
the last equality by \eqref{eq:Valismphi}. This is
\eqref{sol_mfg:def:solution:cond:linkpropagated}; at \(t=T\) it is
trivial, and at \(t=t_0\), combined with condition
\ref{sol_mfg:def:solution:cond:t_0}, it is condition
\ref{sol_mfg:def:solution:cond:link}. The proof is complete.
\end{proof}

\subsection{Equilibrium operators}\label{subsection:equilibriumops}
Let \(0\leq s<T\).
We define the functional
\(J_{s,T}\) on \(\RAdm[C]{s}{T}\times\RAdm[C]{s}{T}\)
by
\begin{align*}
	J_{s,T}[\beta,\gamma] \triangleq I^r[s,T,\beta,F(\cdot,\mathfrak{f}_{\cdot}[\gamma]),
	G(\cdot,\mathfrak{f}_T[\gamma])],
\end{align*}
that is (recall once more the notation \eqref{eq:actionnotation}
introduced in the Preliminaries and the remark around it):
\begin{align*}
  J_{s,T}[\beta,\gamma] \triangleq
  \action{\beta}{F(\cdot,\mathfrak{f}_{\cdot}[\gamma])-L(\cdot,\cdot)}
  + \action{\mathfrak{f}_T[\beta]}{G(\cdot,\mathfrak{f}_T[\gamma])};
\end{align*}
more generally for \(\sigma\in \mathcal{C}(\dtorus)\) the functional \(J_{s,r}^{\sigma}\)
on \(\RAdm[C]{s}{r}\times \RAdm[C]{s}{r}\) (\(s<r\leq T\)) will be given by:
\begin{align*}
	J_{s,r}^{\sigma}[\beta,\gamma] \triangleq I^r[s,r,\beta,F(\cdot,\mathfrak{f}_{\cdot}[\gamma]),
	\sigma],
\end{align*}
that is,
\begin{align*}
  J_{s,r}^{\sigma}[\beta,\gamma] \triangleq
\action{\beta}{F(\cdot,\mathfrak{f}_{\cdot}[\gamma])-L(\cdot,\cdot)}
  + \action{\mathfrak{f}_r[\beta]}{\sigma}.
\end{align*}
With these conventions,
\[
J_{s,T}[\beta,\gamma] = J_{s,T}^{G(\cdot,\mathfrak{f}_T[\gamma])}[\beta,\gamma],
\]
so the functional $J_{s,T}$ is a special case of the functionals
$J_{s,r}^{\sigma}$.
Also, for \(0\leq s<r<T\), any functional \(J^{\sigma}_{s,r}\) can also be
looked at as a functional on
\(
\RAdm[C]{s}{r} \times \RAdm[C]{s}{r}\)
or on \(\RAdm[C]{0}{T}\times
\RAdm[C]{0}{T}
\)
(see the last paragraph of Section \ref{subsection:relaxedcontrolprocesses}).

Due to the continuity of
\(\mathfrak{f}_t[\gamma]\) in both \(t\) and \(\gamma\)
(Propositions \ref{prel:prop:Holder_m} and \ref{prel:prop:continuity})
and the joint continuity of \(F\), \(G\) and \(L\) in all their
variables, \(J_{s,r}^{\sigma}[\beta,\gamma]\) is continuous
simultaneously in \(\beta\), \(\gamma\), \(\sigma\) and the time
endpoints
(for varying endpoints, in the sense made precise in Proposition
\ref{prel:prop:cont_I_r}); we will invoke this in the sections on the
value multifunction, and give the exact statements there.
\begin{definition}\label{defn:relaxedequilibria}
Let \(0\leq s<T\), \(\gamma\in\RAdm{s}{T}\). We will
call \(\gamma\) a \emph{relaxed equilibrium} if
\begin{align*}
	J_{s,T}[\gamma,\gamma] \geq J_{s,T}[\beta,\gamma]
\end{align*}
for all \(\beta\in\RAdm{s}{T}\) such that \(\mathfrak{f}_{s}[\gamma]=
\mathfrak{f}_{s}[\beta]\).
If \(0\leq s< r\leq T\) and \(\sigma\in \mathcal{C}(\dtorus)\), we will
call \(\gamma\) a \emph{relaxed equilibrium for \(\sigma\) on \([s,r]\)}
if
\begin{align*}
	J_{s,r}^{\sigma}[\gamma,\gamma]
	\geq J_{s,r}^{\sigma}[\beta,\gamma]
\end{align*}
for all \(\beta \in \RAdm{s}{r}\)
such that \(\mathfrak{f}_{s}[\gamma]=
\mathfrak{f}_{s}[\beta]\).
\end{definition}
The equilibrium property has a useful reformulation in terms of
the value function \eqref{prel:intro:value}.
Indeed, Lemma \ref{lemma:characterizationofrelaxedequilibria}
below shows that
a relaxed equilibrium is
nothing but a relaxed control that attains the value of the control
problem in which the interaction field is ``frozen'' at the flow generated
by the control itself.
\begin{lemma}\label{lemma:characterizationofrelaxedequilibria}
Let \(0\leq s<r\leq T\), \(\sigma\in\mathcal{C}(\dtorus)\),
\(\gamma\in\RAdm{s}{r}\), and set \(m_0\triangleq\mathfrak{f}_s[\gamma]\),
\[
f_{\gamma}(t,x)\triangleq F(x,\mathfrak{f}_t[\gamma]),
\qquad (t,x)\in[s,r]\times\dtorus .
\]
Then \(\gamma\) is a relaxed equilibrium for \(\sigma\) on \([s,r]\)
if and only if
\begin{equation}\label{eq:equilibriumisvalue}
I^r[s,r,\gamma,f_{\gamma},\sigma]
=\operatorname{Val}[s,r,m_0,f_{\gamma},\sigma],
\end{equation}
the value being the supremum \eqref{prel:intro:value} over the
class \(\Adm{s}{r}\). The same statement holds for relaxed equilibria
in the sense of the first part of Definition
\ref{defn:relaxedequilibria}, with \(r=T\) and
\(\sigma=G(\cdot,\mathfrak{f}_T[\gamma])\).
\end{lemma}
\begin{proof}
Observe first that, by the definition of the functionals,
\(J^{\sigma}_{s,r}[\beta,\gamma]=I^r[s,r,\beta,f_{\gamma},\sigma]\)
for every \(\beta\in\RAdm{s}{r}\); in particular
\(J^{\sigma}_{s,r}[\gamma,\gamma]=I^r[s,r,\gamma,f_{\gamma},\sigma]\).
Observe also that, since every feedback control is bounded,
\(\Adm{s}{r}=\bigcup_{C'>0}\Adm[C']{s}{r}\), so that
\begin{equation}\label{eq:Valassup}
\operatorname{Val}[s,r,m_0,f_{\gamma},\sigma]
=\sup_{C'>0}\,
\sup\big\{I^a[s,r,m_\cdot,w,f_{\gamma},\sigma]:
(m_\cdot,w)\in\Adm[C']{s}{r},\ m_s=m_0\big\}.
\end{equation}
Suppose \eqref{eq:equilibriumisvalue} holds, and let
\(\beta\in\RAdm{s}{r}\) with \(\mathfrak{f}_s[\beta]=m_0\); then
\(\beta\in\RAdm[C']{s}{r}\) for some \(C'>0\), and by Lemma
\ref{lemma:equality:value_RAdm} (applied with the constant \(C'\))
and \eqref{eq:Valassup},
\[
J^{\sigma}_{s,r}[\beta,\gamma]
=I^r[s,r,\beta,f_{\gamma},\sigma]
\ \leq\ \operatorname{Val}[s,r,m_0,f_{\gamma},\sigma]
=I^r[s,r,\gamma,f_{\gamma},\sigma]
=J^{\sigma}_{s,r}[\gamma,\gamma],
\]
so \(\gamma\) is a relaxed equilibrium for \(\sigma\) on \([s,r]\).

Conversely, suppose \(\gamma\) is a relaxed equilibrium for \(\sigma\)
on \([s,r]\). For every \(C'>0\), Lemma
\ref{lemma:equality:value_RAdm} and the equilibrium property give
\begin{align*}
J^{\sigma}_{s,r}[\gamma,\gamma]
 \geq & \
\sup\big\{I^r[s,r,\beta,f_{\gamma},\sigma]:\beta\in\RAdm[C']{s}{r},\,\mathfrak{f}_s[\beta]=m_0\big\}
\\
= & \ \sup\big\{I^a[\ldots]:(m_\cdot,w)\in\Adm[C']{s}{r},\ m_s=m_0\big\},
\end{align*}
and taking the supremum over \(C'>0\) and using \eqref{eq:Valassup},
\(J^{\sigma}_{s,r}[\gamma,\gamma]\geq
\operatorname{Val}[s,r,m_0,f_{\gamma},\sigma]\). The reverse
inequality holds because
if
\(\gamma\in\RAdm[C'']{s}{r}\), then, by Lemma
\ref{lemma:equality:value_RAdm} once more,
\(I^r[s,r,\gamma,f_{\gamma},\sigma]\leq
\operatorname{Val}[s,r,m_0,f_{\gamma},\sigma]\). This proves
\eqref{eq:equilibriumisvalue}.
\end{proof}

\subsubsection{Alternative definition of MFG
solution}\label{subsubsection:alternativedefinitionofMFGsolution}
Using these notions, we can reformulate the definition of solution of
a MFG.
\begin{proposition}\label{prop:solofMFGusingrelaxedequilibrium}
Let \(t_0\in[0,T)\), \(m_*\in\wassonespace\). The following
statements hold true.
\begin{enumerate}[(i)]
\item\label{item:soltoequilibrium} If
\((\hat{\varphi},\hat{m}_\cdot)\) is a solution of the MFG with
initial time \(t_0\) and initial distribution \(m_*\), and
\(\hat{w}\) is a corresponding equilibrium feedback strategy, then
the relaxed control \(\gamma\) obtained from
\((\hat{m}_\cdot,\hat{w})\) by Proposition \ref{prel:prop:link_adm}
is a relaxed equilibrium with \(\mathfrak{f}_{t_0}[\gamma]=m_*\).
\item\label{item:equilibriumtosol} Conversely, if
\(\gamma\in\RAdm{t_0}{T}\) is a relaxed equilibrium with
\(\mathfrak{f}_{t_0}[\gamma]=m_*\), then, putting
\[
\hat{m}_\cdot\triangleq\mathfrak{f}_\cdot[\gamma],
\quad
\hat{w}\triangleq\mathfrak{v}[\gamma],
\quad
\hat{\varphi}(s,\cdot)\triangleq
A^{s,T}_{\hat{m}_\cdot}G(\cdot,\hat{m}_T)\ \ (t_0\leq s<T),
\quad
\hat{\varphi}(T,\cdot)\triangleq G(\cdot,\hat{m}_T),
\]
the pair \((\hat{\varphi},\hat{m}_\cdot)\) is a solution of the MFG
with initial time \(t_0\) and initial distribution \(m_*\), and
\(\hat{w}\) is an equilibrium feedback strategy for it.
\end{enumerate}
In particular, the MFG with initial measure $m_*$ at
time $t_0$ admits a
solution if and only if there exists a relaxed equilibrium
\(\gamma\in\RAdm{t_0}{T}\) with \(\mathfrak{f}_{t_0}[\gamma]=m_*\).
\end{proposition}
\begin{proof}
In both parts we write
\[
f(t,x)\triangleq F(x,\hat{m}_t),\qquad
\sigma\triangleq G(\cdot,\hat{m}_T),
\]
where \(\hat{m}_\cdot\) is the path appearing in the corresponding
part; note that with this notation \(f=f_{\gamma}\) and
\(\sigma=G(\cdot,\mathfrak{f}_T[\gamma])\) in the sense of Lemma
\ref{lemma:characterizationofrelaxedequilibria}, since
\(\hat{m}_\cdot=\mathfrak{f}_\cdot[\gamma]\) in both parts (in part
\ref{item:soltoequilibrium}, because the flow
$\mathfrak{f}_{\cdot}[\gamma]$ coincides a.e.~with the continuous path \(\hat{m}_\cdot\),
hence everywhere).

\underline{\smash{Part \ref{item:soltoequilibrium}}}. Let
\(\gamma\) be the probability measure constructed from
\((\hat{m}_\cdot,\hat{w})\) as in Proposition \ref{prel:prop:link_adm};
then
\(\mathfrak{f}_{t_0}[\gamma]=\hat{m}_{t_0}=m_*\) (condition
\ref{sol_mfg:def:solution:cond:t_0} of Proposition
\ref{sol_mfg:def:solution}) and, by the construction
\eqref{prel:intr:adm_embedding},
\[
I^r[t_0,T,\gamma,f,\sigma]
=I^a[t_0,T,\hat{m}_\cdot,\hat{w},f,\sigma].
\]
By condition \ref{sol_mfg:def:solution:cond:link},
\(I^a[t_0,T,\hat{m}_\cdot,\hat{w},f,\sigma]
=\action{m_*}{\hat{\varphi}(t_0,\cdot)}\); by condition
\ref{sol_mfg:def:solution:cond:value} together with Lemma
\ref{prel:lm:verification}\ref{lm:verification:a} (the pair
\((f,\sigma)\) satisfies \ref{cond:V:f}--\ref{cond:V:sigma}, by the
Observation at the beginning of the proof of Proposition
\ref{sol_mfg:def:solution}), the function
\(\hat{\varphi}(t_0,\cdot)\) coincides with
\(\operatorname{Val}[t_0,T,\cdot,f,\sigma]\), and hence, by
\eqref{eq:Valismphi},
\(\action{m_*}{\hat{\varphi}(t_0,\cdot)}
=\operatorname{Val}[t_0,T,m_*,f,\sigma]\).
Thus,
\[
I^r[t_0,T,\gamma,f,\sigma]=\operatorname{Val}[t_0,T,m_*,f,\sigma],
\]
and Lemma \ref{lemma:characterizationofrelaxedequilibria} shows that
\(\gamma\) is a relaxed equilibrium.

\underline{\smash{Part \ref{item:equilibriumtosol}}}. Let
\(\gamma\in\RAdm[C']{t_0}{T}\) for some \(C'>0\) be a relaxed
equilibrium with \(\mathfrak{f}_{t_0}[\gamma]=m_*\), and let
\(\hat{m}_\cdot\), \(\hat{w}\), \(\hat{\varphi}\) be as in the
statement. The path \(\hat{m}_\cdot\) is continuous (Proposition
\ref{prel:prop:Holder_m}), \(\hat{w}\) is bounded and Borel, and
\((\hat{m}_\cdot,\hat{w})\in\Adm{t_0}{T}\) with
\(\hat{m}_{t_0}=m_*\) (Corollary \ref{prel:corollary:link_adm});
thus, conditions \ref{sol_mfg:def:solution:cond:t_0} and
\ref{sol_mfg:def:solution:cond:admissible} of Proposition
\ref{sol_mfg:def:solution} are met.
Conditions
\ref{sol_mfg:def:solution:cond:T} and
\ref{sol_mfg:def:solution:cond:value} hold by the definition of
\(\hat{\varphi}\) so it only remains to verify condition
\ref{sol_mfg:def:solution:cond:link}. We have:
\[
\operatorname{Val}[t_0,T,m_*,f,\sigma]
\ \geq\ I^a[t_0,T,\hat{m}_\cdot,\hat{w},f,\sigma]
\ \geq\ I^r[t_0,T,\gamma,f,\sigma]
\ =\ \operatorname{Val}[t_0,T,m_*,f,\sigma],
\]
where the first inequality holds by the definition
\eqref{prel:intro:value} of the value, the second is the Jensen
computation preceding Lemma \ref{lemma:equality:value_RAdm}, and the
equality is Lemma \ref{lemma:characterizationofrelaxedequilibria}.
Hence
\begin{equation}\label{eq:equilibriumpairismaximizer}
I^a[t_0,T,\hat{m}_\cdot,\hat{w},f,\sigma]
=\operatorname{Val}[t_0,T,m_*,f,\sigma],
\end{equation}
i.e.~\((\hat{m}_\cdot,\hat{w})\) is a maximizer
from \(m_*\).
Furthermore, \((f,\sigma)\) satisfies
\ref{cond:V:f}--\ref{cond:V:sigma} (the Observation in the proof of
Proposition \ref{sol_mfg:def:solution} is applicable here again because
\((\hat{m}_\cdot,\hat{w})\in\Adm{t_0}{T}\) and \(\hat{m}_\cdot\) is
continuous), so Lemma
\ref{prel:lm:verification}\ref{lm:verification:a}--\ref{lm:verification:b}
applies: the classical solution \(\varphi\) of \eqref{eq:HJBfsigma}
on \([t_0,T]\) satisfies
\(\varphi(t_0,y)=\operatorname{Val}[t_0,T,y,f,\sigma]
=\hat{\varphi}(t_0,y)\) and
\(\operatorname{Val}[t_0,T,m_*,f,\sigma]
=\action{m_*}{\varphi(t_0,\cdot)}\). Combining with
\eqref{eq:equilibriumpairismaximizer},
\[
\action{m_*}{\hat{\varphi}(t_0,\cdot)}
=I^a[t_0,T,\hat{m}_\cdot,\hat{w},f,\sigma],
\]
which is condition \ref{sol_mfg:def:solution:cond:link}. By the
sufficiency part of Proposition \ref{sol_mfg:def:solution}, the pair
\((\hat{\varphi},\hat{m}_\cdot)\) is a solution of the MFG with
initial time \(t_0\) and initial distribution \(m_*\), and
\(\hat{w}\) is an equilibrium feedback strategy for it.
\end{proof}

\begin{corollary}\label{sol_mfg:cor:structure_equilibria}
Let \(\gamma\in\RAdm{t_0}{T}\) be a relaxed equilibrium with
\(\mathfrak{f}_{t_0}[\gamma]=m_*\), and let
\(\hat{m}_\cdot\triangleq\mathfrak{f}_\cdot[\gamma]\),
\(\hat{w}\triangleq\mathfrak{v}[\gamma]\),
\(\hat{\varphi}\) be as in Proposition
\ref{prop:solofMFGusingrelaxedequilibrium}\ref{item:equilibriumtosol}.
Then:
\begin{enumerate}[label=(\alph*)]
\item\label{item:structure:dirac} for a.e.~\(t\in[t_0,T]\) and
\(\hat{m}_t\)-a.e.~\(x\in\td\), the disintegration of \(\gamma\)
satisfies \(\gamma(dv|t,x)=\delta_{\hat{w}(t,x)}(dv)\); that is,
\(\gamma\) coincides with the relaxed control constructed from
\((\hat{m}_{\cdot},\hat{w})\) in
\eqref{prel:intr:adm_embedding} of Proposition \ref{prel:prop:link_adm};
\item\label{item:structure:feedback}
\(\hat{w}(t,x)=D_pH\big(x,D\hat{\varphi}(t,x)\big)\) and
\(|\hat{w}(t,x)|\leq C\) for a.e.~\(t\in[t_0,T]\) and
\(\hat{m}_t\)-a.e.~\(x\in\td\), where \(C\) is the constant fixed at
the beginning of this section;
\item\label{item:structure:class} consequently,
\(\gamma\in\RAdm[C]{t_0}{T}\).
\end{enumerate}
\end{corollary}
\begin{proof}
Let \(f\), \(\sigma\) be defined as in the preceding proof.

Part \ref{item:structure:dirac}: We are
in the context of part \ref{item:equilibriumtosol} of
Proposition \ref{prop:solofMFGusingrelaxedequilibrium}.
The inequalities in the proof of that part that,
because they turned out to
be equalities, led to
\eqref{eq:equilibriumpairismaximizer}, show that
\(I^a[t_0,T,\hat{m}_\cdot,\hat{w},f,\sigma]
=I^r[t_0,T,\gamma,f,\sigma]\). Applying the definitions
\eqref{eq:defnofIa} and \eqref{eq:defnofIr},
\[
0=I^a-I^r=\int_{t_0}^T\int_{\td}
\Big[\int_{\rd}L(x,v)\,\gamma(dv|t,x)
-L\big(x,\hat{w}(t,x)\big)\Big]\hat{m}_t(dx)\,dt.
\]
The expression inside the square parentheses
is non-negative by Jensen's inequality
(as \(\hat{w}(t,x)\) is the the barycenter of \(\gamma(dv|t,x)\)),
hence it vanishes for a.e.~\(t\) and \(\hat{m}_t\)-a.e.~\(x\).
Moreover, at any such
\((t,x)\), by \ref{prel:cond:L_convex} and the inequality
\[
L(x,v)\geq L(x,\hat{w}(t,x))
+D_vL(x,\hat{w}(t,x))\cdot(v-\hat{w}(t,x))
\]
with equality holding if and
only if \(v=\hat{w}(t,x)\) (owing to the strict convexity of $L$),
an integration with respect to the measure
\(\gamma(dv|t,x)\) shows that the non-negative integrand
\(v\mapsto L(x,v)-L(x,\hat{w}(t,x))
-D_vL(x,\hat{w}(t,x))\cdot(v-\hat{w}(t,x))\) has vanishing
\(\gamma(\cdot|t,x)\)-integral, so
\(\gamma(\cdot|t,x)=\delta_{\hat{w}(t,x)}\).

Part \ref{item:structure:feedback}: by
\eqref{eq:equilibriumpairismaximizer}, the pair
\((\hat{m}_\cdot,\hat{w})\) is a maximizer
among admissible controls starting from \(m_*\); Lemma
\ref{prel:lm:verification}\ref{lm:verification:c} identifies
\(\hat{w}(t,x)=D_pH(x,D\varphi(t,x))=D_pH(x,D\hat{\varphi}(t,x))\)
a.e.~(with \(\varphi=\hat{\varphi}\) as in the preceding proof), and
Proposition \ref{prop:prel:boundness} gives
\(|\hat{w}|\leq C_*[\widehat{K}_{\sigma,f}]\leq C\), with the last
inequality explained by the facts that
\(K_f=K_F\), \(K_\sigma=K_G\), \(T-t_0\leq T\) and both
\(C_*[\cdot]\), \(\widehat{K}_{\sigma,f}\) are nondecreasing in their
arguments.

Part \ref{item:structure:class}:
by \ref{item:structure:dirac} and
\ref{item:structure:feedback}, the set
\(\{(t,x,v):|v|>C\}\) has zero $\gamma$ measure, so
\(\spt\gamma\subset[t_0,T]\times\dtorus\times\mathbb{B}_C\); the
remaining conditions of Definition
\ref{prel:def:relaxed_control_procs} hold because \(\gamma\) is
a relaxed control on $[t_0,T]$.
\end{proof}

\begin{remark}\label{sol_mfg:prop:continuity_sol}
For later use we collect the bounds available for an arbitrary
solution \((\hat{\varphi},\hat{m}_\cdot)\) of the MFG with initial
data \((t_0,m_*)\) and equilibrium feedback strategy \(\hat{w}\); all
the constants below depend only on the data \textbf{(HS)} and on
\(T\), and not on the particular solution nor on \((t_0,m_*)\).
\begin{itemize}[-~~]
\item \(\hat{\varphi}\in\mathcal{C}^0([t_0,T]\times\dtorus)\), by the
definition of a classical solution (Section
\ref{subsection:notionofclassicalsolutionofmfgsystem}); moreover
\(|\hat{\varphi}|\leq M_G+T(M_F+M^0_L+\sup_{|v|\le C}|L(\cdot,v)-L(\cdot,0)|)\),
since \(\hat\varphi\) is the value function of a problem with data
bounded through \eqref{sol_mfg:ineq:F_bound},
\eqref{sol_mfg:ineq:G_bound} and \eqref{sol_mfg:intro:M_L}.
\item For each \(s\in[t_0,T]\), the function
\(\hat{\varphi}(s,\cdot)\) is Lipschitz continuous with constant
\[
K[s]\triangleq K_G+(T-s)(K_F+K_L)\leq
\widehat{K}_{G(\cdot,\hat{m}_T),F} \ ,
\]
indeed this is condition
\ref{sol_mfg:def:solution:cond:value} combined with Proposition
\ref{prel:prop:val_lip}.
\item \(|\hat{w}(t,x)|\leq C\) for a.e.~\(t\) and
\(\hat{m}_t\)-a.e.~\(x\): Proposition \ref{prop:prel:boundness}
applied to the maximizer \((\hat{m}_\cdot,\hat{w})\) (that it is one
was established in the proofs of Propositions
\ref{sol_mfg:def:solution} and
\ref{prop:solofMFGusingrelaxedequilibrium}).
\item \(\hat{m}_\cdot\) is \(\tfrac12\)-H\"older continuous with
constant \(C\sqrt{T}+\sqrt{2d}\): the
relaxed control \(\gamma\) obtained from
\((\hat{m}_\cdot,\hat{w})\) by Proposition \ref{prel:prop:link_adm},
used in Proposition
\ref{prop:solofMFGusingrelaxedequilibrium}\ref{item:soltoequilibrium},
lies in \(\RAdm[C]{t_0}{T}\) by Corollary
\ref{sol_mfg:cor:structure_equilibria}, and Proposition
\ref{prel:prop:Holder_m} applies to it; the H\"older representative
coincides with the continuous path \(\hat{m}_\cdot\).
\end{itemize}
In particular, in the alternative definition of solution of the MFG
provided by Proposition \ref{prop:solofMFGusingrelaxedequilibrium},
one may equivalently quantify over relaxed equilibria
\(\gamma\in\RAdm[C]{t_0}{T}\), since by Corollary
\ref{sol_mfg:cor:structure_equilibria}\ref{item:structure:class}
every relaxed equilibrium lies in this
class.
\end{remark}

\subsection{Value multifunction for the
MFG}\label{subsection:valuemultifunctionfortheMFG}
Now let us introduce the concept of value multifunction,
the main object of this work.
\begin{definition}\label{sol_mdf:def:value_multifunction} (Value multifunction) We say that a
	multivalued mapping $\mathcal{V}:[0,T]\times
	\mathcal{P}^1(\td)\rightrightarrows C(\td)$
	is a \emph{value multifunction} if:
	\begin{enumerate}
		\item $\mathcal{V}(T,\mu)=\{G(\cdot,\mu)\}$;
		\item for each $s\in [0,T)$,
		$\mu\in\mathcal{P}^1(\td)$
		and $\psi\in \mathcal{V}(s,\mu)$
		there exists a solution \((\hat{\varphi},\hat{m}(\cdot))\)
		of the mean field game with the initial time $s$
		and initial distribution $\mu$ at time \(s\)
		(i.e.~\(\hat{m}_s = \mu\))
		 such that
		\[
		\hat{\varphi}(s,\cdot)
		=\psi(\cdot).
		\]
	\end{enumerate}
\end{definition}
The next statement gives the definition of value multifunction
through relaxed equilibria.
\begin{corollary}\label{sol_mfg:corollary:value_mfg}
	A multivalued mapping
	$\mathcal{V}:[0,T]\times \mathcal{P}^1(\td)\rightrightarrows C(\td)$
 is a value multifunction if and only if:
	\begin{enumerate}
		\item $\mathcal{V}(T,\mu)=\{G(\cdot,\mu)\}$ for all
		\(\mu\in\wassonespace\);
		\item for each \(s\in [0,T)\), \(\mu\in\wassonespace\),
		\(\psi\in\mathcal{V}(s,\mu)\) there
		exists \(\gamma \in \RAdm[C]{s}{T}\) such that
		\begin{enumerate}[(a)]
			\item \(\mathfrak{f}_s[\gamma] = \mu\) and
			\(\gamma\) is a relaxed equilibrium,
			\item
			\(\psi(y) =
			(A^{s,T}_{\mathfrak{f}_{\cdot}[\gamma]}G(\cdot,
			\mathfrak{f}_T[\gamma]))(y)\),
			\(y\in\dtorus\).
		\end{enumerate}
	\end{enumerate}
\end{corollary}
\begin{proof}
If \(\mathcal{V}\) is a value multifunction and
\(\psi\in\mathcal{V}(s,\mu)\), \(s<T\), let
\((\hat{\varphi},\hat{m}(\cdot))\) be a solution as in Definition
\ref{sol_mdf:def:value_multifunction} and \(\hat{w}\) an equilibrium
feedback strategy for it. By Proposition
\ref{prop:solofMFGusingrelaxedequilibrium}\ref{item:soltoequilibrium},
the relaxed control \(\gamma\) constructed from
\((\hat{m}_\cdot,\hat{w})\) is a relaxed
equilibrium with \(\mathfrak{f}_s[\gamma]=\mu\), and it lies in
\(\RAdm[C]{s}{T}\) by Corollary
\ref{sol_mfg:cor:structure_equilibria}\ref{item:structure:class};
condition
\ref{sol_mfg:def:solution:cond:value} of Proposition
\ref{sol_mfg:def:solution} together with
\(\mathfrak{f}_\cdot[\gamma]=\hat{m}_\cdot\) gives (b)
of the present corollary.
Conversely,
given \(\gamma\) as in condition 2, Proposition
\ref{prop:solofMFGusingrelaxedequilibrium}\ref{item:equilibriumtosol}
produces a solution \((\hat{\varphi},\hat{m}_\cdot)\) with
\(\hat{m}_s=\mu\) and
\(\hat{\varphi}(s,\cdot)=A^{s,T}_{\mathfrak{f}_\cdot[\gamma]}
G(\cdot,\mathfrak{f}_T[\gamma])=\psi\).
\end{proof}
\begin{remark}\label{sol_mfg:prop:bound_Lipschitz}
If \(\mathcal{V}\) is a value multifunction,
\(t_0\in[0,T]\), \(m_*\in\wassonespace\) then any value
\(\psi\in\mathcal{V}(t_0,m_*)\) possesses these bounds and regularity:
	\begin{itemize}[-~~]
		\item $\|\psi\|\leq M_G+(T-t_0)(M_F+L_C)$: indeed
		\(\psi(y)=\operatorname{Val}[t_0,T,y,F(\cdot,\hat{m}_\cdot),
		G(\cdot,\hat{m}_T)]\) for a solution \(\hat{m}_\cdot\) as in
		Definition \ref{sol_mdf:def:value_multifunction}, the value may
		be computed over \(\Adm[C]{t_0}{T}\) (Proposition
		\ref{prop:prel:boundness}), and over that class
		\(|I^a|\leq(T-t_0)(M_F+L_C)+M_G\), by
		\eqref{sol_mfg:ineq:F_bound}, \eqref{sol_mfg:ineq:G_bound} and
		\eqref{sol_mfg:intro:M_L};
		\item $\psi$ is Lipschitz continuous with constant
		\(K[t_0]\) (defined in Remark
		\ref{sol_mfg:prop:continuity_sol}): this is Proposition
		\ref{prel:prop:val_lip} with \(K_f=K_F\), \(K_\sigma=K_G\),
		\(r-s=T-t_0\);
		\item for \(t_0<T\), the function \(\psi\) is, in fact, twice
		continuously differentiable: \(\psi=\hat{\varphi}(t_0,\cdot)\)
		for a classical solution \(\hat{\varphi}\), whose spatial
		derivatives up to second order are continuous on the closed
		cylinder (Section
		\ref{subsection:notionofclassicalsolutionofmfgsystem}). We
		record this for completeness; it will not be used below.
	\end{itemize}
\end{remark}

\begin{lemma}\label{lemma:stability_RAdm}
Let \(0\leq s_n<r_n\leq T\), \(s_n\rightarrow s\),
\(r_n\rightarrow r\) with \(s<r\), and let
\(\gamma_n\in\RAdm[C]{s_n}{r_n}\) converge, in the topology
\(\mathcal{T}\) of Section \ref{subsubsection:topologyRAdm} (all the
sets involved being subsets of \(B^{\ast}_{0,T}\)), to a limit
\(\gamma\). Then:
\begin{enumerate}[label=(\alph*)]
\item\label{item:stability:limit}
the limit \(\gamma\in\RAdm[C]{s}{r}\);
\item\label{item:stability:flows} the flows converge uniformly, in
the sense that
\[
\sup_{t\in[s,r]}
W_1\big(\mathfrak{f}_{\pi_n(t)}[\gamma_n],\mathfrak{f}_t[\gamma]\big)
\ \longrightarrow\ 0,
\qquad
\pi_n(t)\triangleq\min\{\max\{t,s_n\},r_n\};
\]
in particular
\(\mathfrak{f}_{s_n}[\gamma_n]\rightarrow\mathfrak{f}_s[\gamma]\) and
\(\mathfrak{f}_{r_n}[\gamma_n]\rightarrow\mathfrak{f}_r[\gamma]\) in
\(W_1\);
\item\label{item:stability:payoffs} if
\(f_n:[s_n,r_n]\times\td\rightarrow\R\), \(\sigma_n:\td\to\R\) are
continuous and converge uniformly to \(f\) and \(\sigma\)
respectively (with the \(f_n\) being extended to a common interval by
their values at the endpoints), then
\[
I^r[s_n,r_n,\gamma_n,f_n,\sigma_n]\ \longrightarrow\
I^r[s,r,\gamma,f,\sigma].
\]
\end{enumerate}
\end{lemma}
Notice that, by design, $\pi_n(t)$ in
part \ref{item:stability:flows} is defined
so as to make the case $f_n(t,x)=F(x,\mathfrak{f}_{\pi_n(t)}
[\gamma_n])$ consistent with part
\ref{item:stability:payoffs} (extension of $f_n$
to a common interval by their values at the endpoints).
\begin{proof}
\ref{item:stability:limit}. The limit is represented by a
non-negative Borel measure \(\gamma\) on
\([0,T]\times\dtorus\times\mathbb{B}_C\), by the argument of Section
\ref{subsubsection:topologyRAdm}. Since for every
\(g\in\mathcal{C}([0,T])\),
\(\action{\gamma_n}{g\circ\pj^1}=\int_{s_n}^{r_n}g(t) dt\)
converges to \(\int_s^r g(t) dt\), we have
\(\pj^1\sharp\gamma=\mathscr{L}^1_{[s,r]}\); in particular the total
mass is \(r-s\). With the time projection being supported on \([s,r]\)
while the remaining coordinates supported in the compact
\(\dtorus\times\mathbb{B}_C\), we have
\(\spt\gamma\subset[s,r]\times\dtorus\times\mathbb{B}_C\). Finally,
if \(\varphi\in\mathcal{C}^{1,2}_c((s,r)\times\dtorus)\), then
\(\spt\varphi\subset[s+\delta,r-\delta]\times\dtorus\) for some
\(\delta>0\); for \(n\) large,
\([s+\delta,r-\delta]\subset(s_n,r_n)\), so
\(\varphi\in\mathcal{C}^{1,2}_c((s_n,r_n)\times\dtorus)\) and
\eqref{eq:weakcontrolFPK} holds for \(\gamma_n\); passing to the
limit (note that
the integrand is a fixed continuous function) yields
\eqref{eq:weakcontrolFPK} for \(\gamma\).

\ref{item:stability:flows}. Fix \(\eps\in(0,(r-s)/2)\) and set
\(J_{\eps}\triangleq[s+\eps,r-\eps]\), so that
\(J_{\eps}\subset(s_n,r_n)\) for \(n\) large. We claim that
\(\gamma_n|_{J_{\eps}}\rightarrow\gamma|_{J_{\eps}}\) in the topology of
\(\RAdm[C]{s+\eps}{r-\eps}\): for
\(g\in\mathcal{C}([0,T]\times\dtorus\times\mathbb{B}_C)\) and
\(\delta>0\), let \(\chi_\delta:[0,T]\to[0,1]\) be continuous, equal
to \(1\) on \([s+\eps+\delta,r-\eps-\delta]\) and to \(0\) outside
\(J_{\eps}\); since the time projections of all the measures involved are
dominated by Lebesgue measure,
\(|\action{\gamma_n|_{J_{\eps}}}{g}-\action{\gamma_n}{g\chi_\delta}|
\leq2\delta\|g\|\) uniformly in \(n\), and likewise for \(\gamma\),
while \(\action{\gamma_n}{g\chi_\delta}\rightarrow
\action{\gamma}{g\chi_\delta}\); letting \(\delta\to 0\) proves
the claim. By Proposition \ref{prel:prop:continuity} on the fixed
interval \(J_{\eps}\) (and since the flow of a restriction coincides
with the flow of the original control on the subinterval),
\[
d_n\triangleq
\sup_{t\in J_{\eps}}W_1\big(\mathfrak{f}_t[\gamma_n],
\mathfrak{f}_t[\gamma]\big)\ \longrightarrow\ 0 .
\]
Now let \(t\in[s,r]\) and let \(t'\) be the point of \(J_{\eps}\)
nearest to \(t\), so that \(|t-t'|\leq\eps\) and
\(|\pi_n(t)-t'|\leq|s_n-s|+|r_n-r|+\eps\). By the triangle
inequality and the uniform \(\tfrac12\)-H\"older bound of
Proposition \ref{prel:prop:Holder_m}, with the constant
\(d\triangleq C\sqrt{T}+\sqrt{2d}\) valid for all the
controls involved,
\[
W_1\big(\mathfrak{f}_{\pi_n(t)}[\gamma_n],\mathfrak{f}_t[\gamma]\big)
\leq
d\big(|s_n-s|+|r_n-r|+\eps\big)^{1/2}
+d_n+d\,\eps^{1/2},
\]
and from this, \(\limsup_n\sup_{t\in[s,r]}
W_1(\mathfrak{f}_{\pi_n(t)}[\gamma_n],\mathfrak{f}_t[\gamma])
\leq2d\,\eps^{1/2}\) (a soon as
\(|s_n-s|+|r_n-r|\leq\eps\)); as \(\eps\) was arbitrary, the claim
follows.

\ref{item:stability:payoffs}. Denoting by \(\tilde{f}_n\),
\(\tilde{f}\) the extensions to a common interval \([s',r']\)
containing all the intervals involved, we note that
\begin{align*}
	& \ \action{\gamma_n}{f_n-L}=\action{\gamma_n}{\tilde{f}_n-L}
\\
= & \ \action{\gamma_n}{\tilde{f}-L}
+ O\big(\|\tilde f_n-\tilde f\|\,(r'-s')\big)
\ \longrightarrow \ \action{\gamma}{\tilde{f}-L}=\action{\gamma}{f-L},
\end{align*}
where the middle convergence is just by the definition of the topology,
while
\(\action{\mathfrak{f}_{r_n}[\gamma_n]}{\sigma_n}\rightarrow
\action{\mathfrak{f}_r[\gamma]}{\sigma}\) by part
\ref{item:stability:flows} and the uniform convergence
\(\sigma_n\rightarrow\sigma\).
\end{proof}

\begin{definition}\label{sol_mfg:def:closure}
If $\mathcal{V}:[0,T]\times\mathcal{P}^1(\td)\rightrightarrows C(\td)$,
then its closure $\operatorname{cl}\mathcal{V}$ is defined by the rule:
$\psi\in (\operatorname{cl}\mathcal{V})(t_0,m_*)$ if there exist sequences
$\{t_{0,n}\}_{n=1}^\infty$, $\{m_{*,n}\}_{n=1}^\infty$
and $\{\psi_n\}_{n=1}^\infty$
such that
\begin{itemize}[-~~]
	\item $\psi_n\in\mathcal{V}(t_{0,n},m_{*,n})$;
	\item $t_{0,n}\rightarrow t_0$, $m_{*,n}\rightarrow
	m_*$ in \(W_1\), and $\psi_n\rightarrow\psi$ uniformly on
	\(\dtorus\),
	as $n\rightarrow\infty$.
\end{itemize}
We say that a multifunction $\mathcal{V}$ is closed if
$\mathcal{V}=\operatorname{cl}\mathcal{V}$.
\end{definition}

\begin{proposition}\label{sol_mfg:prop:closure} If $\mathcal{V}:[0,T]
	\times\mathcal{P}^1(\td)\rightrightarrows C(\td)$ is a
	value multifunction, then its closure $\operatorname{cl}\mathcal{V}$
	is also a value multifunction.
\end{proposition}
\begin{proof}
	We first verify condition 1 of Corollary
	\ref{sol_mfg:corollary:value_mfg} for
	\(\operatorname{cl}\mathcal{V}\). The inclusion
	\(G(\cdot,\mu)\in(\operatorname{cl}\mathcal{V})(T,\mu)\) is seen by
	taking constant sequences. Conversely, let
	\(\psi\in(\operatorname{cl}\mathcal{V})(T,\mu)\), with sequences
	\(t_{0,n}\rightarrow T\), \(m_{*,n}\rightarrow\mu\),
	\(\psi_n\in\mathcal{V}(t_{0,n},m_{*,n})\), \(\psi_n\rightarrow\psi\).
	For the indices with \(t_{0,n}<T\), Corollary
	\ref{sol_mfg:corollary:value_mfg} provides relaxed equilibria
	\(\gamma_n\in\RAdm[C]{t_{0,n}}{T}\) with
	\(\mathfrak{f}_{t_{0,n}}[\gamma_n]=m_{*,n}\) and
	\(\psi_n=A^{t_{0,n},T}_{\mathfrak{f}_\cdot[\gamma_n]}
	G(\cdot,\mathfrak{f}_T[\gamma_n])\); then
	\(W_1(\mathfrak{f}_T[\gamma_n],m_{*,n})\leq
	(C\sqrt{T}+\sqrt{2d})(T-t_{0,n})^{1/2}\rightarrow0\)
	(Proposition \ref{prel:prop:Holder_m}), so
	\(\mathfrak{f}_T[\gamma_n]\rightarrow\mu\), and, by the elementary
	estimate recorded in Remark \ref{rem:on_the_statement} together
	with \ref{sol_mfg:cond:G},
	\[
	\|\psi_n-G(\cdot,\mu)\|
	\leq
	\big\|A^{t_{0,n},T}_{\mathfrak{f}_\cdot[\gamma_n]}
	G(\cdot,\mathfrak{f}_T[\gamma_n])
	-G(\cdot,\mathfrak{f}_T[\gamma_n])\big\|
	+K_G\,W_1\big(\mathfrak{f}_T[\gamma_n],\mu\big)
	\ \longrightarrow\ 0
	\]
	(for the indices with \(t_{0,n}=T\), directly
	\(\psi_n=G(\cdot,m_{*,n})\rightarrow G(\cdot,\mu)\)); hence
	\(\psi=G(\cdot,\mu)\), and
	\((\operatorname{cl}\mathcal{V})(T,\mu)=\{G(\cdot,\mu)\}\).

	Now let $t_0\in [0,T)$, $m_*\in\mathcal{P}^1(\td)$ and
	$\psi\in\operatorname{cl}\mathcal{V}(t_0,m_*)$, with sequences
	\(t_{0,n}\rightarrow t_0\), \(m_{*,n}\rightarrow m_*\),
	\(\psi_n\in\mathcal{V}(t_{0,n},m_{*,n})\),
	\(\psi_n\rightarrow\psi\); discarding finitely many indices we may
	assume \(t_{0,n}<T\).
	Since $\mathcal{V}$ is a value multifunction, by
	Corollary \ref{sol_mfg:corollary:value_mfg}, for each $n$, one can find a relaxed
	equilibrium \(\gamma_n\) in \(\RAdm[C]{t_{0,n}}{T}\) such that
	\begin{equation}\label{sol_mfg:equality:m_star_n}
		m_{*,n}=\mathfrak{f}_{t_{0,n}}[\gamma_n]
	\end{equation}
	and
	\begin{equation}\label{sol_mfg:equality:val_n}
		\psi_n=A^{t_{0,n},T}_{\mathfrak{f}_\cdot[\gamma_n]}
		G(\cdot,\mathfrak{f}_T[\gamma_n]).
	\end{equation}
	All the \(\gamma_n\) lie in the compact metrizable set
	\(B^{\ast}_{0,T}\) (see the comments at the end of Section
	\ref{subsubsection:topologyRAdm}); passing to a subsequence,
	which we do not relabel, \(\gamma_n\rightarrow\gamma\) in
	\(\mathcal{T}\), and Lemma
	\ref{lemma:stability_RAdm}\ref{item:stability:limit} gives
	\(\gamma\in\RAdm[C]{t_0}{T}\). By Lemma
	\ref{lemma:stability_RAdm}\ref{item:stability:flows} and
	\eqref{sol_mfg:equality:m_star_n},
	\begin{equation*}
		m_*=\mathfrak{f}_{t_0}[\gamma],
	\end{equation*}
	and also
	\(\mathfrak{f}_T[\gamma_n]\rightarrow\mathfrak{f}_T[\gamma]\),
	whence, by \ref{sol_mfg:cond:F} and \ref{sol_mfg:cond:G}, the data
	\[
	f_n(t,x)\triangleq F(x,\mathfrak{f}_{\pi_n(t)}[\gamma_n]),
	\qquad
	\sigma_n\triangleq G(\cdot,\mathfrak{f}_T[\gamma_n])
	\]
	(\(\pi_n\) as in Lemma \ref{lemma:stability_RAdm}) converge
	uniformly to
	\(f(t,x)\triangleq F(x,\mathfrak{f}_t[\gamma])\) and
	\(\sigma\triangleq G(\cdot,\mathfrak{f}_T[\gamma])\)
	respectively.

	The equilibrium property passes to the limit through Lemma
	\ref{lemma:characterizationofrelaxedequilibria}: for each \(n\),
	\[
	I^r[t_{0,n},T,\gamma_n,f_n,\sigma_n]
	=\operatorname{Val}[t_{0,n},T,m_{*,n},f_n,\sigma_n]
	=\action{m_{*,n}}
	{\operatorname{Val}[t_{0,n},T,\cdot,f_n,\sigma_n]},
	\]
	the second equality by
	\eqref{eq:Valismphi} of Lemma
	\ref{prel:lm:verification}\ref{lm:verification:b}, whose
	hypotheses \ref{cond:V:f}--\ref{cond:V:sigma} hold for
	\((f_n,\sigma_n)\) (the flows are \(\tfrac12\)-H\"older by
	Proposition \ref{prel:prop:Holder_m}, and \ref{sol_mfg:cond:F},
	\ref{sol_mfg:cond:G}, \ref{sol_mfg:cond:G2plusalpha} apply).
	As \(n\rightarrow\infty\): the left-hand side converges to
	\(I^r[t_0,T,\gamma,f,\sigma]\), by Lemma
	\ref{lemma:stability_RAdm}\ref{item:stability:payoffs}; the
	functions
	\(\operatorname{Val}[t_{0,n},T,\cdot,f_n,\sigma_n]\) converge to
	\(\operatorname{Val}[t_0,T,\cdot,f,\sigma]\) uniformly on
	\(\dtorus\), by Proposition \ref{prel:prop:cont_Val}, and are
	equi-Lipschitz (Proposition \ref{prel:prop:val_lip}), so that,
	together with \(m_{*,n}\rightarrow m_*\), the right-hand side
	converges to
	\(\action{m_*}{\operatorname{Val}[t_0,T,\cdot,f,\sigma]}
	=\operatorname{Val}[t_0,T,m_*,f,\sigma]\)
	(\eqref{eq:Valismphi} once more). Hence
	\(I^r[t_0,T,\gamma,f,\sigma]
	=\operatorname{Val}[t_0,T,m_*,f,\sigma]\), and Lemma
	\ref{lemma:characterizationofrelaxedequilibria} shows that
	\(\gamma\) is a relaxed equilibrium.

	Finally, passing to the limit in \eqref{sol_mfg:equality:val_n}
	(Proposition \ref{prel:prop:cont_Val} again, together with
	\(\psi_n\rightarrow\psi\)), we deduce that
	\begin{equation*}
		\label{sol_mfg:equality:val_limit}
		\psi=A^{t_0,T}_{\mathfrak{f}_\cdot[\gamma]}G(\cdot,\mathfrak{f}_T[\gamma]).
	\end{equation*}
	Conditions 2(a) and 2(b) of Corollary
	\ref{sol_mfg:corollary:value_mfg} are thus fulfilled at
	\((t_0,m_*,\psi)\), and \(\operatorname{cl}\mathcal{V}\) is a
	value multifunction.
\end{proof}

\section{Mean field game dynamics}\label{sect:mfg_dyn}
In this section we introduce the concept of mean field
game dynamics (Definition \ref{mfg_dyn:def:dynamics}),
which can be interpreted as realizing
a \textit{forward-forward} mean field
game \cite{Gomes2016OneDimensional,Gomes2018Conservation},
we formulate the MFG dynamic programming principle
(Definition \ref{defn:DPP}) for
a closed multifunction as a viability condition on it
with respect to the dynamics operator, and then show
(Theorem \ref{thm:aetoDPP})
that a closed multifunction that satisfies this principle
must be a value multifunction for the MFG
(Definition \ref{sol_mdf:def:value_multifunction}).

Recall once again the definition of the constant \(K[s]\)
from Remark \ref{sol_mfg:prop:continuity_sol},
and the constants \(M_F\), \(M_G\), and \(L_C\)
from \eqref{sol_mfg:ineq:F_bound},
\eqref{sol_mfg:ineq:G_bound}
and \eqref{sol_mfg:intro:M_L} respectively.

\begin{corollary}\label{mfg_dyn:cor:verification_identification}
Let \(0\leq s<r\leq T\), let \(f\) satisfy \ref{cond:V:f}
of Lemma \ref{prel:lm:verification} on
\([s,r]\), and suppose \(\varphi:[s,r]\times\dtorus\to\R\) possesses the
regularity required of solutions in part \ref{lm:verification:a}
of that lemma (\(\varphi\in\mathcal{C}^0([s,r]\times\dtorus)\), continuously
differentiable in \(t\) on \((s,r)\times\dtorus\), with \(D\varphi\),
\(D_{xx}\varphi\) continuous on the closed cylinder) and satisfies
the PDE
\[
\partial_t\varphi+\triangle\varphi+H(x,D\varphi)+f(t,x)=0
\qquad\text{pointwise on }(s,r)\times\dtorus ;
\]
set
\(\sigma_0\triangleq\varphi(r,\cdot)\in\mathcal{C}(\dtorus)\), so
that \(\varphi\) solves the terminal-value problem
\eqref{eq:HJBfsigma} with datum \(\sigma_0\) (and, by Remark
\ref{rem:dtextension}, \(\partial_t\varphi\) extends continuously to
the closed cylinder). Then \(\varphi\) is the value function of its
own terminal datum:
\[
\varphi(t,y)=\operatorname{Val}[t,r,y,f,\sigma_0]
\qquad\text{for all }t\in[s,r],\ y\in\dtorus.
\]
Moreover
\(\operatorname{Val}[t,r,m_0,f,\sigma_0]
=\action{m_0}{\varphi(t,\cdot)}\) for every \(t\in[s,r]\) and
\(m_0\in\ptd\), and, for every \(t\in[s,r)\), the feedback
\(D_pH(\cdot,D\varphi(\cdot,\cdot))\) is a maximizer on \([t,r]\)
from every initial distribution.
\end{corollary}
\begin{proof}
Steps 1 and 2 of the proof of Lemma \ref{prel:lm:verification} use
only the stated regularity of \(\varphi\), the equation, and property
\ref{cond:V:f} of \(f\); the hypothesis \ref{cond:V:sigma} on the
terminal datum enters that proof only through the existence assertion,
which is here replaced by the given \(\varphi\).
\end{proof}

The significance of the next lemma is that
its properties are thus far valid for terminal data
$\sigma$ fulfilling condition
\ref{cond:V:sigma}
(i.e.~\(\sigma\in\operatorname{Lip}_{K_\sigma}(\td)\)
and \(\sigma\in\mathcal{C}^{2+\alpha'}(\td)\) for some
\(\alpha'\in(0,1)\)), but the lemma
extends them to merely Lipschitz terminal data.

\begin{lemma}\label{mfg_dyn:lm:A_properties}
Let \(0\leq s<r\leq T\), let \(\hat{m}_\cdot:[s,r]\rightarrow\ptd\)
be \(\tfrac12\)-H\"older continuous, \(f(t,x)\triangleq
F(x,\hat{m}_t)\), and let \(\sigma\in\operatorname{Lip}_{K_\sigma}(\td)\)
with \(K_\sigma\leq K[r]\). Then:
\begin{enumerate}[label=(\alph*)]
\item\label{item:A:classes} the value admits the representations
\begin{align*}
\operatorname{Val}[s,r,m_0,f,\sigma]
=&\ \sup\big\{I^a[s,r,m_\cdot,w,f,\sigma]:(m_\cdot,w)\in\Adm[C]{s}{r},
\ m_s=m_0\big\}
\\
=&\ \max\big\{I^r[s,r,\beta,f,\sigma]:\beta\in\RAdm[C]{s}{r},
\ \mathfrak{f}_s[\beta]=m_0\big\}
\end{align*}
for every \(m_0\in\ptd\), the last maximum being attained;
\item\label{item:A:onelip} each of the functionals in
\ref{item:A:classes} is \(1\)-Lipschitz in \(\sigma\) with respect to
the supremum norm; in particular
\(\|A^{s,r}_{\hat{m}}\sigma-A^{s,r}_{\hat{m}}\sigma'\|
\leq\|\sigma-\sigma'\|\);
\item\label{item:A:regularity}
\(A^{s,r}_{\hat{m}}\sigma\in
\operatorname{Lip}_{K_\sigma+(r-s)(K_F+K_L)}(\td)\) and
\(\|A^{s,r}_{\hat{m}}\sigma\|\leq\|\sigma\|+(r-s)(M_F+L_C)\);
\item\label{item:A:decoupling}
\(\operatorname{Val}[s,r,m_0,f,\sigma]
=\action{m_0}{A^{s,r}_{\hat{m}}\sigma}\) for every \(m_0\in\ptd\);
\item\label{item:A:semigroup} if \(s<r<\theta\leq T\),
\(\hat{m}_\cdot\) is \(\tfrac12\)-H\"older on \([s,\theta]\) and
\(\sigma\in\operatorname{Lip}_{K_\sigma}(\td)\) with
\(K_\sigma\leq K[\theta]\), then
\[
A^{s,\theta}_{\hat{m}}\sigma
=A^{s,r}_{\hat{m}}\big(A^{r,\theta}_{\hat{m}}\sigma\big);
\]
\item\label{item:A:continuity} Proposition \ref{prel:prop:cont_Val}
remains valid if the requirement that each \(\sigma_n\) belong to
\(\mathcal{C}^{2+\alpha_n}(\td)\) is dropped and we require
the functions \(\sigma_n\) to be
uniformly Lipschitz.
\end{enumerate}
\end{lemma}
\begin{proof}
Throughout, for \(\eps>0\) let
\(\sigma_\eps\triangleq\rho_\eps\ast\sigma\) be a mollification on
the torus: \(\sigma_\eps\in\mathcal{C}^\infty(\td)\),
\(\operatorname{Lip}(\sigma_\eps)\leq K_\sigma\), and
\(\|\sigma_\eps-\sigma\|\leq K_\sigma\eps\); in particular each
\(\sigma_\eps\) satisfies \ref{cond:V:sigma}, and \((f,\sigma_\eps)\)
satisfies the hypotheses of Lemma \ref{prel:lm:verification}.

Part \ref{item:A:onelip} is immediate: each of \(I^a\), \(I^r\) is affine
in \(\sigma\), and
\(|\action{m}{\sigma-\sigma'}|\leq\|\sigma-\sigma'\|\) for any
probability \(m\); suprema of \(1\)-Lipschitz families are
\(1\)-Lipschitz.

Part \ref{item:A:classes}. We treat
the maximization problem with \(\sigma_\eps\) in place
of $\sigma$ first. The constant of
Proposition \ref{prop:prel:boundness} for the data
\((f,\sigma_\eps)\) is
\[
C_*[K_\sigma+(r-s)(K_F+K_L)]\leq C_*[K[r]+(r-s)(K_F+K_L)]
=C_*[K[s]]\leq C_*[K[0]]=C
\]
(due to monotonicity of \(C_*[\cdot]\), Lemma
\ref{prel:lm:Hamiltonian}(ii) and recall
the constants $K[\cdot]$ were defined in
Remark \ref{sol_mfg:prop:continuity_sol} and $\hat{K}$
was introduced in front of
Lemma \ref{prel:lm:verification}); by Lemma
\ref{prel:lm:verification}\ref{lm:verification:b}, the
supremum is attained by a feedback bounded by this constant, so the full and
the \(C\)-restricted values coincide, and coincide with the supremum
over \(\{\beta\in\RAdm[C]{s}{r}:\mathfrak{f}_s[\beta]=m_0\}\) by
Lemma \ref{lemma:equality:value_RAdm}. By
\ref{item:A:onelip}, all three functionals change by at most
\(K_\sigma\eps\) when \(\sigma_\eps\) is replaced back by \(\sigma\);
letting \(\eps\to 0\) transfers the equalities to \(\sigma\).
The maximum over
\(\{\beta\in\RAdm[C]{s}{r}:\mathfrak{f}_s[\beta]=m_0\}\) is attained
because this set is compact (it is closed in the compact
\(\RAdm[C]{s}{r}\), the map \(\beta\mapsto\mathfrak{f}_s[\beta]\)
being continuous by Proposition \ref{prel:prop:continuity}) and
\(\beta\mapsto I^r[s,r,\beta,f,\sigma]\) is continuous (here, again,
the running integrand is a fixed continuous function, and
\(\beta\mapsto\mathfrak{f}_r[\beta]\) is continuous).

Part \ref{item:A:regularity} is Proposition \ref{prel:prop:val_lip}
together with the bound \(|I^a|\leq\|\sigma\|+(r-s)(M_F+L_C)\), valid
over \(\Adm[C]{s}{r}\), and combined with \ref{item:A:classes}.

Part \ref{item:A:decoupling}. For \(\sigma_\eps\) this is
the decoupling identity
\eqref{eq:Valismphi} of Lemma
\ref{prel:lm:verification}\ref{lm:verification:b};
for $\sigma$ the statement follows since both sides are
\(1\)-Lipschitz in the terminal datum. Namely,
by part \ref{item:A:onelip} above
$\operatorname{Val}[s,r,m_0,f,\sigma_{\eps}]$ converges
to $\operatorname{Val}[s,r,m_0,f,\sigma]$ as $\eps\to 0$;
for the
right-hand side,
\(|\action{m_0}{A\sigma-A\sigma'}|\leq\|A\sigma-A\sigma'\|
\leq\|\sigma-\sigma'\|\). Letting
\(\eps\to 0\) we obtain \ref{item:A:decoupling}.

Part \ref{item:A:semigroup}. First, suppose that \(\sigma\) satisfies
\ref{cond:V:sigma} as well, and let \(\varphi\) be the classical
solution of \eqref{eq:HJBfsigma} on \([s,\theta]\) with terminal
datum \(\sigma\) (Lemma
\ref{prel:lm:verification}\ref{lm:verification:a}). Applying
Corollary \ref{mfg_dyn:cor:verification_identification} on
\([s,\theta]\) and on \([s,r]\) (the restriction
\(\varphi|_{[s,r]\times\dtorus}\) has the required regularity and
solves the same equation, with terminal datum \(\varphi(r,\cdot)\)),
\[
A^{s,\theta}_{\hat{m}}\sigma=\varphi(s,\cdot),\qquad
A^{r,\theta}_{\hat{m}}\sigma=\varphi(r,\cdot),\qquad
A^{s,r}_{\hat{m}}\big(\varphi(r,\cdot)\big)=\varphi(s,\cdot),
\]
but this is exactly the claim.
Now, for \(\sigma\) merely Lipschitz, we may apply the case
just treated to \(\sigma_\eps\) and pass to the limit using
\ref{item:A:onelip} twice:
\(A^{r,\theta}\sigma_\eps\rightarrow A^{r,\theta}\sigma\) uniformly,
hence
\(A^{s,r}(A^{r,\theta}\sigma_\eps)\rightarrow
A^{s,r}(A^{r,\theta}\sigma)\), while
\(A^{s,\theta}\sigma_\eps\rightarrow A^{s,\theta}\sigma\). Notice that the
Lipschitz classes are respected in these steps:
\(\operatorname{Lip}(A^{r,\theta}\sigma_\eps)\leq
K_\sigma+(\theta-r)(K_F+K_L)\leq K[r]\), by
\ref{item:A:regularity}.

Part \ref{item:A:continuity}. In the setting of Proposition
\ref{prel:prop:cont_Val}, take a sequence \(\{\sigma_n\}_{n=1}^{\infty}\),
uniformly \(\bar{K}\)-Lipschitz, and converging to $\sigma$.
For fixed \(\eps>0\), the mollified
family \(\{\sigma_{n,\eps}\}\) satisfies the hypotheses of that
proposition (each \(\sigma_{n,\eps}\in\mathcal{C}^\infty\),
uniformly \(\bar{K}\)-Lipschitz,
\(\sigma_{n,\eps}\rightarrow\sigma_\eps\) uniformly), so
\(\operatorname{Val}[s_n,r_n,\cdot,f_n,\sigma_{n,\eps}]\rightarrow
\operatorname{Val}[s,r,\cdot,f,\sigma_\eps]\) uniformly.
Then, by
\ref{item:A:onelip}, replacing \(\sigma_{n,\eps}\) and \(\sigma_\eps\)
by \(\sigma_n\) and \(\sigma\), respectively,
changes each value function by at most
\(\bar{K}\eps\), uniformly in \(n\). Hence,
\(\limsup_n\|\operatorname{Val}[s_n,r_n,\cdot,f_n,\sigma_n]
-\operatorname{Val}[s,r,\cdot,f,\sigma]\|\leq2\bar{K}\eps\) for
every \(\eps>0\).
\end{proof}

\begin{remark}\label{mfg_dyn:rem:semigroup_general}
We stress that the restriction \(\sigma\in\operatorname{Lip}_{K_\sigma}(\td)\),
\(K_\sigma\leq K[\theta]\), in item \ref{item:A:semigroup}
owes itself to item \ref{item:A:classes}, that is, to the
\(C\)-restricted reading of the operators \(A\) used throughout this
section, and not from the semigroup property itself. Indeed, for the
value function alone, the identity
\(\operatorname{Val}[s,\theta,y,f,\sigma]
=\operatorname{Val}[s,r,y,f,\operatorname{Val}[r,\theta,\cdot,f,\sigma]]\)
holds for every \(\sigma\in\mathcal{C}(\td)\), by the argument of
item \ref{item:A:semigroup} together with item \ref{item:A:onelip}
and the sup-norm density of smooth functions in \(\mathcal{C}(\td)\).
This is consistent with the dynamic programming
principle of the general theory, which holds for bounded continuous
terminal data; see, e.g., \cite[Ch.~4]{YongZhou1999}.
\end{remark}

\begin{definition}\label{mfg_dyn:def:dynamics} (Forward operators) Let $s,r\in [0,T]$, $s<r$, $(\psi,\mu)\in C(\td)\times \mathcal{P}^1(\td)$
	be such that
	$\psi\in \operatorname{Lip}_{K[s]}(\td)$ and $\|\psi\|\leq M_G+(M_F+L_C)(T-s)$.
	We denote by
$$
\Phi^{r,s}(\psi,\mu)
$$
the set of all pairs
	$(\sigma,\nu)\in C(\td)\times \mathcal{P}^1(\td)$,
	with $\sigma\in\operatorname{Lip}_{K[r]}(\td)$ and $\|\sigma\|
		\leq M_G+(M_F+L_C)(T-r)$,
	for which there exists a relaxed equilibrium \(\gamma\)
    for $\sigma$ on $[s,r]$ in
	$\RAdm[C]{s}{r}$
    (see Definition \ref{defn:relaxedequilibria})
    fulfilling the conditions:
	\begin{enumerate}
		\item $\mathfrak{f}_s[\gamma]=\mu$;
		\item $\mathfrak{f}_r[\gamma]=\nu$;
		\item $\psi=A^{s,r}_{\mathfrak{f}_\cdot[\gamma]}\sigma$.
	\end{enumerate}
	For \((\psi,\mu)\) not satisfying the requirements above we set
	\(\Phi^{r,s}(\psi,\mu)\triangleq\emptyset\).
\end{definition}
We emphasize the domain and codomain of the operators
\(\Phi^{r,s}\):
\begin{align*}
	\Phi(r,s): \
	\begin{aligned}
	\{(\psi,\mu)\in C(\dtorus)\times \mathcal{P}^1(\dtorus)\ : & \  \
	\Psi\in \operatorname{Lip}_{K[s]},\|\psi\|\leq M_G+(M_F+L_C)(T-s)\}
	\\
	 \longdowndownarrows & \
	 \\
\{(\sigma,\nu)\in C(\dtorus)\times \mathcal{P}^1(\dtorus)\ : & \ \
	\Psi\in \operatorname{Lip}_{K[r]},\|\psi\|\leq M_G+(M_F+L_C)(T-r)\}.
	\end{aligned}
\end{align*}
The definition of the forward operators $\Phi^{r,s}$ can be reformulated as follows,
which will be convenient below in the work
on the infinitesimal version of the dynamic
programming principle.
\begin{remark}\label{mfg_dyn:remark:eq}
	Given $s,r\in [0,T]$, $s<r$, $(\psi,\mu)\in C(\td)\times \mathcal{P}^1(\td)$ such that
	$\psi\in \operatorname{Lip}_{K[s]}(\td)$ and $\|\psi\|\leq M_G+(M_F+L_C)(T-s)$, the set $\Phi^{r,s}(\psi,\mu)$
    consists of all pairs
	$(\sigma,\nu)\in C(\td)\times \mathcal{P}^1(\td)$ with $\sigma\in\operatorname{Lip}_{K[r]}(\td)$ and $\|\sigma\|
		\leq M_G+(M_F+L_C)(T-r)$,
	satisfying the following conditions for some \(\gamma\in
	\RAdm[C]{s}{r}\):
	\begin{enumerate}
		\item $\mathfrak{f}_s[\gamma]=\mu$;
		\item $\mathfrak{f}_r[\gamma]=\nu$;
		\item $\psi=A^{s,r}_{\mathfrak{f}_\cdot[\gamma]}\sigma$;
        \item
       $\action{\mu}{\psi} = \action{\mathfrak{f}_r[\gamma]}{\sigma} +
       \action{\gamma}{F(\cdot,\mathfrak{f}_{\cdot}[\gamma])
       -L(\cdot,\cdot)}$.
	\end{enumerate}
	In other words, in the presence of conditions 1--3, the fourth
	condition
	is equivalent to \(\gamma\) being a
	relaxed equilibrium for \(\sigma\) on \([s,r]\). This form of the
	equilibrium property is convenient for passing to limits, and it
	is the form that will be used in the sequel to this work, devoted
	to the infinitesimal analysis of the dynamics.

	Let us prove the equivalence. Write
	\(f_\gamma(t,x)\triangleq F(x,\mathfrak{f}_t[\gamma])\); note that
	the right-hand side of condition 4 is
	\(I^r[s,r,\gamma,f_\gamma,\sigma]\), while, by condition 3 and
	Lemma \ref{mfg_dyn:lm:A_properties}\ref{item:A:decoupling}
	(once more, this is applicable because \(\mathfrak{f}_\cdot[\gamma]\) is
	\(\tfrac12\)-H\"older by Proposition \ref{prel:prop:Holder_m},
	and \(\sigma\in\operatorname{Lip}_{K[r]}(\td)\)),
	\[
	\action{\mu}{\psi}
	=\action{\mu}{A^{s,r}_{\mathfrak{f}_\cdot[\gamma]}\sigma}
	=\operatorname{Val}[s,r,\mu,f_\gamma,\sigma].
	\]
	Thus, given conditions 1--3, condition 4 reads
	\(I^r[s,r,\gamma,f_\gamma,\sigma]
	=\operatorname{Val}[s,r,\mu,f_\gamma,\sigma]\), which, by Lemma
	\ref{lemma:characterizationofrelaxedequilibria}, holds if and
	only if \(\gamma\) is a relaxed equilibrium for \(\sigma\) on
	\([s,r]\).
\end{remark}
\subsection{Dynamic programming principle}\label{subsection:dppfordynamics}
Let us point out some key properties of the mean field dynamics.
The superposition operation of two multifunctions (the extension
of the notion of composition of functions) is this:
if $B_1,B_2$ are multifunctions from $X$ to $X$, then,
	for each $x\in X$,
	\[
	(B_2\circ B_1)(x)\triangleq \bigcup_{y\in B_1(x)}B_2(y),
	\]
	i.e.~\(z\in B_2\circ B_1(x)\) if and only if
	\(z\in B_2(y)\) for some \(y\in B_1(x)\).
\begin{proposition}\label{mfg_dyn:prop:superposition} If $s<r<\theta\leq T$, then
	\[
	\Phi^{\theta,s}=\Phi^{\theta,r}\circ\Phi^{r,s}.
	\]
\end{proposition}
\begin{proof}
Throughout the proof, for a relaxed control \(\gamma\) we write
\(f_\gamma(t,x)\triangleq F(x,\mathfrak{f}_t[\gamma])\), and we use
without further comment that flows of restrictions and concatenations
coincide with the original flows on the corresponding subintervals,
so that the frozen data \(f_{\gamma}\) restricts and concatenates
accordingly.

First let us prove that
\begin{equation}\label{mfg_dyn:incl:Phi_1}
	\Phi^{\theta,s}\subset\Phi^{\theta,r}\circ\Phi^{r,s}.
\end{equation}
Let $(\psi,\mu)\in C(\td)\times\mathcal{P}^1(\td)$,
with \(\psi\) in compliance with the required bounds.
If \(\Phi^{\theta,s}(\psi,\mu)
=\emptyset\), there is nothing to prove. Otherwise, suppose
$(\sigma,\nu)\in \Phi^{\theta,s}(\psi,\mu)$:
for some relaxed equilibrium \(\gamma\) for \(\sigma\) on
\([s,\theta]\) in $\RAdm[C]{s}{\theta}$,
	\begin{itemize}[-~~]
	\item $\sigma\in\operatorname{Lip}_{K[\theta]}(\td)$, $\|\sigma\|
	\leq M_G+(M_F+L_C)(T-\theta)$;
	\item $\mathfrak{f}_s[\gamma]=\mu$,
	$\mathfrak{f}_\theta[\gamma]=\nu$,
	$\psi=A^{s,\theta}_{\mathfrak{f}_\cdot[\gamma]}\sigma$.
\end{itemize}
Put
\[
\gamma_1\triangleq \gamma|_{[s,r]},\ \ \gamma_2\triangleq \gamma|_{[r,\theta]},
\qquad
\mu'\triangleq \mathfrak{f}_r[\gamma],
\qquad
\psi'\triangleq A^{r,\theta}_{\mathfrak{f}_\cdot[\gamma_2]}\sigma.
\]
By Lemma \ref{mfg_dyn:lm:A_properties}\ref{item:A:regularity},
\(\psi'\in\operatorname{Lip}_{K_1}(\td)\) with
\(K_1=K[\theta]+(\theta-r)(K_F+K_L)=K[r]\), and
\(\|\psi'\|\leq\|\sigma\|+(\theta-r)(M_F+L_C)\leq
M_G+(T-r)(M_F+L_C)\); by Lemma
\ref{mfg_dyn:lm:A_properties}\ref{item:A:semigroup},
\begin{equation}\label{eq:psiviasemigroup}
\psi=A^{s,\theta}_{\mathfrak{f}_\cdot[\gamma]}\sigma
=A^{s,r}_{\mathfrak{f}_\cdot[\gamma_1]}
\big(A^{r,\theta}_{\mathfrak{f}_\cdot[\gamma_2]}\sigma\big)
=A^{s,r}_{\mathfrak{f}_\cdot[\gamma_1]}\psi' .
\end{equation}
We now show that \(\gamma_1\) is a relaxed equilibrium for \(\psi'\)
on \([s,r]\) and \(\gamma_2\) is a relaxed equilibrium for \(\sigma\)
on \([r,\theta]\); together with the identities already displayed,
this gives \((\psi',\mu')\in\Phi^{r,s}(\psi,\mu)\) and
\((\sigma,\nu)\in\Phi^{\theta,r}(\psi',\mu')\), proving
\eqref{mfg_dyn:incl:Phi_1}. By Remark \ref{mfg_dyn:remark:eq}, it
suffices to verify the identities (condition 4 of the remark)
\begin{equation}\label{eq:twobudgets}
\action{\mu}{\psi}
=\action{\mu'}{\psi'}
+\action{\gamma_1}{f_\gamma-L},
\qquad
\action{\mu'}{\psi'}
=\action{\nu}{\sigma}
+\action{\gamma_2}{f_\gamma-L}.
\end{equation}
Consider the quantities
\[
B_1\triangleq
\action{\mu'}{\psi'}+\action{\gamma_1}{f_\gamma-L}
-\action{\mu}{\psi},
\qquad
B_2\triangleq
\action{\nu}{\sigma}+\action{\gamma_2}{f_\gamma-L}
-\action{\mu'}{\psi'} .
\]
Each one of them is nonpositive: indeed
\(\action{\mu'}{\psi'}+\action{\gamma_1}{f_\gamma-L}
=I^r[s,r,\gamma_1,f_\gamma,\psi']
\leq\operatorname{Val}[s,r,\mu,f_\gamma,\psi']
=\action{\mu}{A^{s,r}_{\mathfrak{f}_\cdot[\gamma_1]}\psi'}
=\action{\mu}{\psi}\),
by Lemma \ref{mfg_dyn:lm:A_properties}\ref{item:A:classes},
\ref{item:A:decoupling} and \eqref{eq:psiviasemigroup}, and likewise
\(B_2\leq0\) using
\(\psi'=A^{r,\theta}_{\mathfrak{f}_\cdot[\gamma_2]}\sigma\). On the
other hand, adding,
\[
B_1+B_2=\action{\nu}{\sigma}
+\action{\gamma}{f_\gamma-L}-\action{\mu}{\psi}=0,
\]
which is condition 4 of Remark \ref{mfg_dyn:remark:eq} for \(\gamma\)
(\(\gamma\) is a relaxed equilibrium for
\(\sigma\) on \([s,\theta]\) satisfying conditions 1--3
of that remark). Hence
\(B_1=B_2=0\), which is \eqref{eq:twobudgets}.

As for the opposite inclusion,
now suppose \((\psi,\mu)\), with the required bounds,
is in \(\mathcal{C}(\dtorus)\times\wassonespace\) and
\((\sigma,\nu)\in (\Phi^{\theta,r}\circ \Phi^{r,s})(\psi,\mu)\):
there is
\((\psi',\mu')\in \mathcal{C}(\dtorus)\times\wassonespace\)
such that \((\psi',\mu')\in \Phi^{r,s}(\psi,\mu)\) and
\((\sigma,\nu) \in \Phi^{\theta,r}(\psi',\mu')\), realized by
relaxed equilibria \(\gamma_1\in \RAdm[C]{s}{r}\) (for \(\psi'\)) and
\(\gamma_2 \in \RAdm[C]{r}{\theta}\) (for \(\sigma\)) with
\begin{align*}
	\mathfrak{f}_s[\gamma_1]=\mu,
	\qquad
	\mathfrak{f}_r[\gamma_1]=\mu',
	\qquad
	\psi = A^{s,r}_{\mathfrak{f}_{\cdot}[\gamma_1]}\psi',
	\\
	\mathfrak{f}_r[\gamma_2]=\mu',
	\qquad
	\mathfrak{f}_{\theta}[\gamma_2] =\nu,
	\qquad
	\psi' = A^{r,\theta}_{\mathfrak{f}_{\cdot}[\gamma_2]}\sigma.
\end{align*}
Let \(\gamma \triangleq \gamma_1 \diamond \gamma_2\) (Definition
\ref{prel:def:concatenation}; the matching condition holds). Then
\(\mathfrak{f}_{s}[\gamma]=\mu\),
\(\mathfrak{f}_{\theta}[\gamma]=\nu\), and, by Lemma
\ref{mfg_dyn:lm:A_properties}\ref{item:A:semigroup},
\[
A_{\mathfrak{f}_{\cdot}[\gamma]}^{s,\theta}\sigma
=
A_{\mathfrak{f}_{\cdot}[\gamma_1]}^{s,r}
(A_{\mathfrak{f}_{\cdot}[\gamma_2]}^{r,\theta}\sigma)
= A_{\mathfrak{f}_{\cdot}[\gamma_1]}^{s,r} \psi'
= \psi .
\]
It remains to show that \(\gamma\) is a relaxed equilibrium for
\(\sigma\) on \([s,\theta]\); once again by Remark \ref{mfg_dyn:remark:eq}, it
suffices to verify its fourth condition,
and indeed, adding the identities of this condition for
 \(\gamma_1\) (corresponding to \(\psi'\)) and for \(\gamma_2\) (corresponding to
\(\sigma\)), valid by the same remark,
\begin{align*}
\action{\nu}{\sigma}+\action{\gamma}{f_\gamma-L}
=& \ \Big(\action{\nu}{\sigma}+\action{\gamma_2}{f_\gamma-L}\Big)
+\action{\gamma_1}{f_\gamma-L}
\\
= & \ \action{\mu'}{\psi'}+\action{\gamma_1}{f_\gamma-L}
=\action{\mu}{\psi}.
\end{align*}
This concludes the proof.
\end{proof}

\begin{definition}\label{defn:DPP} (Dynamic programming principle for the MFG)
Let \(\mathcal{W}:[0,T]\times \wassonespace \rightsquigarrow
	\mathcal{C}(\dtorus)\) be a closed
	multifunction such that \(\mathcal{W}(T,\mu) =
	\{G(\cdot,\mu)\}\) for all \(\mu\in\wassonespace\). If the multifunction
	\(\mathcal{W}\)
	is such that for each \(s\in [0,T)\), \(\mu\in\wassonespace\),
	\(\phi\in\mathcal{W}(s,\mu)\), \(r\in (s,T)\), there exist
	\((\psi,\nu)\in \mathcal{C}(\dtorus)\times \wassonespace\) satisfying:
	\begin{enumerate}
	  \item \(\psi\in \mathcal{W}(r,\nu)\),
	  \item \((\psi,\nu)\in \Phi^{r,s}(\phi,\mu)\),
	\end{enumerate}
then we say that \(\mathcal{W}\) \emph{satisfies the dynamic programming
principle (DPP) for the mean-field game}.
\end{definition}
We note that, by the empty-set convention in Definition
\ref{mfg_dyn:def:dynamics}, the DPP forces every value
\(\phi\in\mathcal{W}(s,\mu)\), \(s<T\), to lie in
\(\operatorname{Lip}_{K[s]}(\td)\) with
\(\|\phi\|\leq M_G+(M_F+L_C)(T-s)\), since, otherwise,
\(\Phi^{r,s}(\phi,\mu)=\emptyset\).

\begin{theorem}\label{thm:aetoDPP}
	Let \(\mathcal{W}:[0,T]\times \wassonespace \rightsquigarrow
	\mathcal{C}(\dtorus)\) be a closed
	multifunction that satisfies the DPP for the mean-field game.
	Then \(\mathcal{W}\) is a value multifunction.
	\end{theorem}
	\begin{proof}
	Property 1 in the definition of value multifunction
    (Definition \ref{sol_mdf:def:value_multifunction})
	holds by hypothesis.
	Fix arbitrary \(s\in [0,T)\),
	\(\mu \in \wassonespace\), \(\phi\in \mathcal{W}(s,\mu)\)
	and consider an increasing sequence \(r_n\in(s,T)\),
	\(r_n \to T\).
	By the DPP and the definition of the forward operators
	(Definition \ref{mfg_dyn:def:dynamics}), for each \(n\) there are
	\((\psi_n,\nu_n)\in
	\mathcal{C}(\dtorus)\times \wassonespace\) and a
	relaxed equilibrium \(\gamma_n\in
	\RAdm[C]{s}{r_n}\) for \(\psi_n\) on \([s,r_n]\) such that
	\begin{enumerate}
	  \item \(\psi_n \in \mathcal{W}(r_n,\nu_n)\),
	  \(\psi_n\in\operatorname{Lip}_{K[r_n]}(\td)\),
	  \(\|\psi_n\|\leq M_G+(M_F+L_C)(T-r_n)\);
	  \item \(\mathfrak{f}_s[\gamma_n] = \mu\),
	  \(\mathfrak{f}_{r_n}[\gamma_n] =\nu_n\);
	  \item \(\phi= A_{\mathfrak{f}_{\cdot}[\gamma_n]}^{s,r_n}
	  \psi_n\).
	\end{enumerate}
	Condition 3 means, by Lemma
	\ref{mfg_dyn:lm:A_properties}\ref{item:A:classes}, that
\begin{align}
	\phi(y) = &\ \max\big\{
    I^r[s,r_n,\beta,F(\cdot,\mathfrak{f}_{\cdot}[\gamma_n]),
    \psi_n]
		\ : \
	\beta \in \RAdm[C]{s}{r_n}, \ \mathfrak{f}_s[\beta]=\delta_y \big\}
	\label{eq:displayphiofyn}
\end{align}
	   for all \(y\in \dtorus\).

	First we extract a limiting relaxed
    equilibrium.The \(\gamma_n\) lie in the compact
	metrizable set \(B^{\ast}_{0,T}\); passing to a subsequence,
	which we do not relabel, \(\gamma_n\rightarrow\gamma\), and Lemma
	\ref{lemma:stability_RAdm}\ref{item:stability:limit} gives
	\(\gamma\in\RAdm[C]{s}{T}\), with \(\mathfrak{f}_s[\gamma]=\mu\)
	and
	\(\nu_n=\mathfrak{f}_{r_n}[\gamma_n]\rightarrow
	\mathfrak{f}_T[\gamma]\triangleq\nu\)
	by Lemma \ref{lemma:stability_RAdm}\ref{item:stability:flows}.
	The functions \(\psi_n\) are equibounded and equi-Lipschitz (with constants
	\(K[r_n]\leq K[0]\)), by condition 1; by the Arzel\`a--Ascoli
	theorem, a further subsequence (not relabeled) converges
	uniformly to some \(\psi\in\mathcal{C}(\dtorus)\). Since
	\(\mathcal{W}\) is closed, \(\psi_n\in\mathcal{W}(r_n,\nu_n)\),
	\(r_n\rightarrow T\), \(\nu_n\rightarrow\nu\) and
	\(\psi_n\rightarrow\psi\) give
	\(\psi\in\mathcal{W}(T,\nu)=\{G(\cdot,\nu)\}\), i.e.
	\begin{equation}\label{eq:psiisG}
	\psi=G(\cdot,\nu)=G\big(\cdot,\mathfrak{f}_T[\gamma]\big).
	\end{equation}
	By Lemma \ref{lemma:stability_RAdm}\ref{item:stability:flows} and
	\ref{sol_mfg:cond:F}, the running data
	\(f_n(t,x)\triangleq F(x,\mathfrak{f}_{\pi_n(t)}[\gamma_n])\)
	(\(\pi_n\) as in that lemma, with $s_n=s$ for all $n$,
    see the comment after the lemma)
	converge uniformly to
	\(f_\gamma(t,x)\triangleq F(x,\mathfrak{f}_t[\gamma])\).

	Now we show that the limit \(\gamma\) is a relaxed
    equilibrium.  By Lemma
	\ref{lemma:characterizationofrelaxedequilibria}, for each \(n\),
	\[
	I^r[s,r_n,\gamma_n,f_n,\psi_n]
	=\operatorname{Val}[s,r_n,\mu,f_n,\psi_n]
	=\action{\mu}{\operatorname{Val}[s,r_n,\cdot,f_n,\psi_n]},
	\]
	the second equality by Lemma
	\ref{mfg_dyn:lm:A_properties}\ref{item:A:decoupling}. As
	\(n\rightarrow\infty\), the left-hand side converges to
	\(I^r[s,T,\gamma,f_\gamma,G(\cdot,\nu)]\) by Lemma
	\ref{lemma:stability_RAdm}\ref{item:stability:payoffs}, while the
	functions \(\operatorname{Val}[s,r_n,\cdot,f_n,\psi_n]\) converge
	uniformly to
	\(\operatorname{Val}[s,T,\cdot,f_\gamma,G(\cdot,\nu)]\) by Lemma
	\ref{mfg_dyn:lm:A_properties}\ref{item:A:continuity} (the
	\(\psi_n\) being uniformly Lipschitz and converging uniformly to
	\(G(\cdot,\nu)\), and the \(f_n\) as required); hence the
	right-hand side converges to
	\(\action{\mu}{\operatorname{Val}[s,T,\cdot,f_\gamma,
	G(\cdot,\nu)]}
	=\operatorname{Val}[s,T,\mu,f_\gamma,G(\cdot,\nu)]\)
	(Lemma \ref{mfg_dyn:lm:A_properties}\ref{item:A:decoupling}
	again). We conclude
	\[
	I^r\big[s,T,\gamma,f_\gamma,G(\cdot,\mathfrak{f}_T[\gamma])\big]
	=\operatorname{Val}\big[s,T,\mu,f_\gamma,
	G(\cdot,\mathfrak{f}_T[\gamma])\big],
	\]
	and Lemma \ref{lemma:characterizationofrelaxedequilibria} (in its
	\(r=T\), \(\sigma=G(\cdot,\mathfrak{f}_T[\gamma])\) form) shows
	that \(\gamma\) is a relaxed equilibrium with
	\(\mathfrak{f}_s[\gamma]=\mu\).

	Last we address the value condition. Passing to the limit in
	\eqref{eq:displayphiofyn} ---that is, in
	\(\phi=A^{s,r_n}_{\mathfrak{f}_\cdot[\gamma_n]}\psi_n\)--- by
	Lemma \ref{mfg_dyn:lm:A_properties}\ref{item:A:continuity}, and
	using \eqref{eq:psiisG},
	\[
	\phi
	=A^{s,T}_{\mathfrak{f}_\cdot[\gamma]}
	G\big(\cdot,\mathfrak{f}_T[\gamma]\big).
	\]
	Conditions 2(a) and 2(b) of Corollary
	\ref{sol_mfg:corollary:value_mfg} are thus verified at
	\((s,\mu,\phi)\). As \((s,\mu,\phi)\) was arbitrary,
	\(\mathcal{W}\) is a value multifunction. The proof is complete.
	\end{proof}

	\begin{remark}\label{rem:closednessnatural}
Let $\mathcal V^{\max}(T,\mu)\triangleq\{G(\cdot,\mu)\}$ and, for $s<T$,
let $\mathcal V^{\max}(s,\mu)$ be the set of the functions $u(s,\cdot)$
over all solutions $(u,m)$ of the MFG on $[s,T]$ with initial
distribution $\mu$. By Definition \ref{sol_mdf:def:value_multifunction},
$\mathcal V^{\max}$ is a valuemultifunction, and every value multifunction $\mathcal W$ satisfies
$\mathcal W(s,\mu)\subset\mathcal V^{\max}(s,\mu)$. By Proposition \ref{sol_mfg:prop:closure}, $\operatorname{cl}\mathcal V^{\max}$ is a
value multifunction, hence $\operatorname{cl}\mathcal V^{\max}\subset
\mathcal V^{\max}$, i.e.~$\mathcal V^{\max}$ is closed. Thus the
closedness required in Theorem \ref{thm:aetoDPP} is a property of the
largest value multifunction, and any value multifunction may be
replaced by its closure without leaving the class of value
multifunctions.
\end{remark}

\begin{corollary}\label{cor:localDPP}
Let $\mathcal W:[0,T]\times\wassonespace\rightsquigarrow\mathcal
C(\dtorus)$ be a closed multifunction with \(\mathcal W(T,\mu) =\{G(\cdot,\mu)\}$ for all $\mu\in\wassonespace$. Suppose that there is
$h>0\) such that the condition of Definition \ref{defn:DPP} holds for all
$s\in[0,T)$, $\mu\in\wassonespace$, $\phi\in\mathcal W(s,\mu)$ and
$r\in(s,T)$ with $r-s<h$. Then $\mathcal W$ satisfies the dynamic programming principle
for the MFG and therefore, by Theorem \ref{thm:aetoDPP}, \(\mathcal{W}\)
is a value multifunction.
\end{corollary}
\begin{proof}
Fix $s\in[0,T)$, $\mu\in\wassonespace$, $\phi\in\mathcal W(s,\mu)$ and
$r\in(s,T)$. Choose an integer $k$ with $(r-s)/k<h$ and put
$r_j\triangleq s+j(r-s)/k$, $j=0,\ldots,k$, and
$(\psi_0,\nu_0)\triangleq(\phi,\mu)$. Suppose that, for some
$j\in\{1,\ldots,k\}$, we have $\psi_{j-1}\in\mathcal
W(r_{j-1},\nu_{j-1})$; this holds for $j=1$. Since $r_{j-1}<T$ and
$r_j-r_{j-1}<h$, the hypothesis provides
$(\psi_j,\nu_j)\in\Phi^{r_j,r_{j-1}}(\psi_{j-1},\nu_{j-1})$ with
$\psi_j\in\mathcal W(r_j,\nu_j)$. Proceeding up to $j=k$ and applying
Proposition \ref{mfg_dyn:prop:superposition} $k-1$ times, we obtain
\[
(\psi_k,\nu_k)\in
\big(\Phi^{r_k,r_{k-1}}\circ\cdots\circ\Phi^{r_1,r_0}\big)(\phi,\mu)
=\Phi^{r,s}(\phi,\mu),
\qquad \psi_k\in\mathcal W(r,\nu_k),
\]
which is the condition of Definition \ref{defn:DPP} for $(s,\mu,\phi,r)$.
\end{proof}

We conclude with a remark that puts the example of
Section \ref{section:nonuniquenessexample} in context.
	\begin{remark}\label{rem:examplevaluemultifunction}
		Let the data be as in
		Lemma \ref{prop:allsolutions} of Section
		\ref{section:nonuniquenessexample},
		with, we recall, the function $\mathcal{G}$ defined as
		$\mathcal G(\mu)=f\big(2\int\cos(2\pi x)\,\mu(dx)\big)$.
Let $\mathcal V(T,\mu)\triangleq\{G(\cdot,\mu)\}\) and, for $t_0<T\), let
$\mathcal V(t_0,\mu)$ be the set of the functions $u(t_0,\cdot)$ over all
solutions $(u,m)$ of the MFG on $[t_0,T]$ with initial distribution
$\mu$. By Definition \ref{sol_mdf:def:value_multifunction}, $\mathcal V$
is a value multifunction, and every value multifunction \(\mathcal W\)
satisfies $\mathcal{W}(t_0,\mu)\subset\mathcal V(t_0,\mu)$. By Lemma
\ref{prop:allsolutions} and \eqref{eq:terminalcoefficient},
\[
\mathcal V(t_0,m_*)=\big\{2\log\big(1+\theta e^{-4\pi^2(T-t_0)}
\cos(2\pi\,\cdot)\big):\
\theta\in[\theta_1,\theta_2],\ \theta=f(B[\theta;m_*])\big\},
\]
where $B[\theta;m_*]$ was introduced in \eqref{eq:terminalcoefficient}.
For instance, for the terminal cost of Corollary \ref{cor:allt0}, Proposition
\ref{prop:nonuniqueness} gives
\[
\mathcal V(t,\mu_\sigma)\supset\Big\{2\log\big(1+\sigma\cos(2\pi\,\cdot)\big),\
2\log\Big(1+\frac{\kappa}{\sigma}\cos(2\pi\,\cdot)\Big)\Big\}
\]
for all $t$ in the interval determined by \eqref{eq:conditiont0sigma},
while, by Proposition \ref{prop:uniquenessnearT}, $\mathcal V(t,\mu)$ is a
singleton for every $\mu$ when $t$ is sufficiently close to $T$.

Let us also describe what the dynamic programming principle (Definition
\ref{defn:DPP}) amounts to in this example. Let $s\in[0,T)$,
$\mu\in\mathcal P^1(\onetorus)$, $\theta\in[\theta_1,\theta_2]$, and let
$\bar m^\theta$ denote the path constructed before Lemma
\ref{prop:allsolutions}, with $(s,\mu)$ in place of $(t_0,m_*)$.
Assume, as we may, that $K_G$ and $M_G$ dominate the Lipschitz constants
and the norms of the functions $2\log(1+\theta\cos(2\pi\,\cdot))$,
$\theta\in[\theta_1,\theta_2]$. Then, for every $r\in(s,T)$,
\begin{equation*}
\Phi^{r,s}(u^\theta(s,\cdot),\mu)
=\big\{\big(u^\theta(r,\cdot),\,\bar m^\theta(r,\cdot)\big)\big\}.
\end{equation*}
Indeed, $u^\theta$ solves the Bellman equation on $[s,r]$, so, by
Corollary \ref{mfg_dyn:cor:verification_identification},
$u^\theta(s,\cdot)=A^{s,r}u^\theta(r,\cdot)$ (since $F=0$, the operators
$A^{s,r}$ do not depend on the flow), and the feedback $Du^\theta$, whose
flow from $\mu$ is $\bar m^\theta$, is a maximizer; by Lemma
\ref{lemma:characterizationofrelaxedequilibria} the corresponding relaxed
control is a relaxed equilibrium for $u^\theta(r,\cdot)$ on $[s,r]$.
Conversely, let $(\sigma,\nu)\in\Phi^{r,s}(u^\theta(s,\cdot),\mu)$. The
change of variables $w=e^{u/2}$ used at the beginning of this subsection
gives $e^{A^{s,r}\sigma/2}=e^{\sigma/2}\star p_{r-s}$ for smooth $\sigma$, and, by Lemma \ref{mfg_dyn:lm:A_properties}\ref{item:A:onelip} and uniform
approximation, for Lipschitz $\sigma$. Hence
$e^{\sigma/2}\star p_{r-s}=q^\theta(s,\cdot)=q^\theta(r,\cdot)\star p_{r-s}$;
since convolution with $p_{r-s}$ multiplies the $k$-th Fourier coefficient by $e^{-4\pi^2k^2(r-s)}\neq0$, $\sigma=u^\theta(r,\cdot)$. The
corresponding relaxed equilibrium is then a maximizer for this datum
(Lemma \ref{lemma:characterizationofrelaxedequilibria}), so, by Lemma
\ref{prel:lm:verification}\ref{lm:verification:c} and the convexity of
$L$, its velocity field is $Du^\theta$ along its flow, and its flow,
being a continuous path solving the FPK
equation with drift $Du^\theta$
and initial distribution $\mu$, coincides
with $\bar m^\theta$ by the uniqueness of such solutions (which do not use the terminal condition). Thus
$\nu=\bar m^\theta(r,\cdot)$.

Consequently, if $\mathcal W$ is a multifunction with
$\mathcal W(T,\cdot)=\{G(\cdot,\cdot)\}$ whose values at times $s<T$ are
among the functions $u^\theta(s,\cdot)$, $\theta\in[\theta_1,\theta_2]$,
then $\mathcal W$ satisfies the dynamic programming principle if and only
if, for all $s<r<T$, $\mu$ and $\theta$,
\[
u^\theta(s,\cdot)\in\mathcal W(s,\mu)\quad\Longrightarrow\quad
u^\theta(r,\cdot)\in\mathcal W\big(r,\bar m^\theta(r,\cdot)\big).
\]
In other words, the dynamic programming principle is the invariance of the
set $\{(s,\mu,\theta):u^\theta(s,\cdot)\in\mathcal W(s,\mu),\ s<T\}$ under
the explicit flow $(s,\mu,\theta)\mapsto(r,\bar m^\theta(r,\cdot),\theta)$,
along which $\theta$ stays constant. This condition does not involve $G$.
The terminal cost enters only through $\mathcal W(T,\cdot)=\{G\}$ and the
closedness of $\mathcal W$: letting $r\to T$ along a trajectory of the
flow, closedness gives $u^\theta(T,\cdot)=G(\cdot,\bar m^\theta(T,\cdot))$,
i.e.~\eqref{eq:generalfixedpoint}, which by
Lemma \ref{prop:allsolutions}
is exactly the statement that the trajectory is a solution.
This illustrates the mechanism of the proof of Theorem \ref{thm:aetoDPP}.
\end{remark}

\begin{appendices}
\section{Properties of relaxed control process}\label{sect_A:relaxed_CS}
\begin{proof}[Proof of Proposition \ref{prel:prop:link_adm}]
	The first part of the proposition directly follows from
	\eqref{prel:intr:adm_embedding}. To prove the converse statement,
    we choose $\gamma\in\RAdm{s}{r}$.
	Further, for each $\phi\in \mathcal{C}_c^{1,2}([s,r]\times \td)$,
	\[
    \begin{split}
		0&=\int_{[s,r]\times \td\times\rd}
		\big[\partial_t\phi(t,x)
		+D\phi(t,x)\cdot v
		+\triangle\phi(t,x)\big]\gamma(dt,dx,dv)\\&=
		\int_s^r\int_{\td}
		\int_{\rd}\big[\partial_t\phi(t,x)
		+D\phi(t,x)\cdot v
		+\triangle\phi(t,x)\big]
		\gamma(dv|t,x)\mathfrak{f}_t[\gamma](dx)dt
		\\
		&=
		\int_s^r\int_{\td}\int_{\rd}
		\big[\partial_t\phi(t,x)+D\phi(t,x)\cdot
		\mathfrak{v}[\gamma](t,x)
		+\triangle\phi(t,x)\big]\mathfrak{f}_t[\gamma](dx)dt.
	\end{split}
	\]
    Thus,
    $(\mathfrak{f}_\cdot[\gamma],\mathfrak{v}[\gamma])
	\in\RAdm{s}{r}$.
\end{proof}

\begin{proof}[Proof of Proposition \ref{prel:prop:Holder_m}]
We mention without proof (it is done through a standard
density argument) that an equivalent definition of
the weak formulation
\eqref{eq:weakcontrolFPK} is
\begin{align*}
    \int\limits_{[s,r]\times \dtorus\times \mathbb{B}_C}
    \psi(x) \eta'(t) \gamma(dtdxdv) =
    -\int_s^r \int\limits_{\dtorus\times\mathbb{B}_C}
    \eta(t) (v\cdot D\psi(x) + \Delta \psi(x) )
    \gamma(dxdv|t) dt,
\end{align*}
for every $\psi\in \mathcal{C}^{\infty}(\dtorus)$,
where $\{\gamma(\,\cdot\,|t)\}_{s\leq t\leq r}$ is the disintegration
of $\gamma$ w.r.t.~Lebesgue measure on $[s,r]$
(equivalently said, the distributional derivative
$\frac{d}{dt} \int_{\dtorus} \psi(x) \gamma(dtdxdv)$ must
be equal to $v\cdot D\psi(x) + \Delta \psi(x)$).
Fix $\psi\in \mathcal{C}^{\infty}(\dtorus)$,
let $h_{\psi}(t) \triangleq \int_{\dtorus} \psi(x) \mathfrak{f}_t[\gamma](dx)$,
which is defined a.e. Applying the weak formulation
equality to $\phi(t,x) \triangleq \eta(t)\psi(x)$,
we obtain
\begin{align*}
    \int_s^r \eta'(t)h_{\psi}(t) dt =
        -\int_s^r \eta(t) g_{\psi}(t) dt,
\end{align*}
where
\begin{equation*}
    g_{\psi}(t) \triangleq \int (D\psi(x)\cdot v + \Delta \psi(x))\gamma(dxdv|t).
\end{equation*}
Since $|v|\leq C$ on the support of $\gamma$, we have
$\|g_{\psi}\|_{\infty}\leq C\|D\psi\|_{\infty} + \|\Delta \psi\|_{\infty}$.
This means that $h_{\psi}$ has a distributional derivative in $L^{\infty}(s,r)$,
so it has a Lipschitz representative (equal to it a.e.), which we
denote by $\tilde{h}_{\psi}$, on the
closed interval $[s,r]$. Let \(\{\psi_k\}_{1}^{\infty}\) a family of smooth
functions dense in \(\dtorus\) and a set $N\subset [s,r]$ of zero Lebesgue
measure such that for every $t\in [s,r]\setminus N$, $\gamma(\,\cdot\,|t)$ is a probability
measure and $h_{\psi_k}(t)=\tilde{h}_{\psi_k}(t)$ for all $k$.
Let $t\in [s,r]$ be arbitrary and $\{t_n\}_{1}^{\infty}\subset [s,r]\setminus N$
such that $t_n\to t$; for each $k\geq 1$ we have
$\tilde{h}_{\psi_{k}}(t_n) \to \tilde{h}_{\psi_k}(t)$ as $n\to\infty$.
The functional $\mathcal{K}_t$, mapping to $\psi_k$ to $\tilde{h}_{\psi_k}(t)$,
as the pointwise limit of positive normalized (defined in the same way)
functionals $\mathcal{K}_{t_n}$, is positive, linear, $\mathcal{K}_t(1)=1$, and
\(|\mathcal{K}_t(\psi_k)| \leq \|\psi_k\|_{\infty}\), so
we extend it by continuity to $\mathcal{C}(\dtorus)$.
By Riesz' representation theorem we obtain the stated
$\mu_t\in\mathcal{P}^1(\dtorus)$
for every $t\in [s,r]$, since $\action{\mu_t}{\psi_k} =
\tilde{h}_{\psi_k}(t)$ is continuous in $t$ for
all $k$ (therefore, as can be easily verified,
for all \(\psi\in\mathcal{C}(\dtorus)\)), which
gives the narrow continuity of $t\mapsto \mu_t$
on $[s,r]$.

Now we re-denote $\mu_t$ back to $\mathfrak{f}_t[\gamma]$
and invoke the superposition principle \cite{Figalli2008}:
	there exists a stochastic process
	$X$ defined on $[s,r]$ that satisfies
	\[
	dX_t=\mathfrak{v}[\gamma](t,X_t)dt+\sqrt{2}dW_t
	\]
	such that $\operatorname{Law}(X_t)=\mathfrak{f}_t[\gamma]$.
    Then,
	for $t_1,t_2\in [s,r]$, $t_1<t_2$,
    \begin{align*}
        W_1(\mathfrak{f}_{t_1}(\gamma),
		\mathfrak{f}_{t_2}(\gamma))
        \leq & \
        \EE|X_{t_2}-X_{t_1}|
\leq C(t_2-t_1)+\sqrt{2}\,\EE|W_{t_2}-W_{t_1}|
\\
\leq & \ C(t_2-t_1)+\sqrt{2d(t_2-t_1)}
\leq  (t_2-t_1)^{1/2}\big(C\sqrt{T}+\sqrt{2d}\big),
    \end{align*}
and the proposition is proved.
\end{proof}
\begin{proof}[Proof of Propsition \ref{prel:prop:continuity}]
	Assume that a sequence $\{\gamma^n\}_{n=1}^\infty\subset
	\RAdm[C]{s}{r}$
	converges narrowly to a relaxed
	control process $\gamma$.
	We are to show that
	a sequence of probability paths
	$\{\mathfrak{f}_\cdot(\gamma^n)\}_{n=1}^{\infty}$
	converges to
	$\mathfrak{f}_\cdot[\gamma]$ in $\mathcal{C}^0([s,r];\mathcal{P}^1(\td))$.
	Denote \[
	m^n_t\triangleq
	\mathfrak{f}_t[\gamma^n],\ \
	m_t\triangleq \mathfrak{f}_t[\gamma].
	\]
	By Proposition \ref{prel:prop:Holder_m}, there exists a constant
	$C'$ such that, for each $n$ and $t_1,t_2\in [s,r]$,
	\[
	W_1(m^n_{t_1},m^n_{t_2})\leq C'|t_1-t_2|^{1/2},
	\ \ W_1(m(t_1),m(t_2))\leq C'|t_1-t_2|^{1/2}.
	\] Thus, for each $1$-Lipschitz continuous function
	$\psi$ from $\td$ to $\mathbb{R}$, one has
	\begin{equation}\label{prel:ineq:psi_m_n}\Bigg|
		\int_{\td}\psi(x)m^n_{t}(dx)-h^{-1}\int_t^{t+h}
		\int_{\td}\psi(x)m^n_{t'}(dx)dt'\Bigg|\leq C'h^{1/2},
	\end{equation} and
	\begin{equation}\label{prel:ineq:psi_m_limit}
		\Bigg|
		\int_{\td}\psi(x)m_{t}(dx)-h^{-1}
		\int_t^{t+h}\int_{\td}\psi(x)m_{t'}(dx)dt'
		\Bigg|\leq C'h^{1/2}.
	\end{equation}
	Further,
	\[
	\begin{split}
		\int_t^{t+h}\int_{\td}\psi&(x)m^n_{t'}(dx)dt'
		=\int_{[t,t+h]\times \td\times\rd}\psi(x)
		\gamma^n(dt',dx,dv)
		\\
		&\longrightarrow
		\int_{[t,t+h]\times \td\times\rd}\psi(x)
		\gamma(dt',dx,dv)
		=\int_t^{t+h}\int_{\td}\psi(x)
		m_{t'}(dx)dt'
	\end{split}
	\]
	as \(n\to\infty\).
	This, together
	with \eqref{prel:ineq:psi_m_n}
	and \eqref{prel:ineq:psi_m_limit},
	implies that, for each $1$-Lipschitz
	continuous function $\psi$,
	\begin{equation}
		\int_{\td}\psi(x)m^n_{t}(dx)
		\rightarrow
		\int_{\td}\psi(x)m_{t}(dx)\text{ as }n\rightarrow\infty.
	\end{equation}
	Therefore, using the tightness of
	Lipschitz functions,
	we conclude that, for each $t$, the sequence
	$\{m^n_t\}_{n=1}^{\infty}$
	converges to $m_t$ narrowly.
	Due to \cite[Theorem 5.9]{santambrogiobook},
	we have that $W(m^n_t,m_t)\rightarrow 0$
	as $n\rightarrow\infty$ for every $s\leq t\leq r$.
	The uniform convergence of $m^n_\cdot$ to $m_\cdot$
	follows from this and Proposition \ref{prel:prop:Holder_m}.
\end{proof}
\end{appendices}

\printbibliography

\end{document}